\documentclass[preprint,3p,12pt]{elsarticle}
\usepackage{mathrsfs}
\usepackage{amsmath}
\usepackage{stmaryrd}
\usepackage{bbding}
\usepackage{dcolumn}
\usepackage{graphicx}
\usepackage{amsfonts}
\usepackage{amssymb}
\usepackage{psfrag}
\usepackage{wrapfig}
\usepackage{subfigure}
\usepackage{makeidx}
\usepackage{bm}
\usepackage{epsf}
\usepackage{epsfig}
\usepackage{setspace}
\usepackage{graphicx}
\usepackage{epstopdf}
\usepackage{psfrag}
\usepackage{subfigure}

\begin{document}

\title{An Efficient and Accurate Two-Stage Fourth-order Gas-kinetic Scheme for the Navier-Stokes Equations}

\author[HKUST1]{Liang Pan}
\ead{panliangjlu@sina.com}

\author[HKUST1,HKUST2]{Kun Xu\corref{cor1}}
\ead{makxu@ust.hk}

\author[tsinghua]{Qibing Li}
\ead{lbq@tsinghua.edu.cn}

\author[iapcm]{Jiequan Li}
\ead{li\_jiequan@iapcm.ac.cn}

\address[HKUST1]{Department of Mathematics, Hong Kong University of Science and Technology, Clear Water Bay, Kowloon, Hong Kong}
\address[HKUST2]{Department of Mechanical and Aerospace Engineering, Hong Kong University of Science and Technology, Clear Water Bay, Kowloon, Hong Kong}
\address[tsinghua]{School of Aerospace, Tsinghua University, Beijing 100084, China}
\address[iapcm]{Institute of Applied Physics and Computational Mathematics, Beijing, 100088, China}
\cortext[cor1]{Corresponding author}

\begin{abstract}
For computational fluid dynamics (CFD), the generalized Riemann problem (GRP) solver and the gas-kinetic kinetic scheme (GKS)
provide a time-accurate flux function starting from a discontinuous piecewise linear flow distributions around each cell interface.
With the use of time derivative of the flux function, a two-stage Lax-Wendroff-type  (L-W for short) time stepping method has been recently proposed
in the design of a fourth-order time accurate method \cite{GRP-high}.
In this paper, based  on the same time-stepping method and the second-order GKS flux function \cite{GKS-Xu2},
a fourth-order gas-kinetic scheme is constructed for the Euler and Navier-Stokes equations.
In comparison with the formal one-stage time-stepping third-order gas-kinetic solver \cite{GKS-high1},
the current fourth-order method not only reduces the complexity of the flux function, but also improves the accuracy of
the scheme, even though the third- and fourth-order schemes have similar computation cost.
Most importantly, the robustness of the fourth-order GKS is as good as the second-order one. Perfect numerical solutions can be obtained from the high Reynolds
number boundary layer solutions to the hypersonic viscous heat conducting flow computations.
Many numerical tests, including many difficult ones for the Navier-Stokes solvers,  have been used to validate the  current fourth-order method.
Following the two-stage time-stepping framework, the one-stage third-order GKS can be easily extended to a fifth-order method with the usage of
both first-order and second-order time derivatives of the  flux function. The use of time-accurate flux function may have great impact on the
development of higher-order CFD methods.

\end{abstract}
\maketitle

\section{Introduction}

To develop third and higher-order numerical methods has attracted great attention in recent decades.
In comparison with  second-order schemes, which were mostly developed in the 70s and 80s,
the higher-order methods can provide more accurate solutions, but they are less robust and more complicated.
There are many review papers and monographs about the current status of higher-orders schemes,
which include the discontinuous Galerkin (DG) \cite{DG1}, essential
non-oscillatory (ENO) \cite{ENO1}, weighted essential
non-oscillatory (WENO) \cite{WENO2},$P_NP_M$
\cite{PNPM1}, multi-moment constrained method \cite{Ii}, and many others. Most of those methods use the
Runge-Kutta  time-stepping approach to achieve higher order temporal
accuracy \cite{TVD-rk}. Based on the time-independent flux function of the Riemann solver \cite{Riemann-appro},
in order to achieve a fourth-order time accuracy, four-stage Runge-Kutta time stepping method has to be used. Moreover,
the CFL number for those methods strongly depend on the order of the scheme, such as the DG method.

Recently, based on the time-dependent flux function of the generalized Riemann problem (GRP) \cite{GRP1,GRP2,GRP3},
a two-stage fourth order time-accurate discretization was
developed for Lax-Wendroff type (L-W for short) flow solvers, particularly applied  for the hyperbolic conservation laws \cite{GRP-high}. The reason for the success of a two-stage L-W type time stepping
method in achieving a fourth-order time accuracy is solely due to the use of both flux function and its time derivative.
In terms of the gas evolution model, the gas-kinetic scheme (GKS) provides a time accurate flux function as well, even though it depends on
time through a much more complicated relaxation process from the kinetic to the hydrodynamic scale physics than the time-dependent flux function of GRP.
 This paper is about to construct a fourth-order time accurate gas-kinetic
scheme (GKS) with the two-stage temporal discretization for the Euler and Navier-Stokes (NS)  equations.

For the NS solutions, second-order and third-order gas-kinetic schemes have been constructed in the past years \cite{GKS-Xu2,GKS-high1,GKS-Kumar,GKS-DG}.
The flux evaluation in the scheme is based on the time evolution of flow
variables from an initial piece-wise discontinuous polynomials
around each cell interface, where high-order spatial and temporal
evolutions of a gas distribution function are coupled nonlinearly.
In comparison with other high-order schemes, the GKS integrates the flux function over a time step analytically without
employing the multi-stage Runge-Kutta time stepping techniques.
However, with the one-stage gas evolution model, the formulation of GKS can become
very complicated for the further improvement of the order of the scheme, such as the fourth-order scheme \cite{liuna},
especially for multidimensional computations.
The two-stage L-W time stepping method in \cite{GRP-high} provides a reliable framework
to develop a fourth-order GKS with a second-order flux function.
In this paper, we are going to present such a fourth-order GKS for the Euler
and Navier-Stokes solutions. The current scheme can use a time step with CFL number on the order of $0.5$.
Most importantly,
the current fourth-order GKS is as robust as the second-order method, which works perfectly from the subsonic to the
hypersonic viscous heat conducting flows.
 Numerical tests show that the current
scheme not only has the expected order of accuracy for the smooth flow,
but also has favorable shock capturing property for the discontinuous solutions.
As a further extension, the third-order flux function \cite{Qian-Li, GKS-high1,GKS-high2} can be also used to construct  two-stage  fifth-order temporal accurate methods
with the inclusion of both first-order and second-order time derivatives of the flux function. The detailed formulation is presented in the Appendix of this paper.
Theoretically, this process for constructing even higher-order schemes can go forward continuously.

This paper is organized as follows. In Section 2, the general
formulation for the two-stage temporal discretization is introduced.
In Section 3, a fourth-order gas-kinetic scheme is presented based
on the two-stage time discretization. Section 4 includes numerical
examples to validate the current algorithm. The last section is drawing  the conclusion.
 The extension for  the construction of two-stage fifth-order schemes  is given in Appendix.

\section{Fourth-order temporal discretization}

A two-stage fourth-order time-accurate discretization was developed
for  Lax-Wendroff flow solvers, particularly applied for  hyperbolic equations with the generalized Riemann problem (GRP) solver \cite{GRP-high}.
Consider  the following time-dependent equation,
\begin{align}\label{pde}
\frac{\partial \textbf{w}}{\partial t}=\mathcal {L}(\textbf{w}),
\end{align}
with the initial condition at $t_n$, i.e.,
\begin{align}\label{pde2}
\textbf{w}(t=t_n)=\textbf{w}^n,
\end{align}
where $\mathcal {L}$ is an operator for spatial derivative of flux.
The time derivatives are obtained using the Cauchy-Kovalevskaya method,
\begin{align*}
\frac{\partial \textbf{w}^n}{\partial t}=\mathcal{L}(\textbf{w}^n),~ \ \ \
\frac{\partial }{\partial t}\mathcal
{L}(\textbf{w}^n)=\frac{\partial }{\partial \textbf{w}}\mathcal
{L}(\textbf{w}^n)\mathcal {L}(\textbf{w}^n).
\end{align*}
Introducing an intermediate state at $t^*=t_n+\Delta t/2$,
\begin{align}\label{step1}
\textbf{w}^*=\textbf{w}^n+\frac{1}{2}\Delta t\mathcal
{L}(\textbf{w}^n)+\frac{1}{8}\Delta t^2\frac{\partial}{\partial
t}\mathcal{L}(\textbf{w}^n),
\end{align}
the corresponding time  derivatives are obtained as well for the intermediate stage state,
\begin{align*}
\frac{\partial \textbf{w}^*}{\partial t}=\mathcal{L}(\textbf{w}^*),~
\frac{\partial }{\partial t}\mathcal
{L}(\textbf{w}^*)=\frac{\partial }{\partial \textbf{w}}\mathcal
{L}(\textbf{w}^*)\cdot \mathcal {L}(\textbf{w}^*).
\end{align*}
Then, the state $\textbf{w}$ can be updated with the following formula,
\begin{align}\label{step2}
\textbf{w}^{n+1}=\textbf{w}^n+\Delta t\mathcal
{L}(\textbf{w}^n)+\frac{1}{6}\Delta t^2\big(\frac{\partial}{\partial
t}\mathcal{L}(\textbf{w}^n)+2\frac{\partial}{\partial
t}\mathcal{L}(\textbf{w}^*)\big).
\end{align}
It can be proved that the above time stepping method with Eq.\eqref{step1} and
Eq.\eqref{step2}  provides a fourth-order time accurate solution for $\textbf{w}(t)$ at $t=t_n +\Delta t$.
The details of the analysis can be
found in \cite{GRP-high}.
Thus, based on a time accurate solution
${\partial\mathcal{L}}/{\partial t}$,
a fourth-order temporal
accuracy can be achieved from the two-stage discretization of Eq.\eqref{pde} through Eq.\eqref{step1} and
Eq.\eqref{step2}.

 We apply this approach for  conservation laws
\begin{align}\label{euler}
\frac{\partial \textbf{w}}{\partial t}+\frac{f(\textbf{w})}{\partial
x}=0,
\end{align}
where $\textbf{w}$ is a conservative variable and
$f(\textbf{w})$ is the corresponding flux, which includes all terms related to the viscous heat conducting flow.
The semi-discrete form
of a finite volume scheme for equations Eq.\eqref{euler}
can be written as
\begin{align}\label{conser}
\frac{\partial \textbf{w}_i}{\partial t}=\mathcal
{L}_i(\textbf{w})=-\frac{1}{\Delta
x_i}(\textbf{f}_{i+1/2}-\textbf{f}_{i-1/2}),
\end{align}
where $\textbf{w}_i$ are the cell averaged conservative variables of
the cell $I_i=[x_{i-1/2}, x_{i+1/2}]$,  $\textbf{f}_{i+1/2}$ are the
fluxes at the cell interface $x=x_{i+1/2}$, and $\Delta
x_i=x_{i+1/2}-x_{i-1/2}$.
A similar finite volume formulation can be obtained in two- and three-dimensional cases. Then \eqref{conser} falls into the framework of the two-stage L-W time stepping.

\section{A fourth-order gas-kinetic scheme}
The similarity between the generalized Riemann problem (GRP) solver and the
gas-kinetic scheme has been studied in \cite{GRP-GKS}. In both
schemes, the spatial and temporal accuracy are coupled through a generalized
Lax-Wendroff-type procedure for the discontinuous cases, and a single stage time integration is used for the
flux transport across a cell interface for the second-order schemes.
In this section, a fourth-order gas-kinetic scheme from a second-order flux function will be constructed
through a two-stage time discretization framework of Eq.\eqref{step1}
and Eq.\eqref{step2} for the Euler and Navier-Stokes solutions.

\subsection{Second-order gas-kinetic flux solver}
The two-dimensional BGK equation can be written as \cite{BGK-1},
\begin{equation}\label{bgk}
f_t+\textbf{u}\cdot\nabla f=\frac{g-f}{\tau},
\end{equation}
where $f$ is the gas distribution function, $g$ is the corresponding
equilibrium state, and $\tau$ is the collision time. The collision
term satisfies the compatibility condition
\begin{equation}\label{compatibility}
\int \frac{g-f}{\tau}\psi d\Xi=0,
\end{equation}
where $\psi=(1,u,v,\displaystyle \frac{1}{2}(u^2+v^2+\xi^2))$,
$d\Xi=dudvd\xi^1...d\xi^{K}$, $K$ is the number of internal freedom,
i.e.  $K=(4-2\gamma)/(\gamma-1)$ for two-dimensional flows, and
$\gamma$ is the specific heat ratio. The conservative variables are
denoted as $W=(\rho, \rho U, \rho V, \rho E)$.  In the smooth
region, the gas distribution function can be expanded as
\begin{align*}
f=g-\tau D_{\textbf{u}}g+\tau D_{\textbf{u}}(\tau
D_{\textbf{u}})g-\tau D_{\textbf{u}}[\tau D_{\textbf{u}}(\tau
D_{\textbf{u}})g]+...,
\end{align*}
where $D_{\textbf{u}}={\partial}/{\partial
t}+\textbf{u}\cdot \nabla$. By truncating on different orders of
$\tau$, the corresponding macroscopic equations can be derived. For
the Euler equations, the zeroth order truncation is taken, i.e.
$f=g$. For the Navier-Stokes equations, the first order truncation
is used,
\begin{align}\label{ns}
f=g-\tau (ug_x+vg_y+g_t).
\end{align}
Based on the higher order truncations, the Burnett and super-Burnett
equations can be obtained \cite{BGK-3,GKS-B,GKS-SB}.

In order to update the flow variables, the flux is based on the integral solution of gas distribution function from the BGK
equation at a cell interface,
\begin{equation}\label{integral1}
f(x_{i+1/2},t,u,v,\xi)=\frac{1}{\tau}\int_0^t g(x',y',t',u,v,\xi)e^{-(t-t')/\tau}dt'\\
+e^{-t/\tau}f_0(-ut,y-vt,u,v,\xi),
\end{equation}
where $x_{i+1/2}=0$ is the location of the cell interface,
$x=x'+u(t-t')$ and $y=y'+v(t-t')$ are the trajectory of particles,
$f_0$ is the initial gas distribution function, and $g$ is the
corresponding equilibrium state. According to Eq.\eqref{integral1},
the time dependent gas distribution function $f(x_{i+1/2},t,u,\xi)$
at the cell interface $x_{i+1/2}$ can be expressed as \cite{GKS-Xu2,GKS-Xu1}
\begin{align}\label{flux}
f(x_{i+1/2},t,u,v,\xi)=&(1-e^{-t/\tau})g_0+((t+\tau)e^{-t/\tau}-\tau)(\overline{a}_1u+\overline{a}_2v)g_0\nonumber\\
+&(t-\tau+\tau e^{-t/\tau}){\bar{A}} g_0\nonumber\\
+&e^{-t/\tau}g_r[1-(\tau+t)(a_{1r}u+a_{2r}v)-\tau A_r)]H(u)\nonumber\\
+&e^{-t/\tau}g_l[1-(\tau+t)(a_{1l}u+a_{2l}v)-\tau A_l)](1-H(u)).
\end{align}

Based on the spatial reconstruction of  macroscopic flow variables, which will be given in the next
subsection, the conservative variables $W_l$ and $ W_r$ on the left and right hand sides of a cell interface, and the
corresponding equilibrium states $g_l$ and $g_r$, can be determined. Their spatial derivatives in both normal and tangential directions,
such as $(a_{1l},a_{1r},a_{2l},a_{2r})$,
are related to the normal and tangential derivatives of the initial macroscopic flow variables.
The time derivatives $(A_l, A_r)$ can be obtained from the requirement on the first-order Chapman-Enskog expansion,
such as
$$ \int g_l (a_{1l}u+a_{2l}v+A_l) \psi d\Xi =0 ,$$ and
$$ \int g_r (a_{1r}u+a_{2r}v+A_r) \psi d\Xi = 0.$$

Through the compatibility condition Eq.\eqref{compatibility}, the
conservative variables $W_{0}$ and the equilibrium state $g_{0}$ at
the cell interface can be determined as follows,
\begin{align}\label{compatibility2}
\int\psi g_{0}d\Xi=W_0=\int_{u>0}\psi g_{l}d\Xi+\int_{u<0}\psi
g_{r}d\Xi.
\end{align}
Then, with the spatial derivatives of macroscopic flow variables across and along a cell interface and the compatibility condition,
 the coefficients related to the spatial derivatives in the equilibrium state in Eq.\eqref{flux}, such as $(\bar{a_1},\bar{a_2})$, and
  its time derivative $\bar{A}$, can be fully obtained by,
\begin{align}\label{kinetic-lw}
\displaystyle\langle{\overline{a}_1} \rangle =\frac{\partial {
\overline{W} }}{\partial x}, \langle \overline{a}_2 \rangle
=\frac{\partial {\overline{W}}}{\partial y},  \langle \overline{a}_1
u+ \overline{a}_2 v + \overline{ A} \rangle=0,
\end{align}
where $\langle...\rangle$ are the moments of the equilibrium gas
distribution function $g_0$, and defined by
\begin{align*}
\langle...\rangle=\int g_0 (...)\psi d\Xi.
\end{align*}
More details of the gas-kinetic scheme can be found in \cite{GKS-Xu1}.

\subsection{Spatial reconstruction}
The above time evolution solution is based on the high-order initial
reconstruction for macroscopic flow variables, and  the fifth-order
WENO reconstruction  is adopted in this study \cite{WENO}.

For one dimensional computation, $W_l, W_r$ and $W_0$ corresponding
to the equilibrium states $g_l, g_r$ and $g_0$ in Eq.\eqref{flux}
can be constructed at the cell interface $x_{i+1/2}$. The spatial
derivatives ${\partial W}/{\partial x}$ are also
given based on the reconstruction. Especially, for the determination of the
equilibrium state $g_0$ across the cell interface with a fifth-order
of accuracy, the conservative variables around the cell interface
can be expanded as
\begin{align*}
\overline{W}(x)=W_{0}+S_1(x-x_*)+\frac{1}{2}S_2(x-x_*)^2+\frac{1}{6}S_3(x-x_*)^3+\frac{1}{24}S_4(x-x_*)^4.
\end{align*}
With the following conditions,
\begin{align*}
\int_{I_{i+k}} \overline{W}(x)=W_{i+k}, k=-1,...,2,
\end{align*}
the derivatives are given by
\begin{align*}
\overline{W}_x=S_1=\big[-\frac{1}{12}(W_{i+2}-W_{i-1})+\frac{5}{4}(W_{i+1}-W_{i})\big]/\Delta
x.
\end{align*}
For two dimensional computation, the fifth-order Gauss quadrature is
used to achieve the accuracy in space
\begin{align}\label{gauss}
\frac{1}{\Delta
y}\int_{y_{j-1/2}}^{y_{j+1/2}}F(W(x_{i+1/2},y,t))dy=\sum_{\ell=1}^k\omega_\ell F(W(x_{i+1/2},y_\ell,t)),
\end{align}
where $y_l\in [y_{j-1/2}, y_{j+1/2}], \ell= 1,...,3$ are the Gauss
quadrature points, and $\omega_\ell$ are corresponding weights. Based on
the tangential reconstruction, the tangential derivatives at each
Gauss quadrature points can be obtained.

\subsection{Two-stage gas-kinetic scheme}
In this section, a two-stage fourth-order gas-kinetic scheme will be
presented based on the time-dependent gas distribution function \eqref{flux} at each cell interface.

For the gas-kinetic scheme, the gas evolution is a relaxation process from kinetic to hydrodynamic scale through the exponential function,
and the corresponding flux is a complicated function of time.
In order to obtain the time derivatives of the
flux function at $t_n$ and $t_*=t_n + \Delta t/2$ with the correct physics,
the flux function should be approximated as a linear function of time within a time interval.
Let's first introduce the following notation,
\begin{align*}
\mathbb{F}_{i+1/2}(W^n,\delta)
=\int_{t_n}^{t_n+\delta}F_{i+1/2}(W^n,t)dt&=\int_{t_n}^{t_n+\delta}\int
uf(x_{i+1/2},t,u, v,\xi)dud\xi dt.
\end{align*}
In the time interval $[t_n, t_n+\Delta t]$, the flux is expanded as
the following linear form
\begin{align}\label{expansion}
F_{i+1/2}(W^n,t)=F_{i+1/2}^n+ \partial_t F_{j+1/2}^n(t-t_n).
\end{align}
The coefficients $F_{j+1/2}^n$ and $\partial_tF_{j+1/2}^n$ can be
determined as follows,
\begin{align*}
F_{i+1/2}(W^n,t_n)\Delta t&+\frac{1}{2}\partial_t
F_{i+1/2}(W^n,t_n)\Delta t^2 =\mathbb{F}_{i+1/2}(W^n,\Delta t) , \\
\frac{1}{2}F_{i+1/2}(W^n,t_n)\Delta t&+\frac{1}{8}\partial_t
F_{i+1/2}(W^n,t_n)\Delta t^2 =\mathbb{F}_{i+1/2}(W^n,\Delta t/2).
\end{align*}
By solving the linear system, we have
\begin{align}\label{second}
F_{i+1/2}(W^n,t_n)&=(4\mathbb{F}_{i+1/2}(W^n,\Delta t/2)-\mathbb{F}_{i+1/2}(W^n,\Delta t))/\Delta t,\nonumber\\
\partial_t F_{i+1/2}(W^n,t_n)&=4(\mathbb{F}_{i+1/2}(W^n,\Delta t)-2\mathbb{F}_{i+1/2}(W^n,\Delta t/2))/\Delta
t^2.
\end{align}
Similarly, $\displaystyle F_{i+1/2}(W^*,t_*), \partial_t
F_{i+1/2}(W^*,t_*)$ for the intermediate state can be constructed.  For the two-dimensional
computation, the corresponding fluxes in the $y$-direction can be
obtained as well.

With these notations, the two-stage algorithm for
both Euler and Navier-Stokes equations is given as follows
\begin{enumerate}
\item[(i)] With the initial reconstruction, update $W_{ij}^*$ at $\displaystyle t_*=t_n+ \Delta
t/2$ by
\begin{align*}
W_{ij}^*=W_{ij}^n \displaystyle&-\frac{1}{\Delta x}
\big[\mathbb{F}_{i+1/2,j}(W^n,\Delta
t/2)-\mathbb{F}_{i-1/2,j}(W^n,\Delta t/2)\big]\\&-\frac{1}{\Delta y}
\big[\mathbb{G}_{i,j+1/2}(W^n,\Delta
t/2)-\mathbb{G}_{i,j-1/2}(W^n,\Delta t/2)\big],
\end{align*}
and compute the fluxes and their derivatives by Eq.\eqref{second} for future use,
\begin{align*}
F_{i+1/2,j}(W^n,t_n),~G_{i,j+1/2,}(W^n,t_n),~\partial_t
F_{i+1/2,j}(W^n,t_n),~\partial_t G_{i,j+1/2}(W^n,t_n).
\end{align*}
\item[(ii)]  Reconstruct intermediate value $W^*_{ij}$, and compute
\begin{align*}
\partial_t F_{i+1/2,j}(W^*,t_*),~~\partial_t G_{i,j+1/2}
(W^*,t_*),
\end{align*}
where the derivatives are determined by Eq.\eqref{second} in the
time interval $[t_*, t_*+\Delta t]$.
\item[(iii)]  Update $W_{ij}^{n+1}$ by
\begin{align*}
\displaystyle W_{ij}^{n+1}=W_{ij}^n-\frac{\Delta t}{\Delta
x}[\mathscr{F}_{i+1/2,j}^n-\mathscr{F}_{i-1/2,j}^n]-\frac{\Delta
t}{\Delta y}[\mathscr{G}_{i,j+1/2}^n-\mathscr{G}_{i,j-1/2}^n],
\end{align*}
where $\mathscr{F}_{i+1/2,j}^n$ and $\mathscr{G}_{i,j+1/2}^n$ are
the numerical fluxes and expressed as
\begin{align*}
\mathscr{F}_{i+1/2,j}^n&=F_{i+1/2,j}(W^n,t_n)+\displaystyle\frac{\Delta
t}{6}\big[\partial_t F_{i+1/2,j}(W^n,t_n)+2\partial_t F_{i+1/2,j}
(W^*,t_*)\big],\\
\mathscr{G}_{i,j+1/2}^n&=G_{i,j+1/2,}(W^n,t_n)+\displaystyle\frac{\Delta
t}{6}\big[\partial_t G_{i,j+1/2}(W^n,t_n)+2\partial_t G_{i,j+1/2}
(W^*,t_*)\big].
\end{align*}
For each flux, the Gaussian quadratures Eq.\eqref{gauss} are used.
\end{enumerate}

\section{Numerical tests}
In this section, numerical tests for both inviscid and viscous
flows will be presented to validate our numerical scheme. For the
inviscid flow, the collision time $\tau$ takes
\begin{align*}
\tau=\epsilon \Delta t+C\displaystyle|\frac{p_l-p_r}{p_l+p_r}|\Delta
t,
\end{align*}
where $\varepsilon=0.05$ and $C=1$. For the viscous flow, we have
\begin{align*}
\tau=\frac{\mu}{p}+C \displaystyle|\frac{p_l-p_r}{p_l+p_r}|\Delta t,
\end{align*}
where $p_l$ and $p_r$ denote the pressure on the left and right
sides of the cell interface, $\mu$ is the viscous coefficient, and
$p$ is the pressure at the cell interface. In  smooth flow
regions, it will reduce to $\tau=\mu/p$. The ratio of specific heats
takes $\gamma=1.4$.
The reason for including artificial dissipation through the additional term in the particle collision time
is to enlarge the kinetic scale physics in the discontinuous region for the construction of a numerical shock structure through the particle free transport and
inadequate particle collision.

For the smooth flow, the WENO reconstruction can be used directly on the
conservative flow variables. For the flow with strong discontinuity, the characteristic variables can be used in the reconstruction.
Based on $A_{i+1/2,j}=(\partial F/\partial W)_{W=W^*}$, where $W$ are
the conservative variables, $F(W)$ are the corresponding fluxes, and
$W^*=(W_{i,j}+W_{i+1,j})/2$, the cell averaged and point
conservative values can be projected into the characteristic field
by $\omega=RW$, where $R$ is the matrix corresponding to right
eigenvectors of $A$. The reconstruction scheme is applied
on the characteristic variables $\omega$. With the reconstructed
polynomials for characteristic variables, the
conservative flow variables can be recovered by the inverse projection.

\subsection{Accuracy tests} The first case is the advection of
density perturbation, and the initial condition is set as follows
\begin{align*}
\rho(x)=1+0.2\sin(\pi x),\ \  U(x)=1,\ \ \  p(x)=1, x\in[0,2].
\end{align*}
The periodic boundary condition is adopted, and the analytic
solution is
\begin{align*}
\rho(x,t)=1+0.2\sin(\pi(x-t)),\ \ \  U(x,t)=1,\ \ \  p(x,t)=1.
\end{align*}
In the computation, a uniform mesh with $N$ points is used. As analyzed in the
section before, with the fifth-order spatial reconstruction, the
leading truncation error for the fourth-order GKS  is $O(\Delta x^5+\Delta t^4)$.
In these tests, a fixed CFL number $CFL=0.4$ is used for different
meshes. The $L^1$ and $L^2$ errors and orders at $t=2$ are presented
in Table.\ref{tab1}. The fifth order accuracy can be kept until the
mesh $N=640$. As a comparison, with the original second-order GKS,
the leading error is on the order of $O(\Delta x^5+\Delta t^2)$.
With the identical spatial reconstruction and CFL number
$CFL=0.4$, only a second-order accuracy can be achieved. To show the
order of accuracy, a small CFL number $CFL=0.1$ is used. The $L^1$
and $L^2$ errors and orders at $t=2$ are presented in
Table.\ref{tab2}. The fifth order accuracy can be kept at the
beginning. With the mesh refinement, the temporal error becomes the
dominant one and the accuracy reduces to a second order method.

\begin{table}[!h]
\begin{center}
\def\temptablewidth{0.85\textwidth}
{\rule{\temptablewidth}{0.5pt}}
\begin{tabular*}{\temptablewidth}{@{\extracolsep{\fill}}c|cc|cc}
mesh & $L^1$ error & convergence order ~ & $L^2$ error & convergence order  \\
\hline
20  &  4.4759E-004 &  ~        &  3.7653E-004 &  ~        \\
40  &  1.3764E-005 &  5.0231   &  1.1504E-005 &  5.0324  \\
80  &  4.2791E-007 &  5.0075   &  3.4744E-007 &  5.0493  \\
160 &  1.3354E-008 &  5.0018   &  1.0644E-008 &  5.0286  \\
320 &  4.1722E-010 &  5.0003   &  3.2940E-010 &  5.0140  \\
640 &  1.3039E-011 &  4.9998   &  1.0250E-011 &  5.0060   \\
1280&  4.5156E-013 &  4.8517   &  3.5536E-013 &  4.8502
\end{tabular*}
{\rule{\temptablewidth}{0.5pt}}
\end{center}
\vspace{-4mm} \caption{\label{tab1} Accuracy test for the advection
of density perturbation by the fourth-order GKS.}
\end{table}

\begin{table}[!h]
\begin{center}
\def\temptablewidth{0.85\textwidth}
{\rule{\temptablewidth}{0.5pt}}
\begin{tabular*}{\temptablewidth}{@{\extracolsep{\fill}}c|cc|cc}
mesh & $L^1$ error & convergence order ~ & $L^2$ error & convergence order  \\
\hline
20  &4.5797E-004  &   ~       &  3.7856E-004  &   ~            \\
40  &1.3994E-005  &   5.0322  &  1.1735E-005 &5.0116    \\
80  &1.0709E-006  &   3.7078  &  8.5971E-007 &3.7708   \\
160 &2.5659E-007  &   2.0613  &  2.0167E-007 &2.0918   \\
320 &6.4243E-008  &   1.9978  &  5.0455E-008 &1.9989
\end{tabular*}
{\rule{\temptablewidth}{0.5pt}}
\end{center}
\vspace{-4mm} \caption{\label{tab2} Accuracy test for the advection
of density perturbation by the second-order GKS.}
\end{table}

\begin{table}[!h]
\begin{center}
\def\temptablewidth{0.85\textwidth}
{\rule{\temptablewidth}{0.6pt}}
\begin{tabular*}{\temptablewidth}{@{\extracolsep{\fill}}c|cc|cc}
mesh & $L^1$ error & convergence order ~ & $L^\infty$ error & convergence order  \\
\hline
20$\times$20   &1.98E-3       &          &3.79E-2              &         \\
40$\times$40   &1.69E-4       &3.55   &8.08E-3              &2.23  \\
80$\times$80   &8.92E-6       &4.24   &4.10E-4              &4.30  \\
160$\times$160 &2.31E-7       &5.27   &5.29E-6              &6.28  \\
320$\times$320 &7.40E-9       &4.96   &2.09E-7              &4.66  \\
640$\times$640 &2.76E-10     &4.74   &7.09E-9              &4.88
\end{tabular*}
{\rule{\temptablewidth}{0.5pt}}
\end{center}
\vspace{-4mm} \caption{\label{tab3} Accuracy of the fourth-order GKS for the isentropic
vortex propagation at time $t=10$}
\end{table}

The next test is the isotropic vortex propagation problem. The mean
flow is $(\rho,U,V,p)=(1,1,1,1)$, and an isotropic vortex is added
to the mean flow, i.e., with perturbation in $u, v$ and temperature
$T=p/\rho$, and no perturbation in entropy $S=p/\rho^\gamma$. The
perturbation is given by
\begin{align*}
(\delta U,\delta
V)=\frac{\epsilon}{2\pi}e^{\frac{(1-r^2)}{2}}(-y,x), ~ \ \ \delta
T=-\frac{(\gamma-1)\epsilon^2}{8\gamma\pi^2}e^{1-r^2},  \ \ \ \ ~\delta S=0,
\end{align*}
where $r^2=x^2+y^2$ and the vortex strength $\epsilon=5$. The
computational domain is $[-5,5]\times[-5,5]$, the periodic boundary
conditions are imposed on the boundaries in both $x$ and $y$
directions. The exact solution is the perturbation which propagates
with the velocity $(1,1)$.  The $L^1$ and $L^\infty$ errors and
orders at $t=10$ with $N\times N$ uniform mesh cells are presented
in Table \ref{tab3}, which shows that the expected accuracy can be
also achieved for the two dimensional computation.

\begin{figure}[!htb]
\centering
\includegraphics[width=0.444\textwidth]{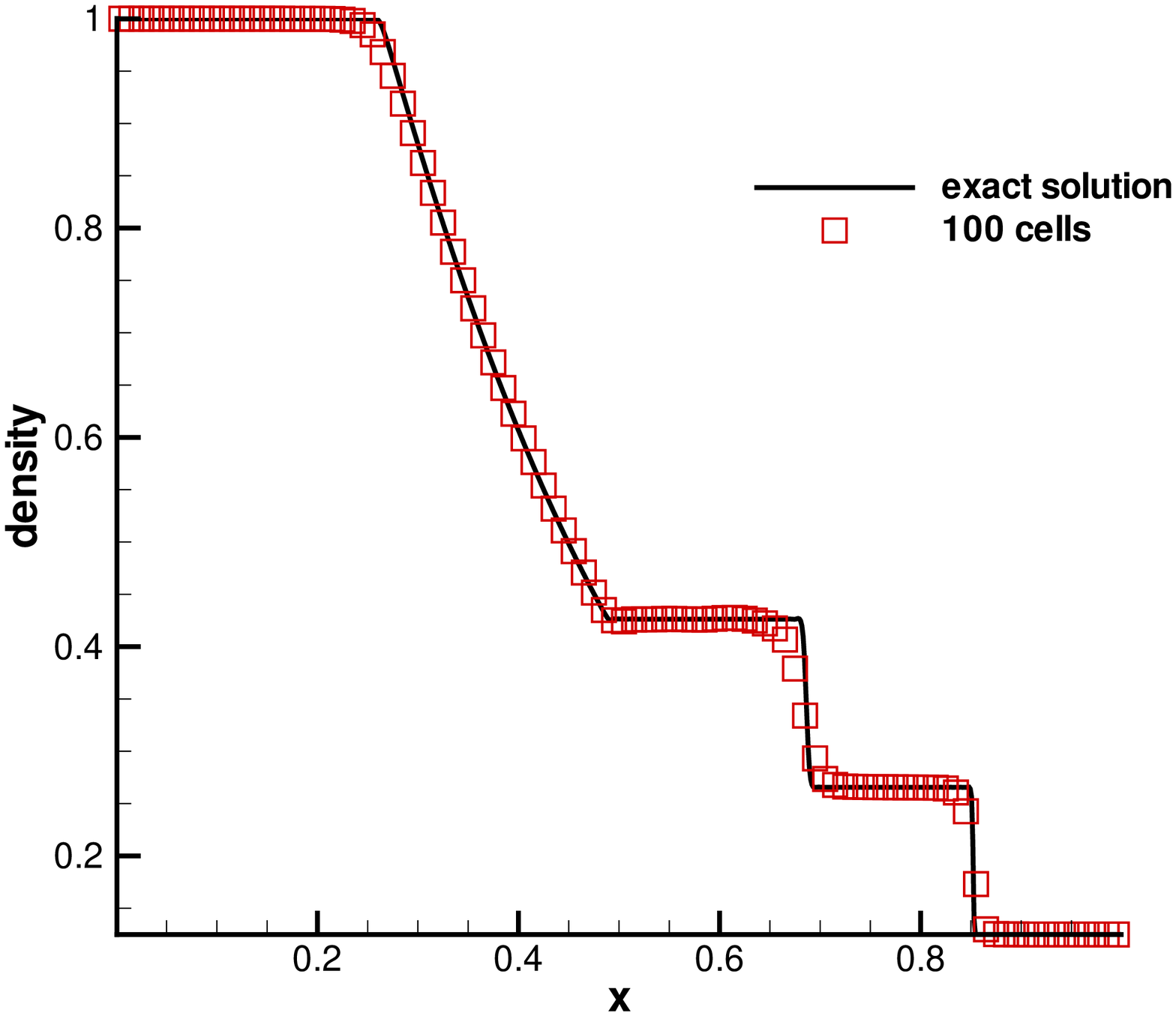}\includegraphics[width=0.444\textwidth]{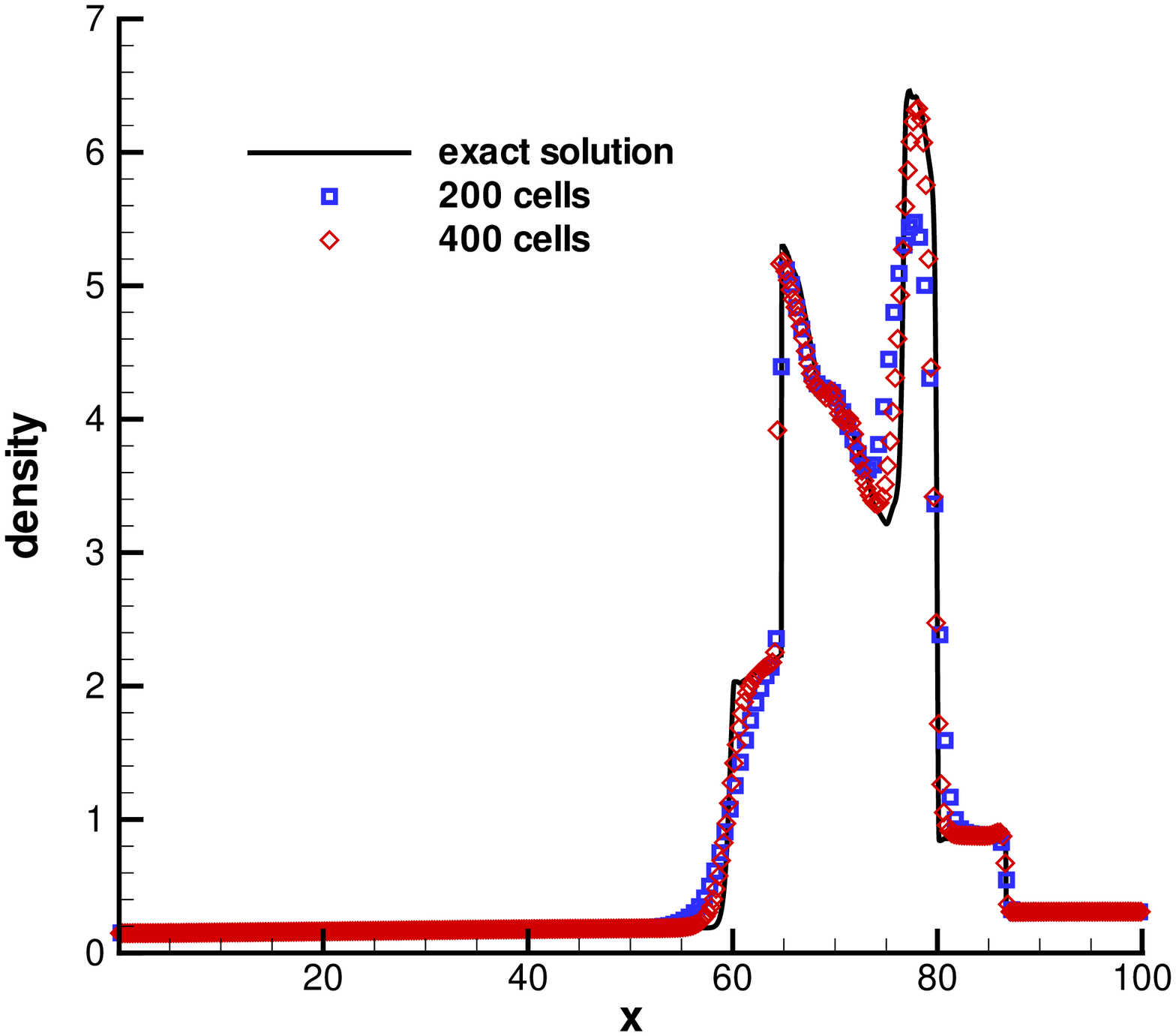}
\includegraphics[width=0.444\textwidth]{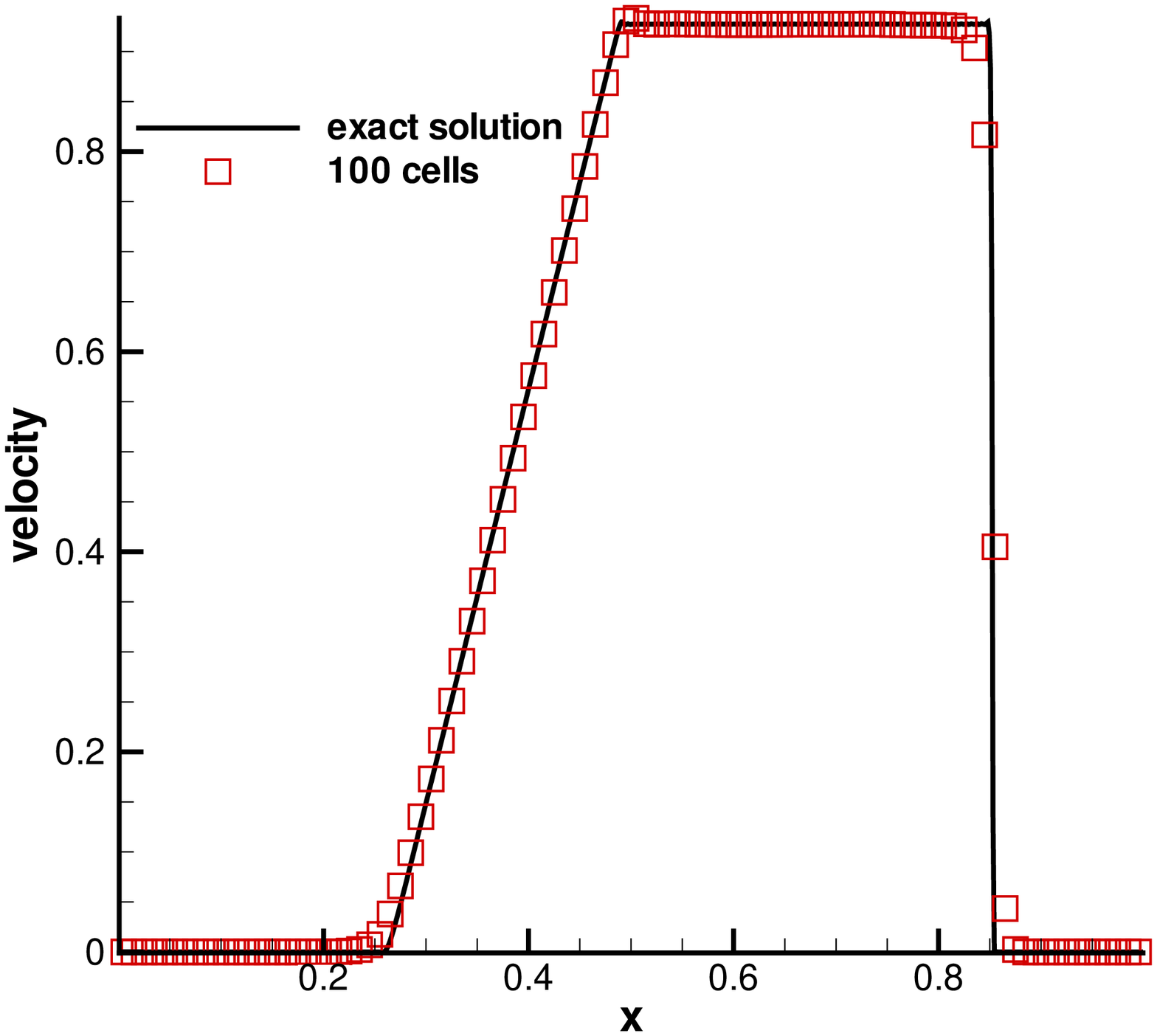}\includegraphics[width=0.444\textwidth]{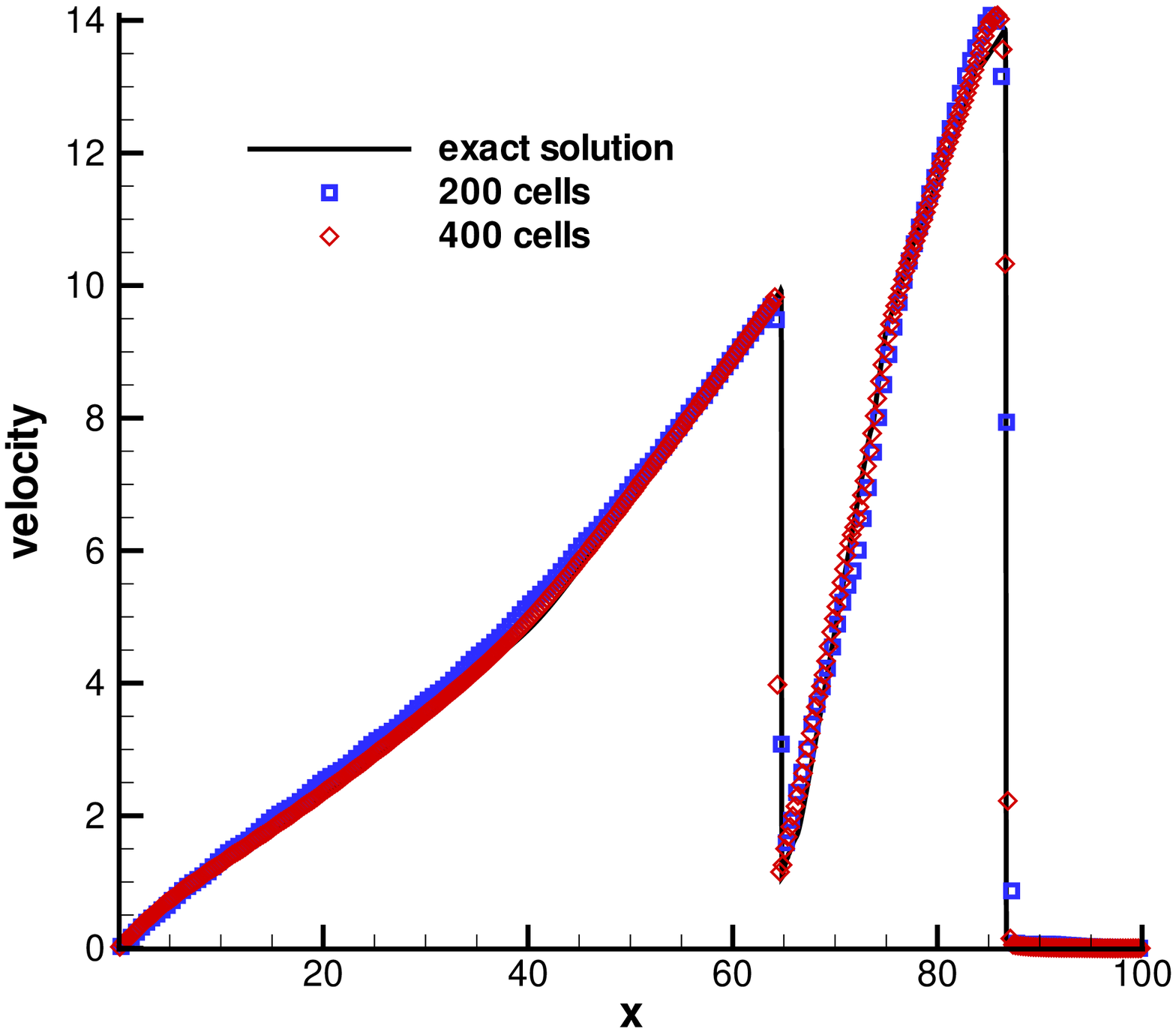}
\includegraphics[width=0.444\textwidth]{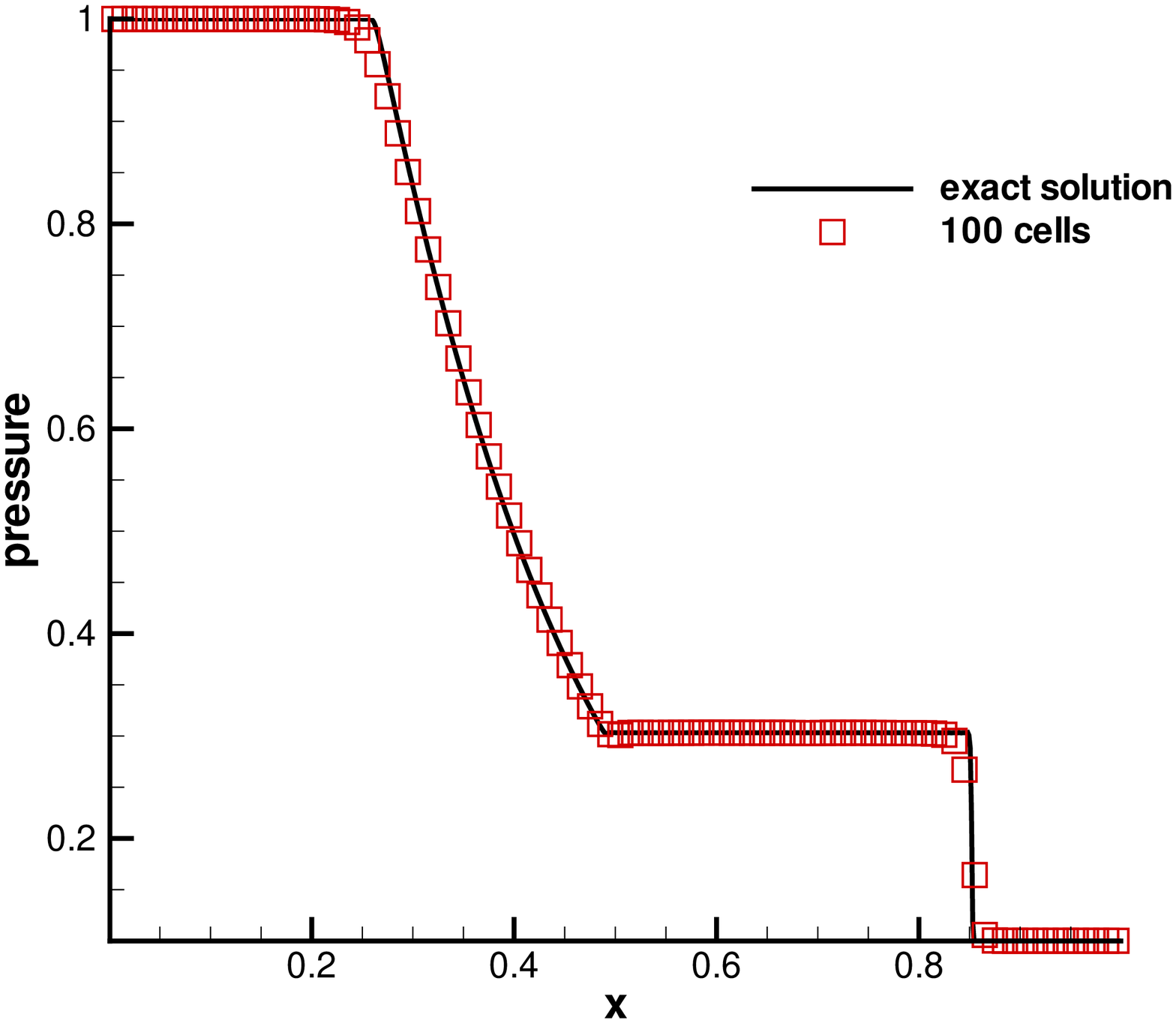}\includegraphics[width=0.444\textwidth]{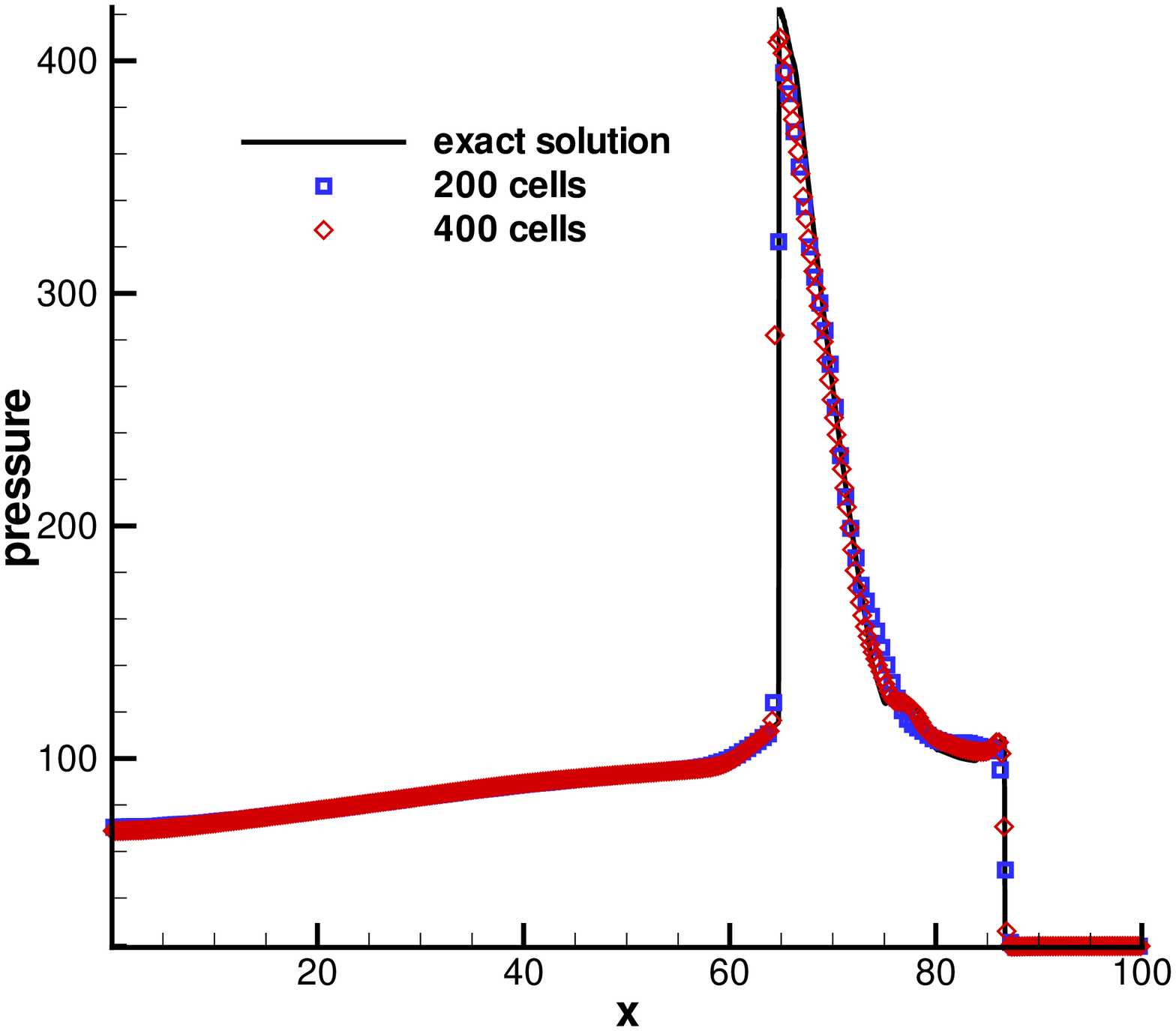}
\caption{\label{riemann-2} Sod problem (left): the density, velocity
and pressure distributions at t=0.2 with $100$ cells, and blast wave
problem (right): the density, velocity and pressure distributions at
$t=3.8$ with $200$ and $400$ cells.}
\end{figure}

\begin{figure}[!htb]
\centering
\includegraphics[width=0.444\textwidth]{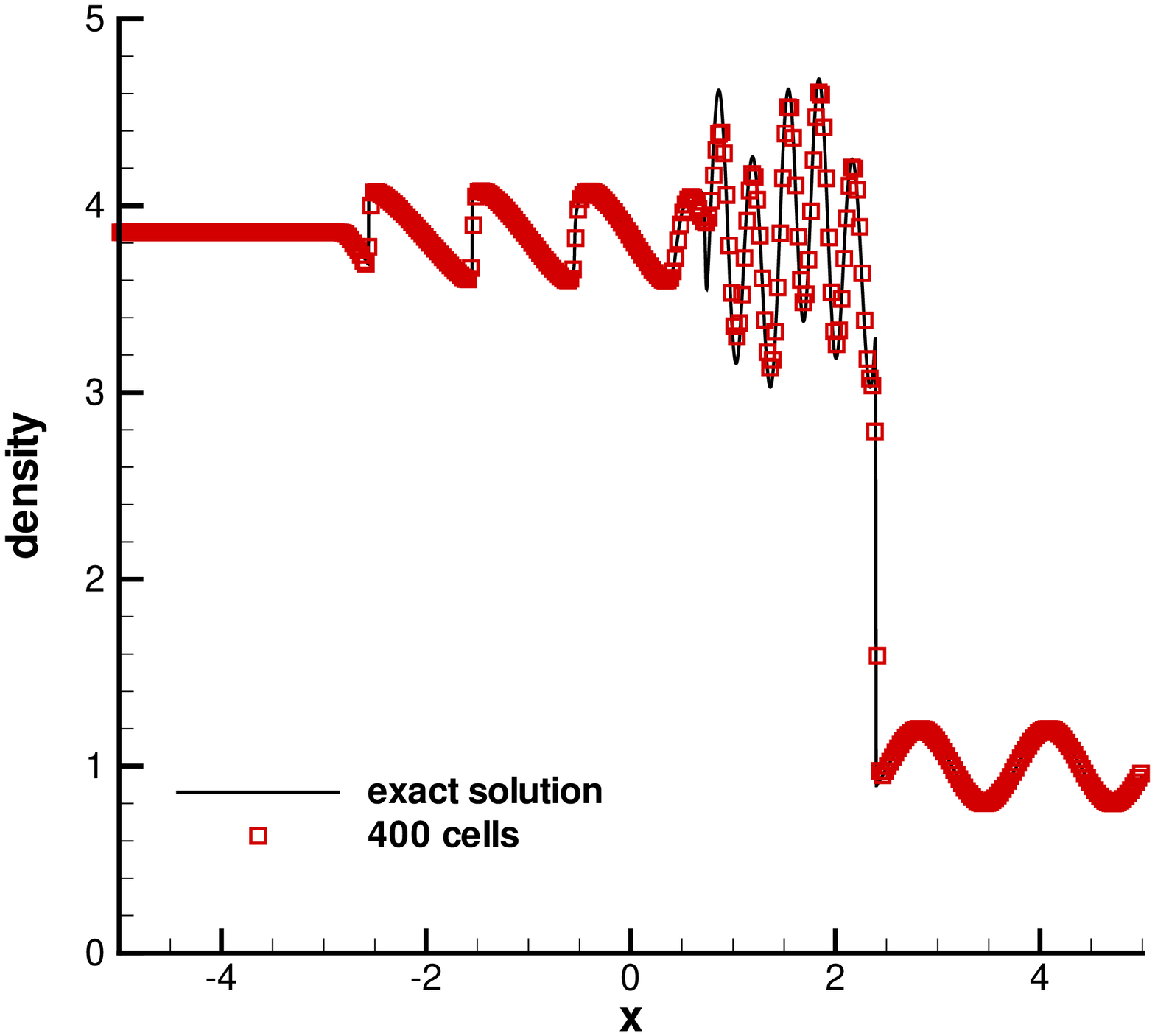}\includegraphics[width=0.444\textwidth]{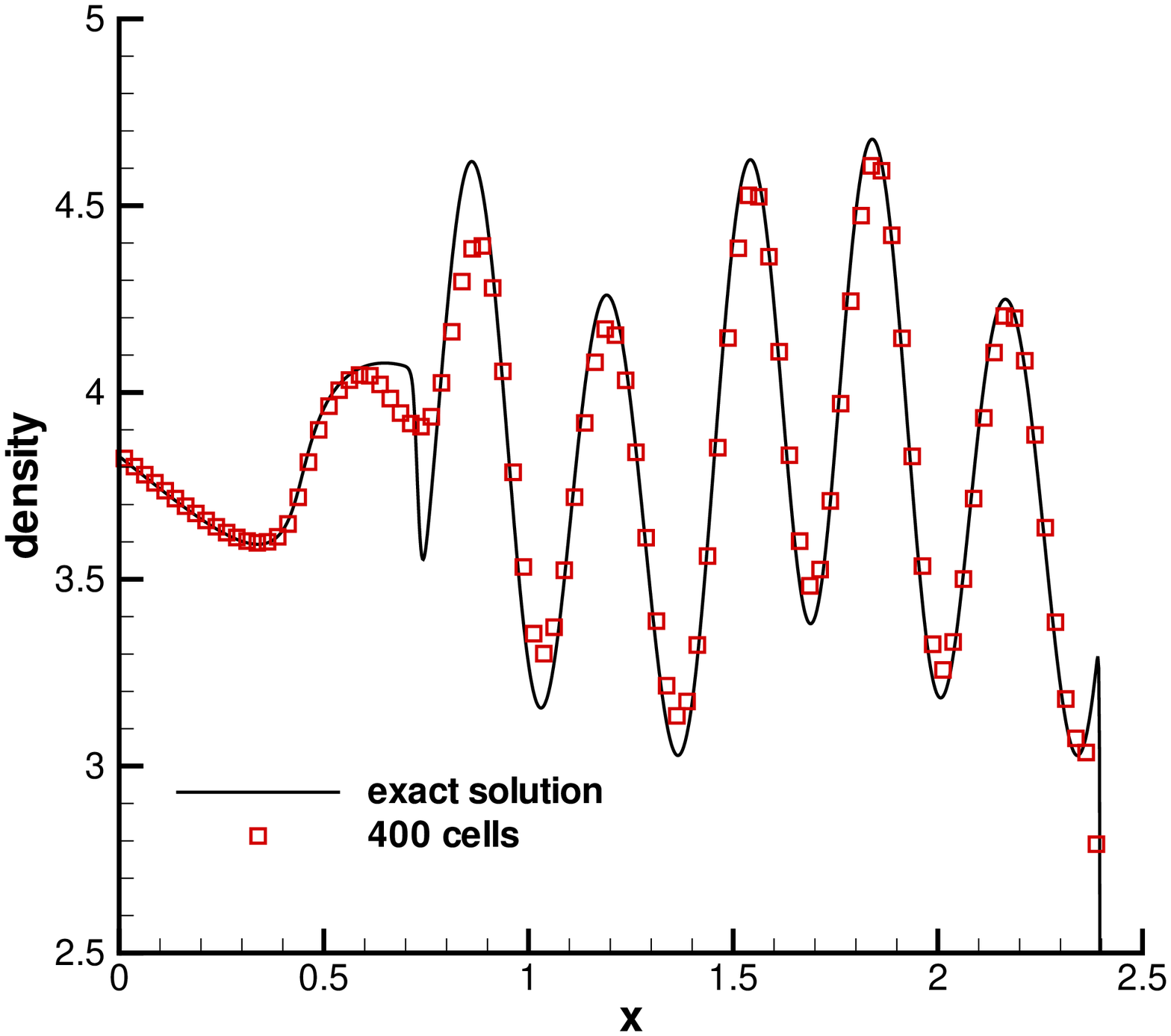}
\caption{\label{shu} Shu-Osher shock acoustic-wave interaction. Density distributios at $t=1.8$ with $400$ cells.}
\end{figure}

\subsection{One dimensional Riemann problems}
For one-dimensional case, two Riemann problems are
considered. The first one is  the Sod problem. The computational domain
is $[0,1]$ with $100$ uniform mesh points and with non-reflected
boundary condition on both ends. The initial condition is given by
\begin{equation*}
(\rho,U,p)=\left\{\begin{array}{ll}
(1, 0, 1),  \ \ \ \ &  0<x<0.5,\\
(0.125,0,0.1),  & 0.5<x<1.
\end{array} \right.
\end{equation*}
The second one is the Woodward-Colella blast wave problem \cite{Case-Woodward}. The
computational domain is $[0,100]$ with $200$ and $400$ uniform mesh
points. The reflected boundary conditions are imposed on both ends.
The initial conditions are given as follows,
\begin{equation*}
(\rho,U,p) =\left\{\begin{array}{ll}
(1, 0, 1000), \ \ \ \ & 0\leq x<10,\\
(1, 0, 0.01), & 10\leq x<90,\\
(1, 0, 100), &  90\leq x\leq 100.
\end{array} \right.
\end{equation*}
The density, velocity, and pressure distributions for the fourth-order GKS
and the exact solutions are presented in Fig.\ref{riemann-2}
for the Sod problem at $t=0.2$ and for the blast wave problem at
$t=3.8$. The numerical results agree well with the exact solutions.
The scheme can resolve the wave profiles well,
particularly for the local extreme values.

In the one-dimensional case, another standard test case is the Shu-Osher shock acoustic interaction \cite{shu-osher}.
The computational domain is $[-5,5]$ and the flow field is initialized as
\begin{equation*}
(\rho,U,p)=\left\{\begin{array}{ll}
(3.857134, 2.629369, 10.33333),  \ \ \ \ &  x \leq -4,\\
(1 + 0.2\sin (5x), 0, 1),  &  -4 <x.
\end{array} \right.
\end{equation*}
The computed density profile with $400$ mesh points at $t = 1.8$ is shown in Fig.\ref{shu}.

\begin{figure}[!h]
\centering
\includegraphics[width=0.75\textwidth]{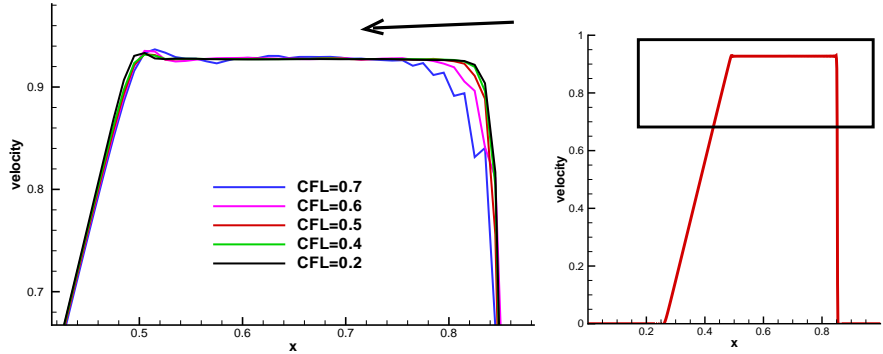}
\caption{\label{stability} The stability test for the Sod problem.}
\end{figure}

The stability for the current scheme is tested by the Sod problem.
The velocity profiles with different $CFL$ numbers from $0.2$ to
$0.7$, are shown in Fig.\ref{stability}. The scheme is basically stable under the conventional CFL condition.
The waves profiles can be well resolved at a $CFL$ number around $0.5$. In the following numerical tests, without spacial statement,
the $CFL$ number takes a fixed value of $0.4$.

\begin{figure}[!h]
\centering
\includegraphics[width=0.46\textwidth]{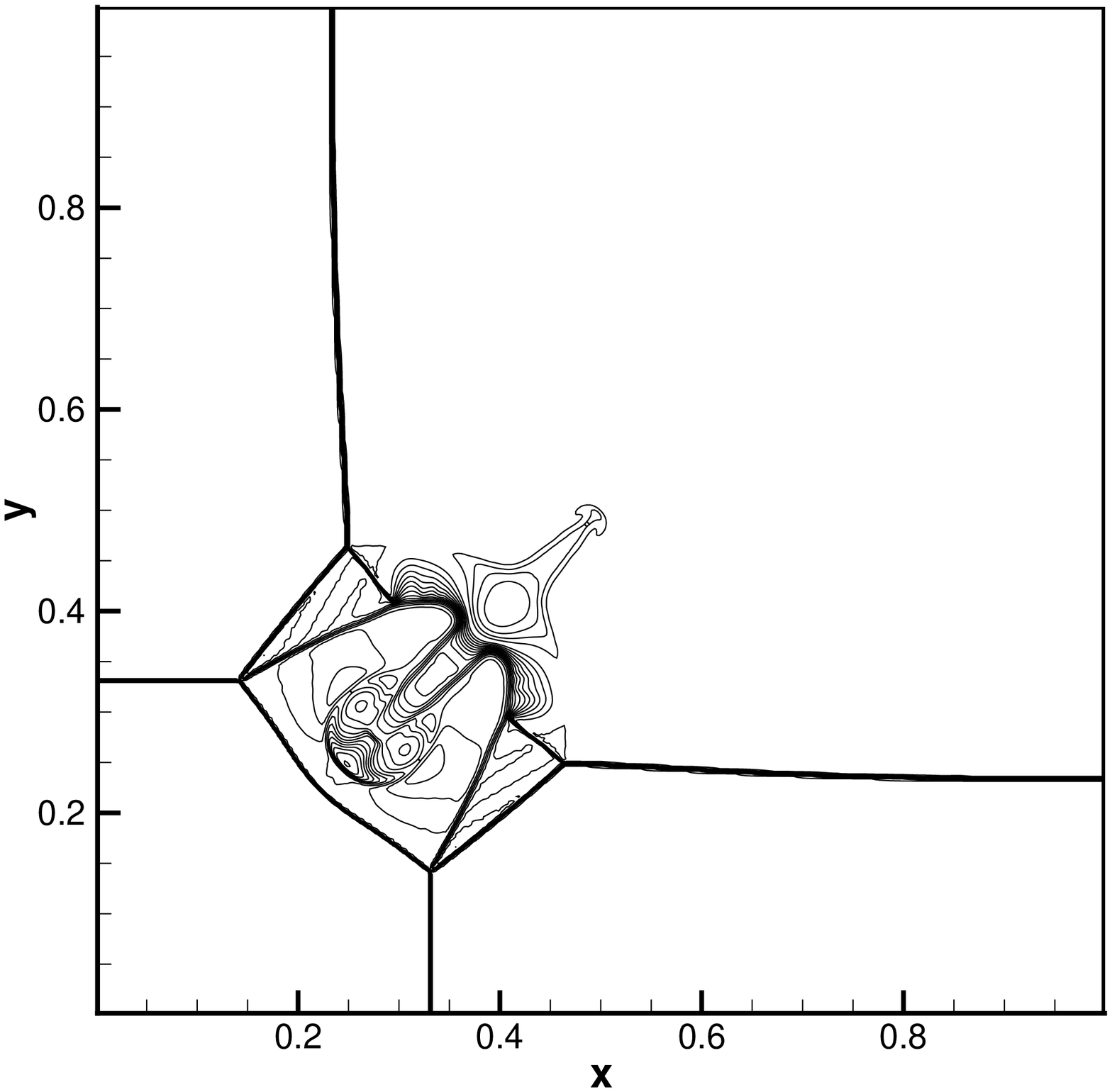}
\includegraphics[width=0.46\textwidth]{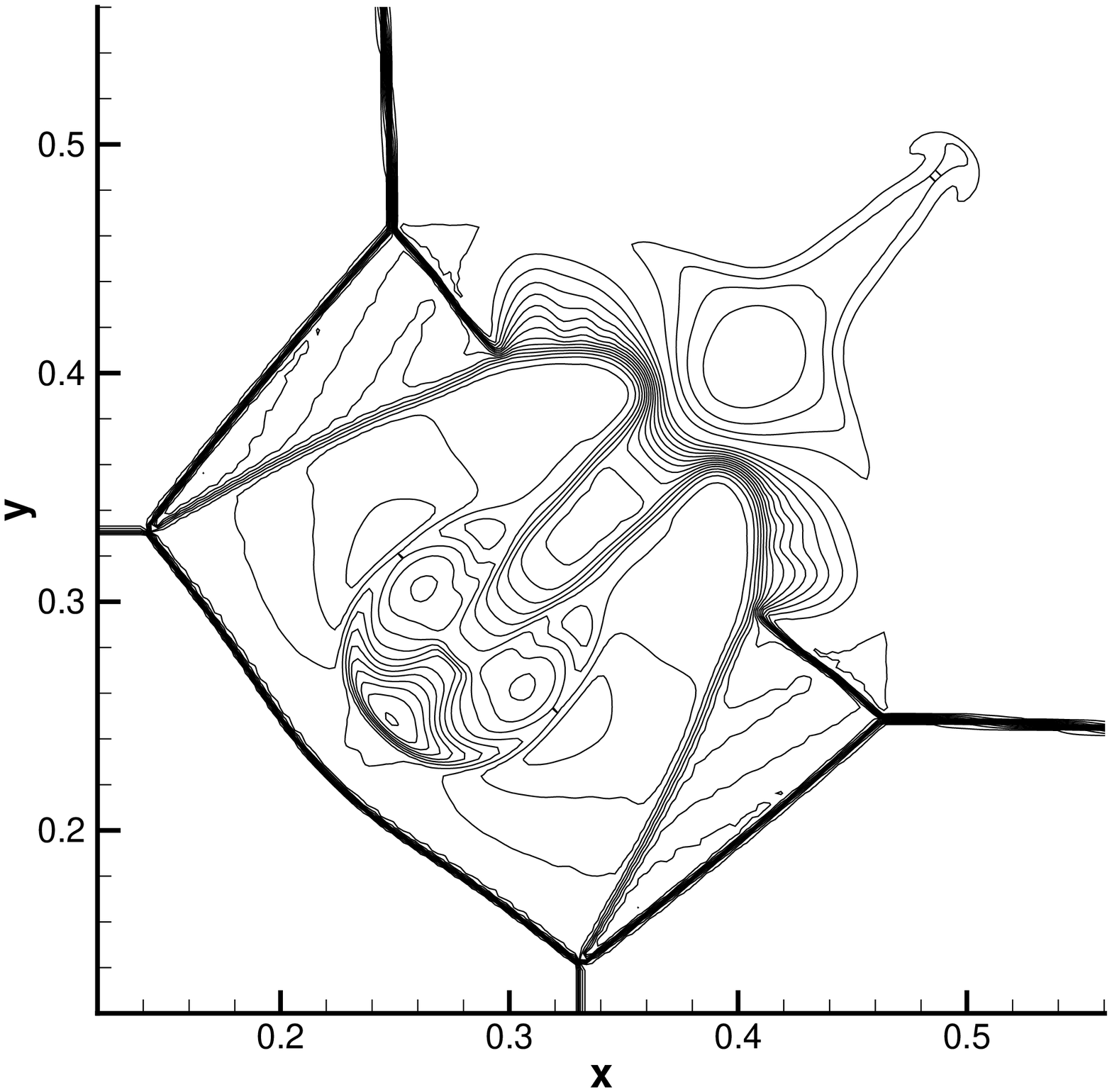}
\includegraphics[width=0.46\textwidth]{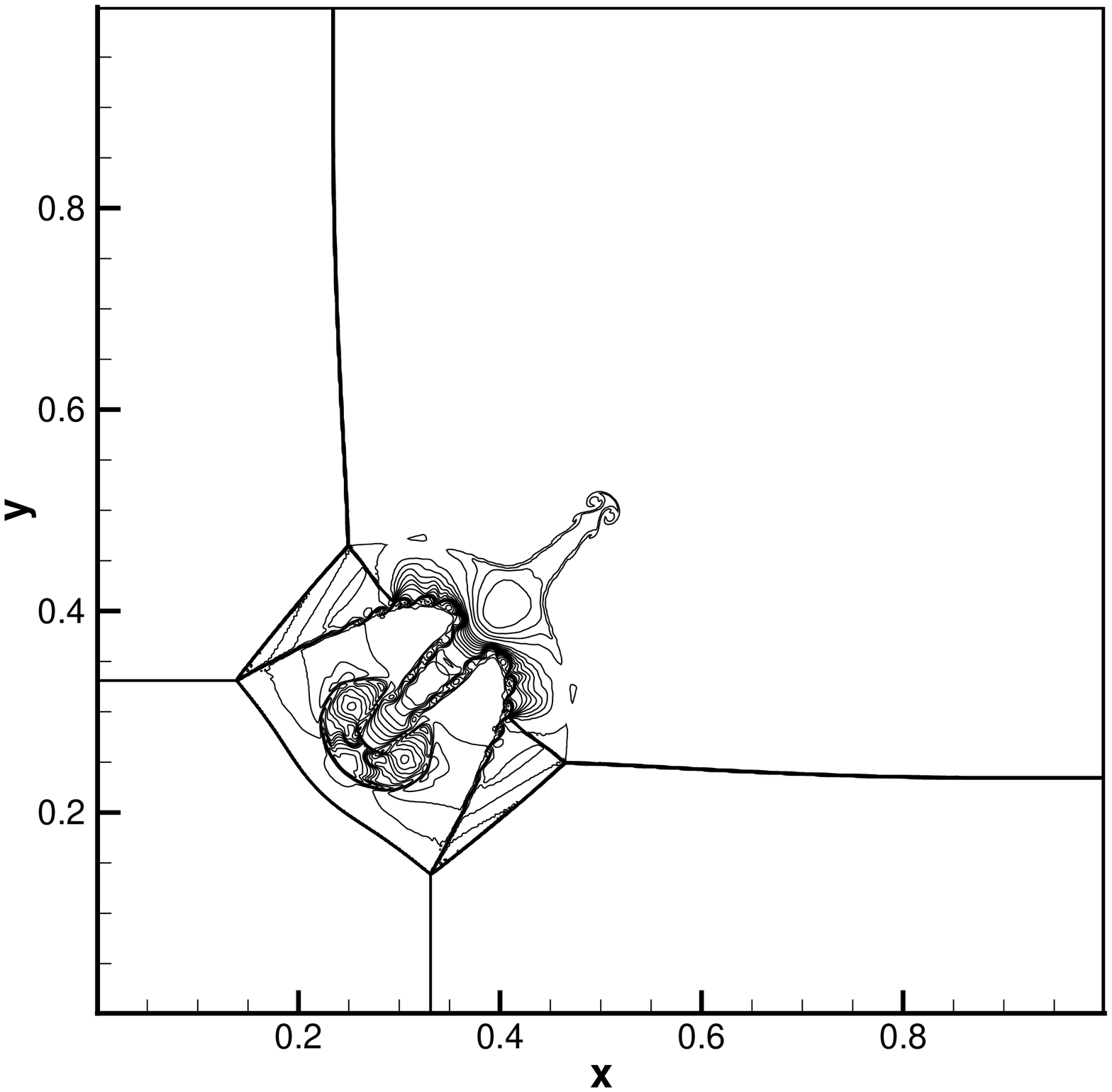}
\includegraphics[width=0.46\textwidth]{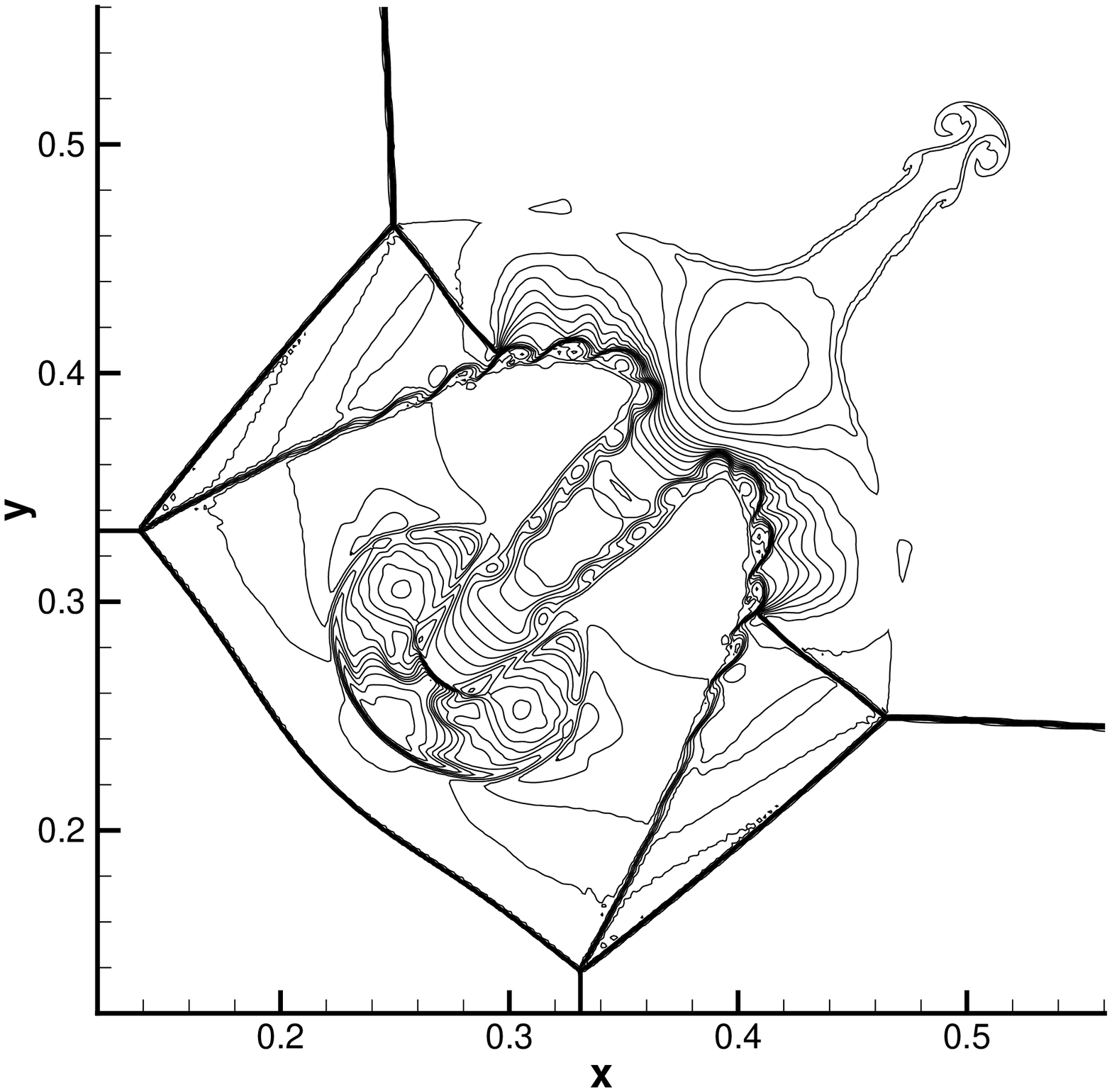}
\caption{\label{2d-riemann-1} The density distribution for the first
two-dimensional Riemann problem at $t=0.3$ with $400\times400$ (top)
and $800\times800$ (bottom) mesh points.}
\end{figure}

\begin{figure}[!h]
\centering
\includegraphics[width=0.46\textwidth]{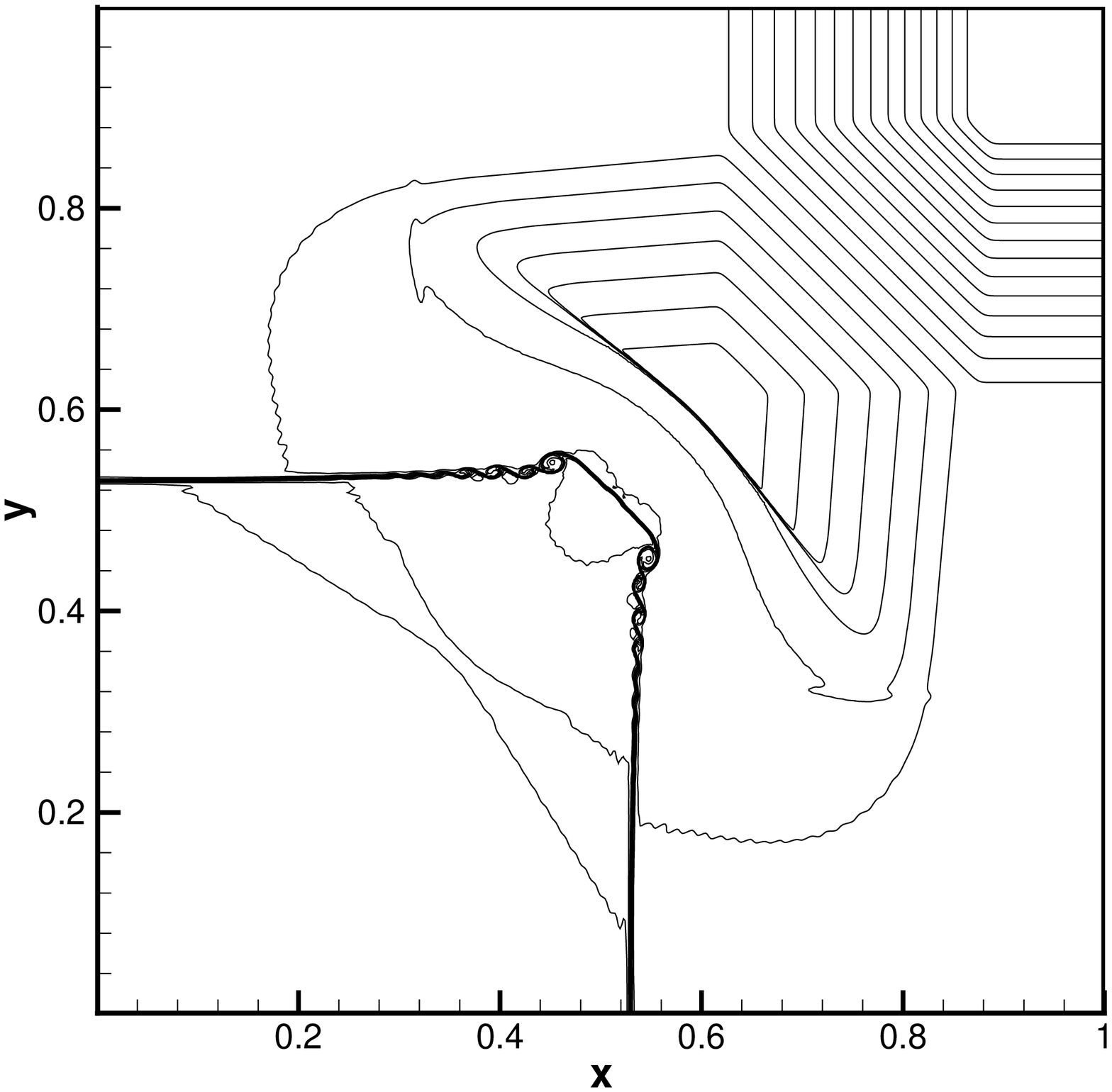}
\includegraphics[width=0.46\textwidth]{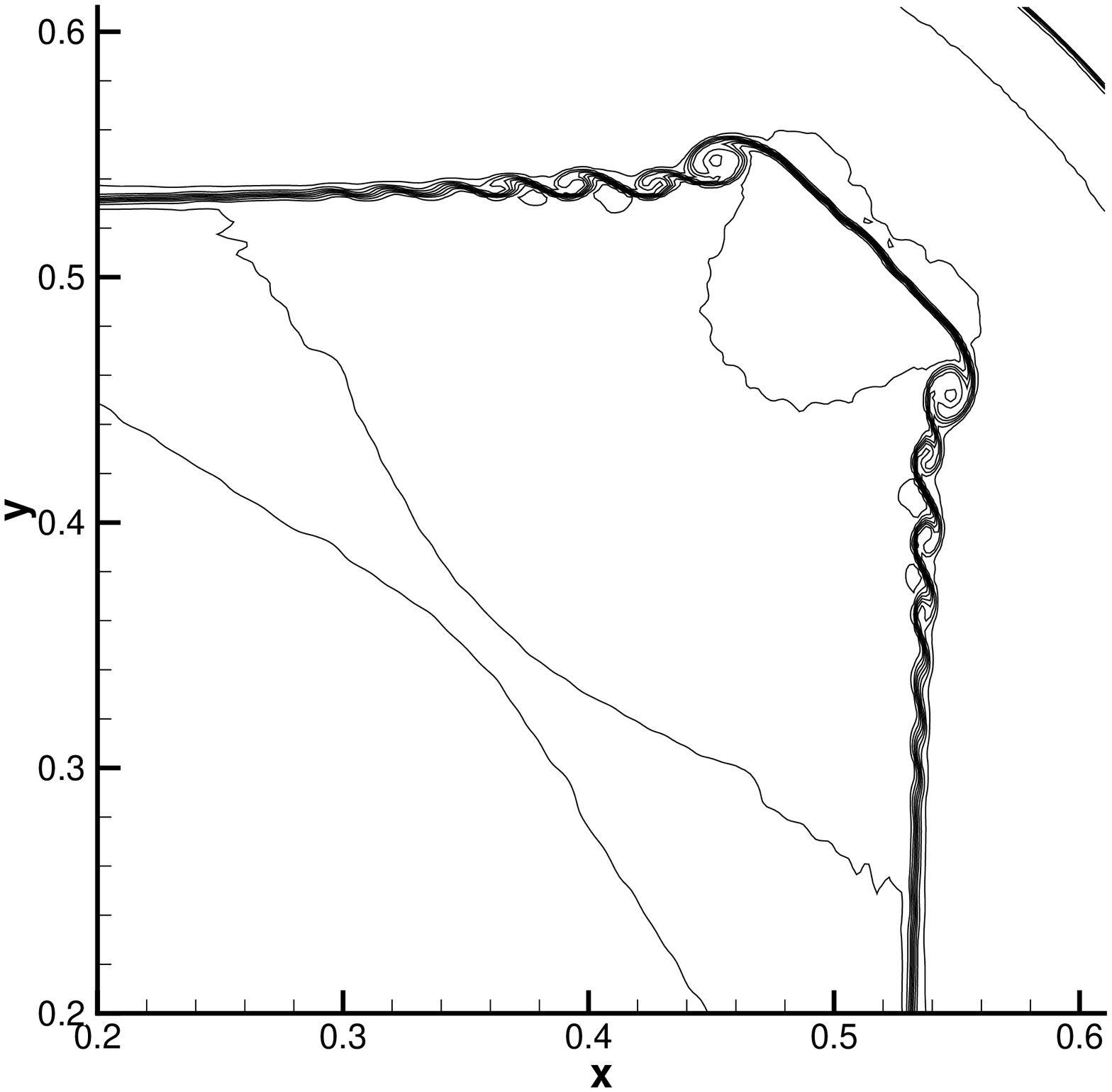}
\includegraphics[width=0.46\textwidth]{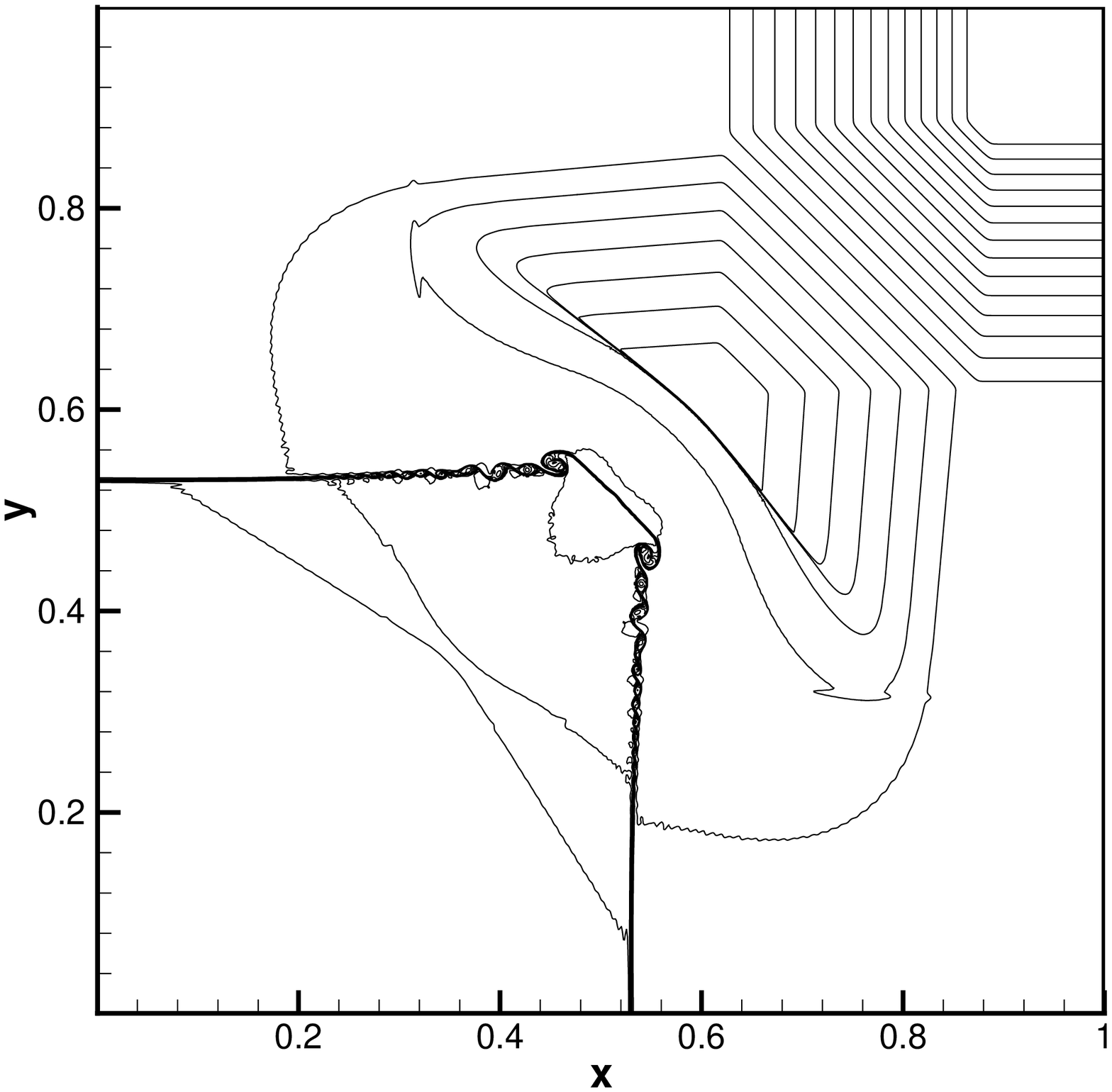}
\includegraphics[width=0.46\textwidth]{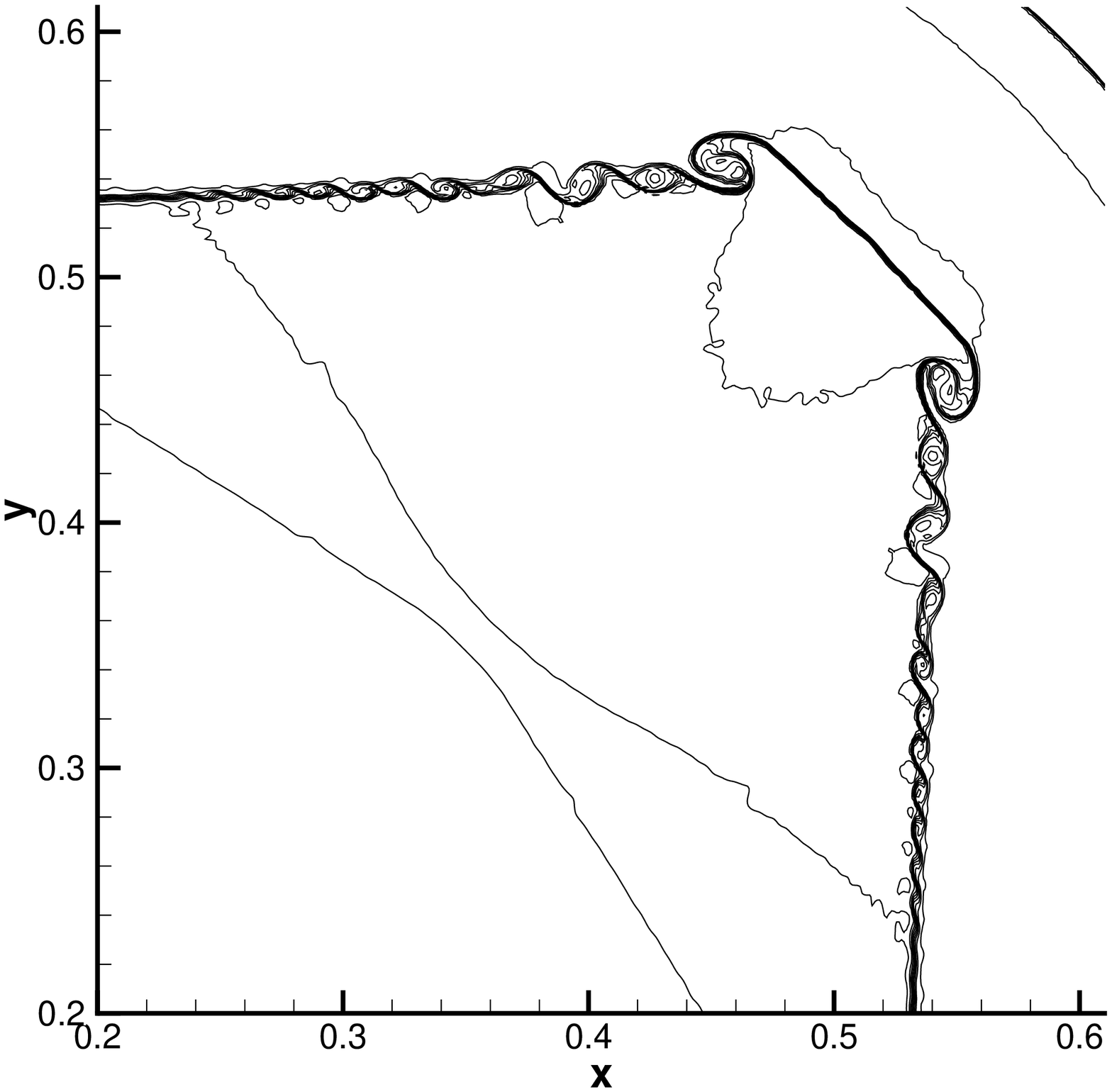}
\caption{\label{2d-riemann-2} The density distribution for the
second two-dimensional Riemann problem at $t=0.25$ with
$600\times600$ (top) and $1000\times1000$ (bottom) mesh points.}
\end{figure}

\subsection{Two-dimensional Riemann problems}

In the following, two examples of two-dimensional Riemann problems are
considered, which involve the interactions of shocks, the
interaction of shocks with vortex sheets, and the interaction of
vortices \cite{Case-lax, Li-Zhang, Han-Li}. The computational domain is
$[0,1]\times[0,1]$, and the non-reflecting boundary conditions are
used in all boundaries. The initial conditions for the first problem
are
\begin{equation*}
(\rho,U,V,p) =\left\{\begin{array}{ll}
         (1.5,0,0,1.5), \ \ \ &x>0.5,y>0.5,\\
         (0.5323,1.206,0,0.3), &x<0.5,y>0.5,\\
         (0.138,1.206,1.206,0.029), &x<0.5,y<0.5,\\
         (0.5323,0,1.206,0.3), &x>0.5,y<0.5.
                          \end{array} \right.
                          \end{equation*}
Four initial shock waves interact with each other and result in a more
complicated pattern. The density distribution and  the local
enlargement are given at $t=0.4$ in Fig.\ref{2d-riemann-1} with
$400\times400$ and $800\times800$ mesh points. From the analysis in
\cite{Case-lax}, the initial shock wave $S^-_{23}$ bifurcates at the
trip point into a reflected shock wave, a Mach stem, and a slip line.
The reflected shock wave interacts with the shock wave $S^-_{12}$ to
produce a new shock. The small scale flow structures are well
captured by the current scheme.

The initial conditions for the second case are
\begin{equation*}
(\rho,U,V,p)=\left\{\begin{aligned}
         &(1, 0.1, 0.1, 1), &x>0.5,y>0.5,\\
         &(0.5197,-0.6259, 0.1, 0.4), &x<0.5,y>0.5,\\
         &(0.8, 0.1, 0.1, 0.4), &x<0.5,y<0.5,\\
         &(0.5197,0.1,-0.6259, 0.4), &x>0.5,y<0.5.
                          \end{aligned} \right.
                          \end{equation*}
This case is to simulate the interaction of the rarefaction waves
and the vortex-sheets. The density distribution at $t=0.4$ and the local
enlargement are given  in Fig.\ref{2d-riemann-2} with
$600\times600$ and $1000\times1000$ mesh points. The roll-up is well
captured by the current scheme.
\vspace{0.2cm}

On the computational cost, the above two-dimensional Riemann problems are
tested again. As a reference, the CPU times for different schemes are
obtained with $100\times 100$ cells and $10$ time steps with Intel
Core i7-4770 CPU $@$ 3.40GHz. Based on the same WENO reconstruction,
the CPU times for the second-order GKS \cite{GKS-Xu2}, the single
stage third-order GKS \cite{GKS-high1,GKS-high2}, and the current
two-stage fourth-order GKS are given in Table \ref{tab0}, where both conservative
and characteristic reconstructions are used.

\begin{table}[!h]
\begin{center}
\def\temptablewidth{0.9\textwidth}
{\rule{\temptablewidth}{0.5pt}}
\begin{tabular*}{\temptablewidth}{@{\extracolsep{\fill}}c|cccc}
variable & 2nd-order GKS & 3rd-order GKS ~ & 4th-order GKS  \\
\hline
conservative & 0.704893$s$ & 1.24681$s$&  1.95370$s$\\
characteristic & 0.842873$s$ & 1.38178$s$ &  2.20566$s$
\end{tabular*}
{\rule{\temptablewidth}{0.5pt}}
\end{center}
\vspace{-4mm} \caption{\label{tab0} The test of the computational
cost for different schemes.}
\end{table}

For the fourth-order GKS, the computational cost is about $3$ times of that of second-order scheme. Since the fourth-order GKS has three Gauss points for flux evaluation and
two stages, the three times computational time difference means that the reconstruction makes great contribution to the computational cost as well,
because in terms of reconstruction cost the fourth-order scheme only takes about two times of computational cost of second-order scheme.
However, in order to get the same accuracy as that of the fourth-order scheme, the second-order method needs refine the mesh, at least once.
In the two-dimensional calculation, the computational cost for one mesh refinement will be increased by $8$ times.
Therefore, the fourth-order scheme is more efficient than the second order method.
Since higher-order scheme does have advantages in comparison with lower order method, it is worth
to construct the two-stage fifth-order GKS from the third-order GKS flux function, see Appendix.
Besides the computational cost, another important property of the fourth-order GKS is its accuracy and robustness. As tested in the above cases  and all cases in the following,
it clearly indicates that the fourth-order scheme is very accurate and is as robust as the second-order one. The accuracy of the scheme is closely related to the
higher-order gas evolution model,  multi-dimensionality for the inclusion of both normal and tangential derivatives around a cell interface, and the
unified treatment of the inviscid and viscous terms. A fundamental reason for the robustness of the scheme is its first-order relaxation model,
where the system is fully hyperbolic with local source term.

\begin{figure}[!h]
\centering
\includegraphics[width=0.449\textwidth]{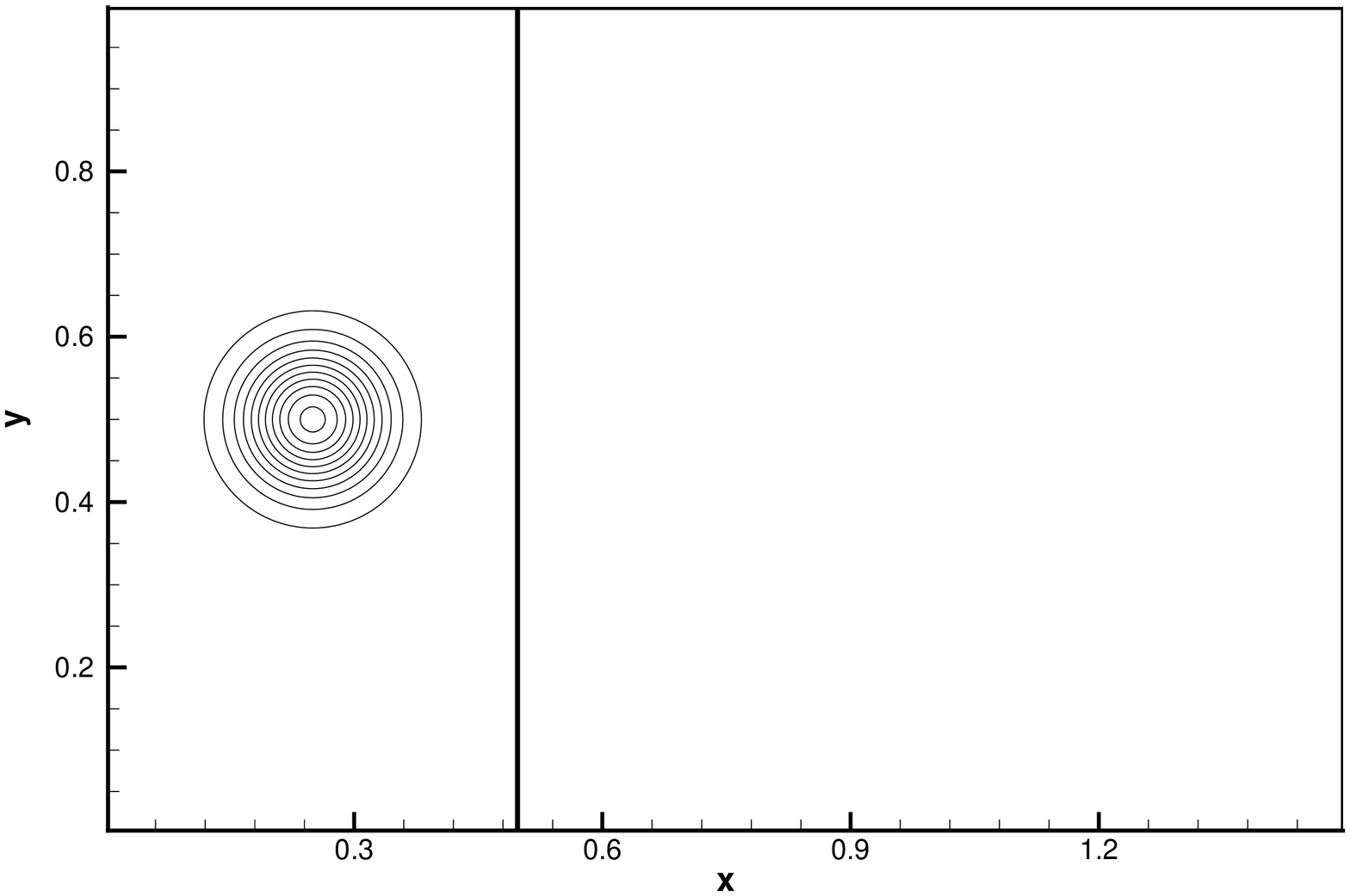}
\includegraphics[width=0.449\textwidth]{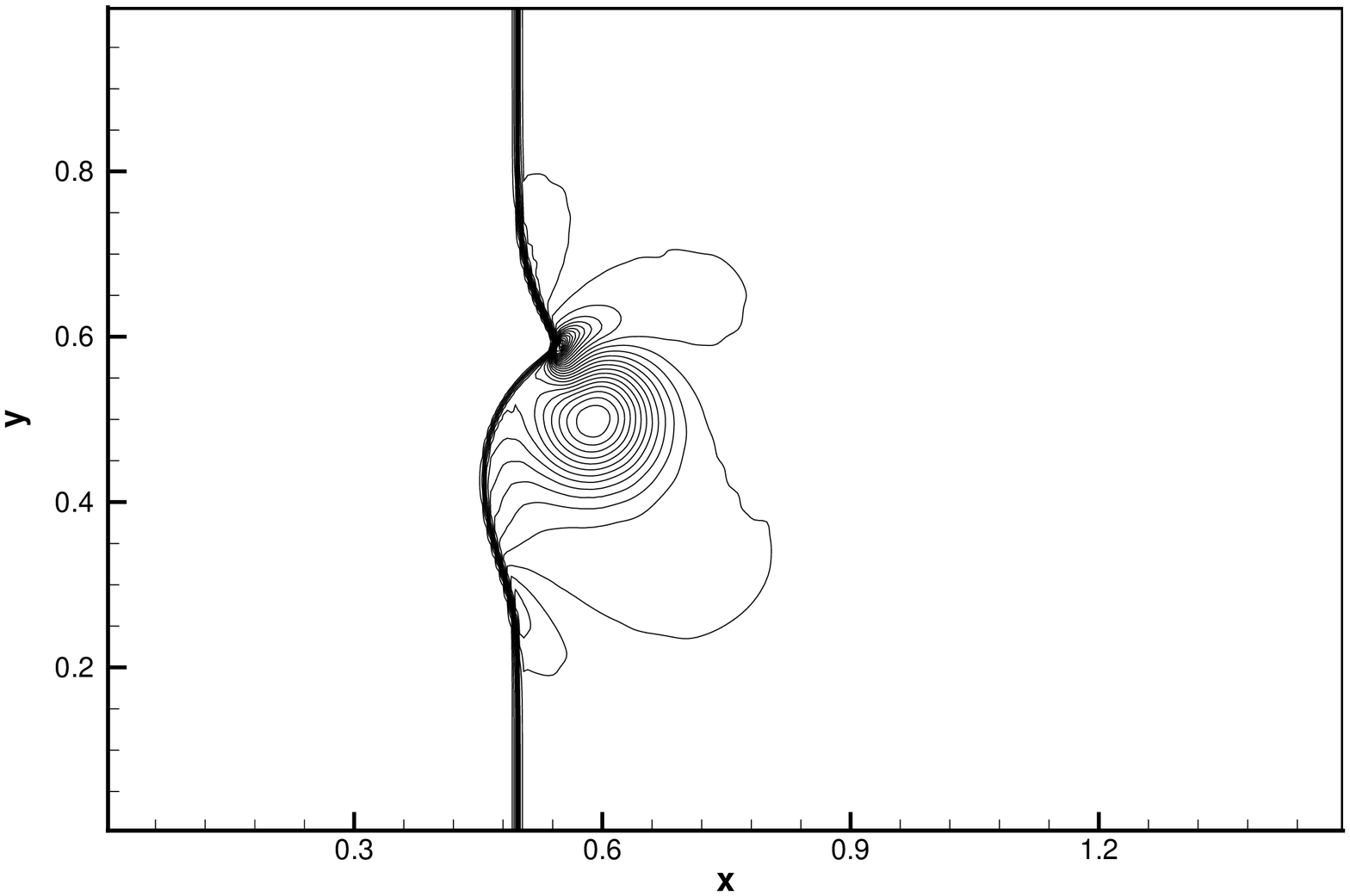}
\includegraphics[width=0.449\textwidth]{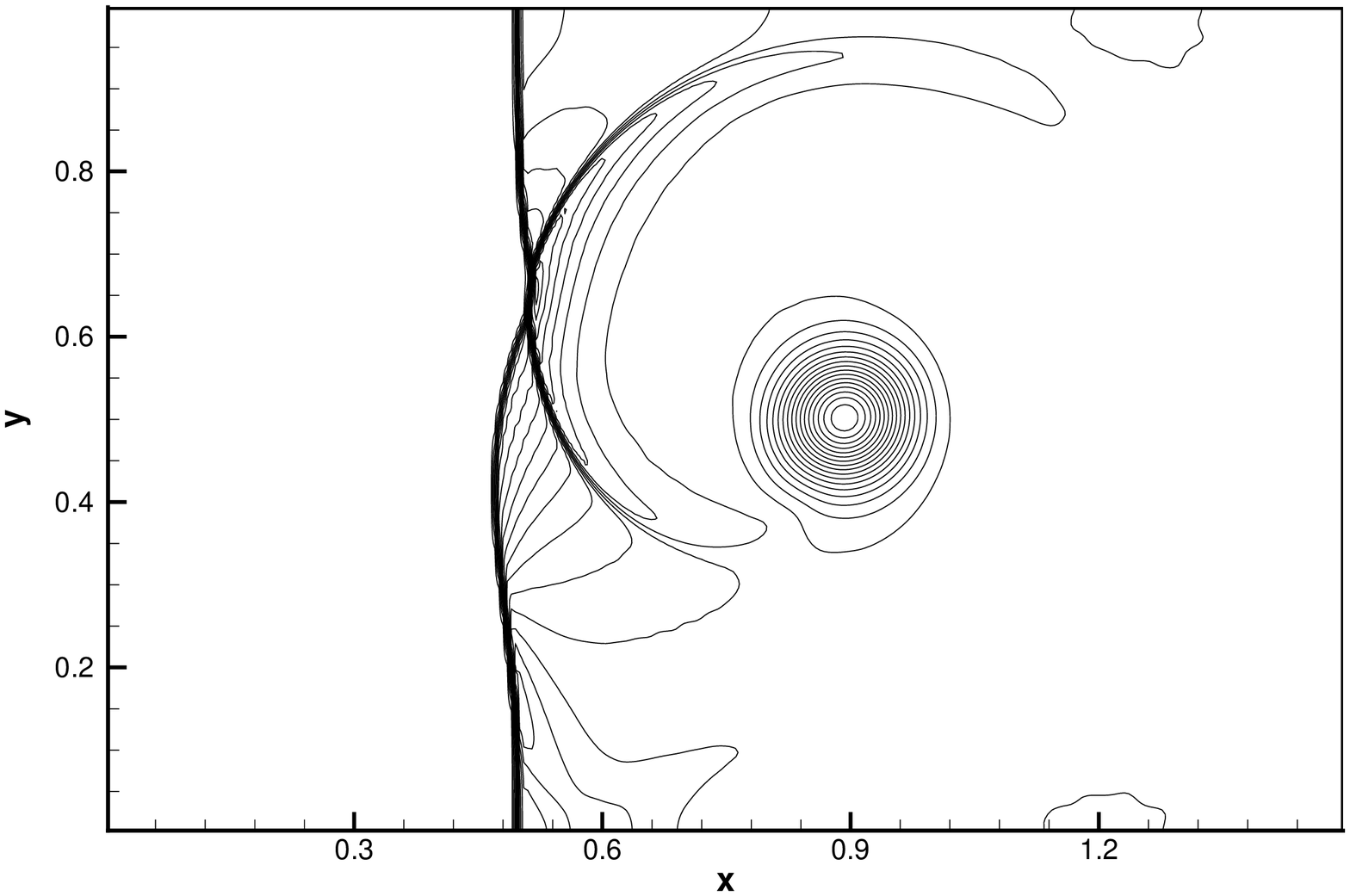}
\includegraphics[width=0.449\textwidth]{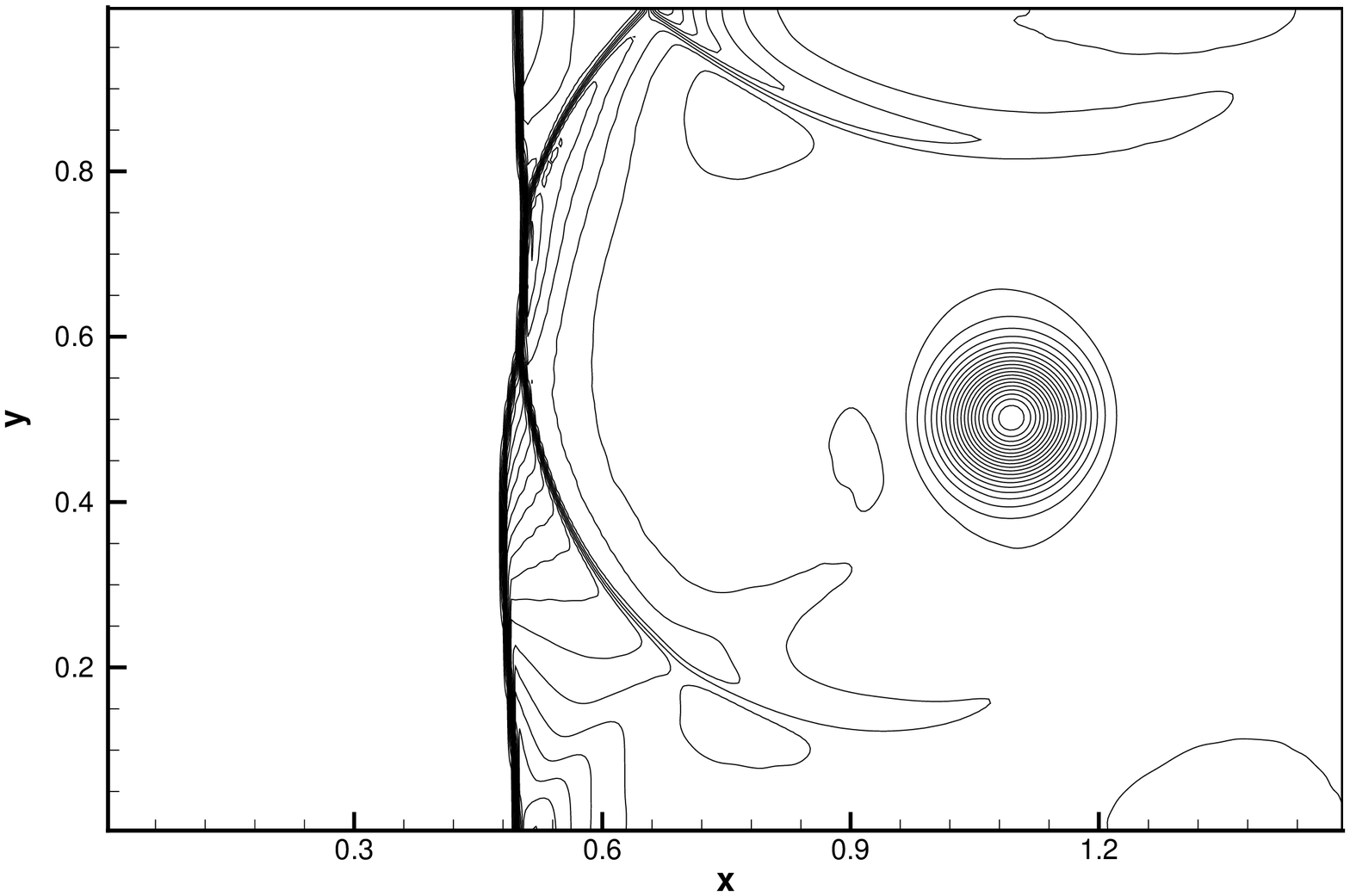}
\caption{\label{shock-vortex1}Shock vortex interaction: the pressure
distribution at $t=0, 0.3, 0.6$ and $0.8$ with mesh size $\Delta
x=\Delta y=1/200$.}  \centering
\includegraphics[width=0.6\textwidth]{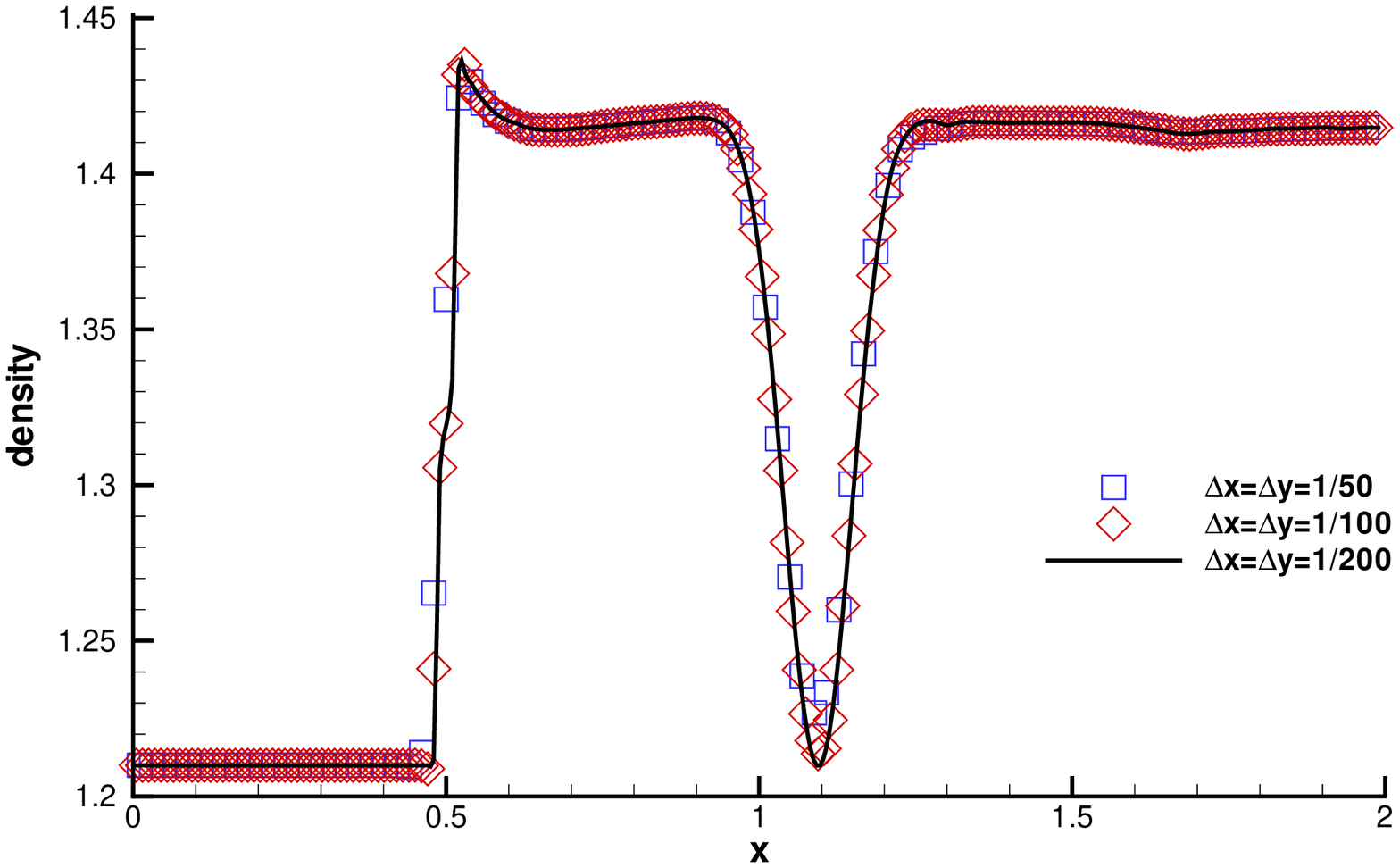}
\caption{\label{shock-vortex2}Shock vortex interaction: the density
distribution at $t=0.8$ along the horizontal symmetric line $y=0.5$
with mesh size $\Delta x=\Delta y=1/50, 1/100$ and $1/200$.}
\end{figure}

\subsection{Shock vortex interaction}
The interaction between a stationary shock and a vortex for the
inviscid flow \cite{WENO} is presented. The computational domain is
taken to be $[0, 2]\times[0, 1]$. A stationary Mach $1.1$ shock is
positioned at $x=0.5$ and normal to the $x$-axis. The left upstream
state is $(\rho, U, V, p) = (Ma^2,\sqrt{\gamma}, 0, 1)$, where $Ma$
is the Mach number. A small vortex is obtained through a
perturbation on the mean flow with the velocity $(U, V)$,
temperature $T=p/\rho$, and entropy $S=\ln(p/\rho^\gamma)$.  The
perturbation is expressed as
\begin{align*}
&(\delta U,\delta V)=\kappa\eta e^{\mu(1-\eta^2)}(\sin\theta,-\cos\theta),\\
&\delta
T=-\frac{(\gamma-1)\kappa^2}{4\mu\gamma}e^{2\mu(1-\eta^2)},\delta
S=0,
\end{align*}
where $\eta=r/r_c$, $r=\sqrt{(x-x_c)^2+(y-y_c)^2}$, and $(x_c,
y_c)=(0.25, 0.5)$ is the center of the vortex. Here $\kappa$ indicates
the strength of the vortex, $\mu$ controls the decay rate of the
vortex, and $r_c$ is the critical radius for which the vortex has
the maximum strength. In the computation, $\kappa=0.3$, $\mu=0.204$,
and $r_c=0.05$. The reflected boundary conditions are used on the
top and bottom boundaries. The pressure distributions with mesh size
$\Delta x=\Delta y=1/200$ at $t=0, 0.3, 0.6$ and $0.8$ are shown in
Fig.\ref{shock-vortex1}. By $t=0.8$, one branch of the shock
bifurcations has reached the top boundary and been reflected. The reflection is well captured. The detailed density distributions
along the center horizontal line with mesh size $\Delta x=\Delta
y=1/50, 1/100$ and $1/200$ at $t=0.8$ are shown in
Fig.\ref{shock-vortex2}. The accuracy of the scheme is well demonstrated.

\begin{figure}[!ht]
\centering
\includegraphics[width=0.6\textwidth]{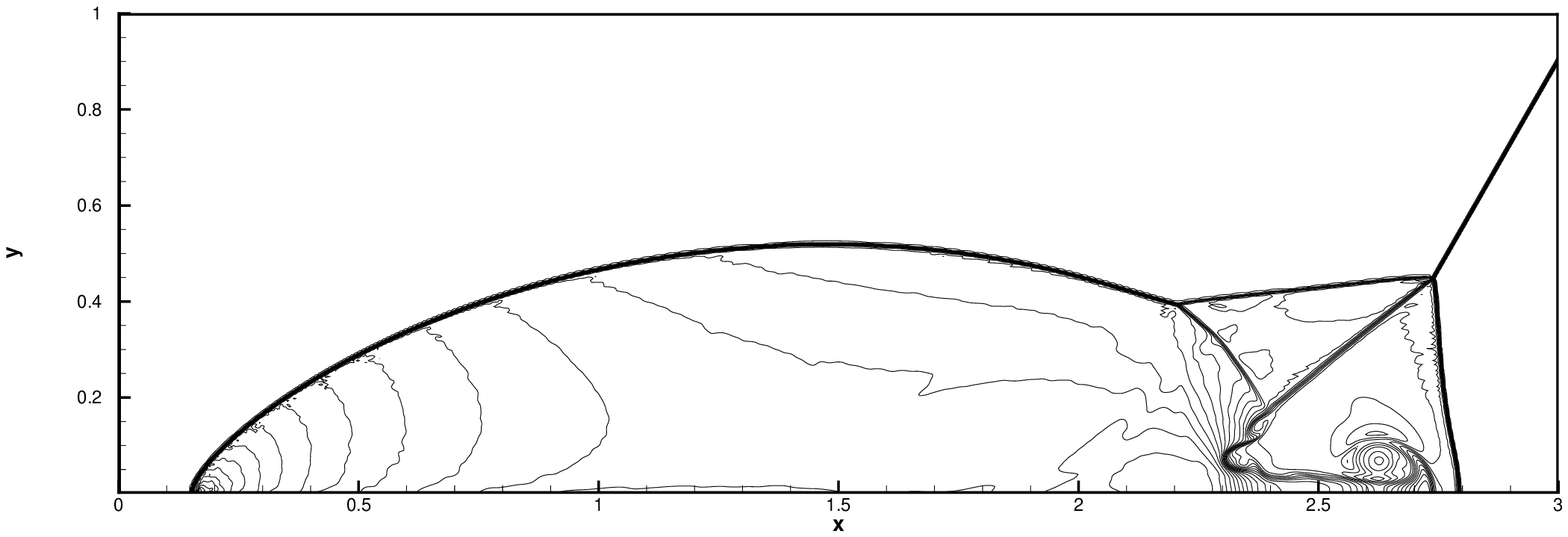}
\includegraphics[width=0.6\textwidth]{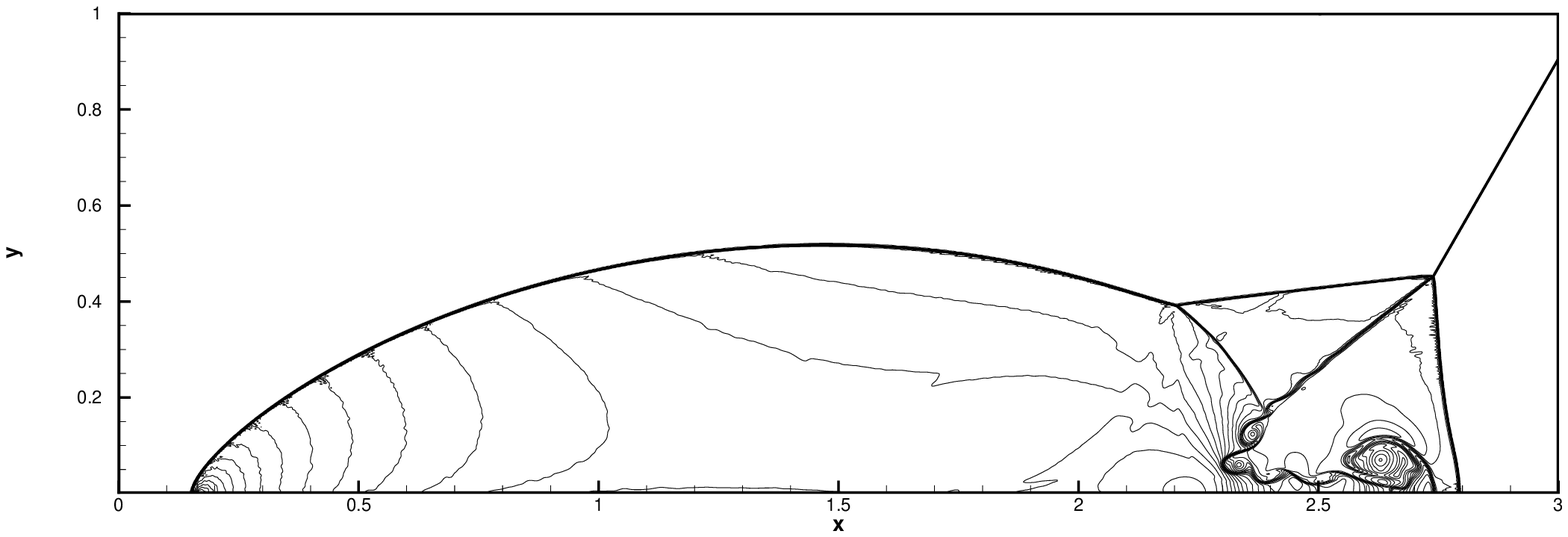}
\caption{\label{double-mach-1} Double Mach reflection: density
contours with the $720\times240$ (top) and $1440\times480$ (bottom)
mesh points.}
\includegraphics[width=0.4\textwidth]{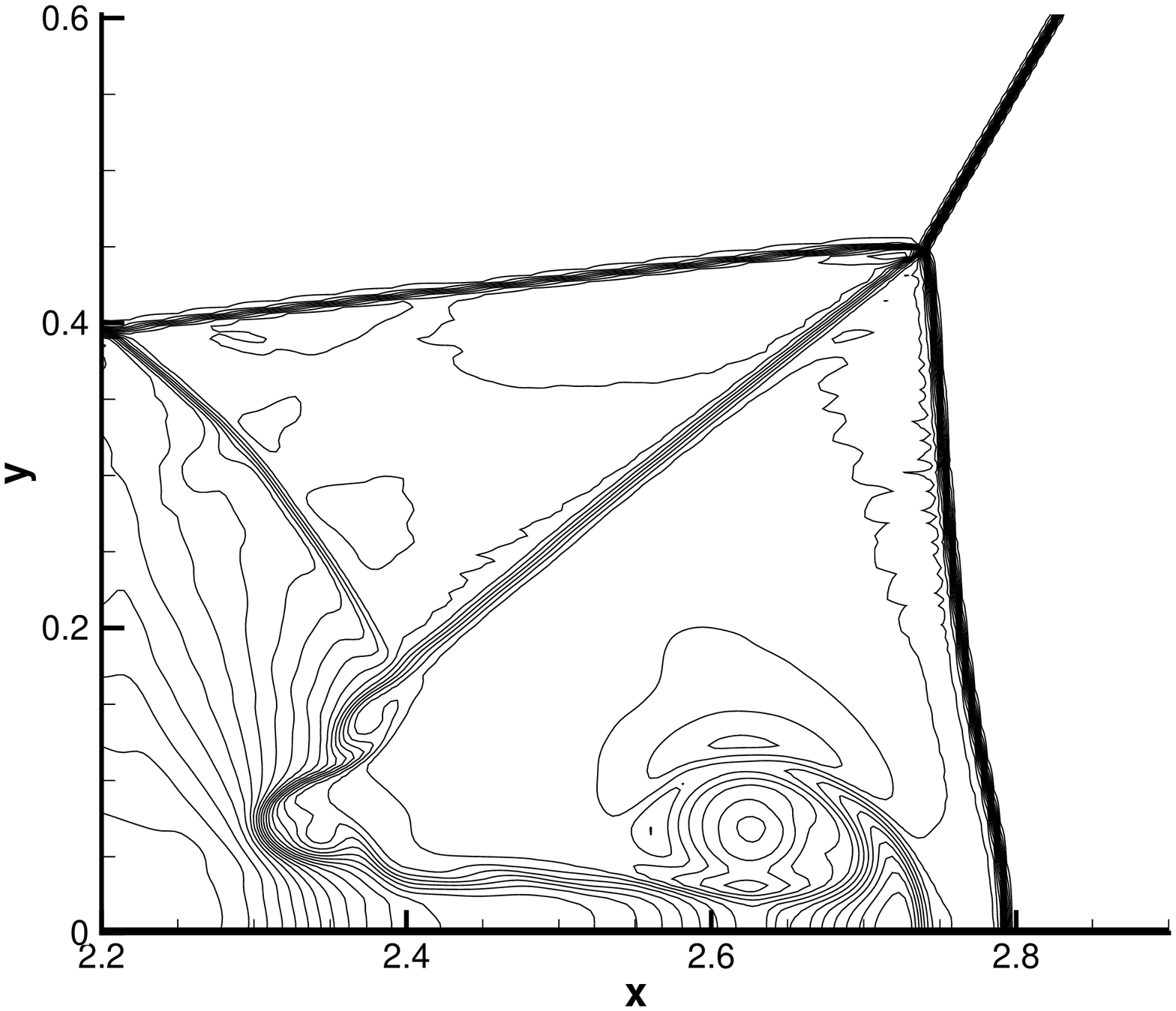}
\includegraphics[width=0.4\textwidth]{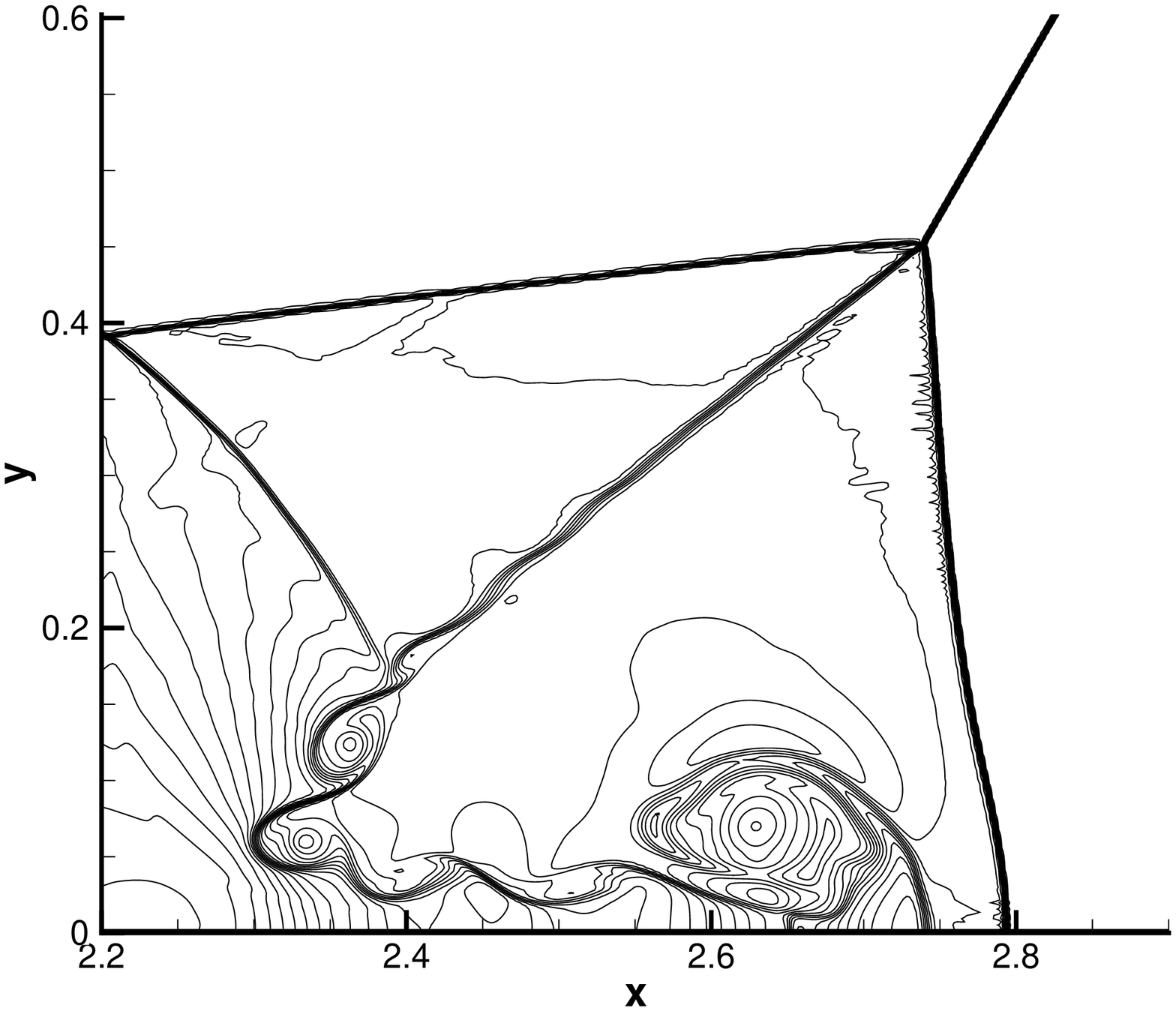}
\caption{\label{double-mach-2} Double Mach reflection: enlarged
density distributions around the triple point with the
$720\times240$ (left) and $1440\times480$ (right) mesh points.}
\end{figure}

\subsection{Double Mach reflection problem}
This problem was extensively studied by Woodward and Colella \cite{Case-Woodward} for the inviscid flow. The computational domain
is $[0,4]\times[0,1]$, and a solid wall lies at the bottom of the
computational domain starting from $x =1/6$. Initially a
right-moving Mach 10 shock is positioned at $(x,y)=(1/6, 0)$, and
makes a $60^\circ$ angle with the x-axis. The initial pre-shock and
post-shock conditions are
\begin{align*}
(\rho, U, V, p)&=(8, 4.125\sqrt{3}, -4.125,
116.5),\\
(\rho, U, V, p)&=(1.4, 0, 0, 1).
\end{align*}
The reflective boundary condition is used at the wall, while for the
rest of bottom boundary, the exact post-shock condition is imposed.
At the top boundary, the flow values are set to describe the exact
motion of the Mach $10$ shock. The density distributions with
$720\times240$ and $1440\times480$ uniform mesh points at $t=0.2$
are shown in Fig.\ref{double-mach-1} and Fig.\ref{double-mach-2}, respectively.
The current scheme resolves the flow structure under the
triple Mach stem clearly.
The amplitude of the oscillation of the slip line from the triple point is not as large as that from many other higher-order schemes,
because the GKS is intrinsically solving the NS equations and the physical dissipation will stabilize the shear instability.

\begin{figure}[!h]
\centering
\includegraphics[width=0.14\textwidth]{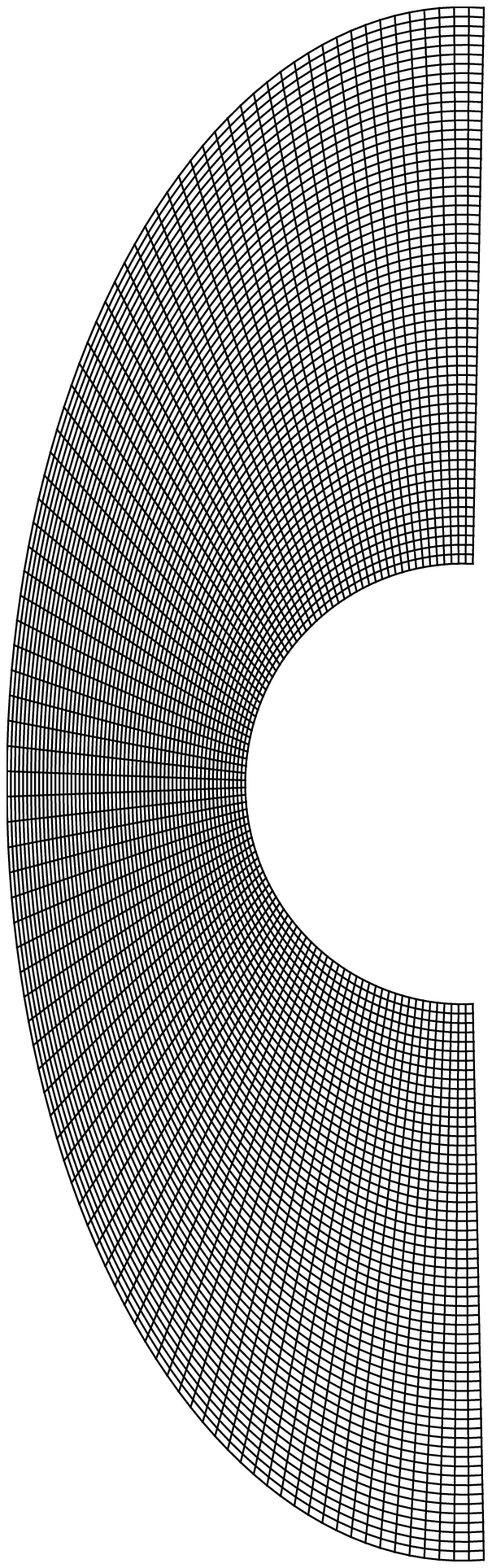}
\includegraphics[width=0.14\textwidth]{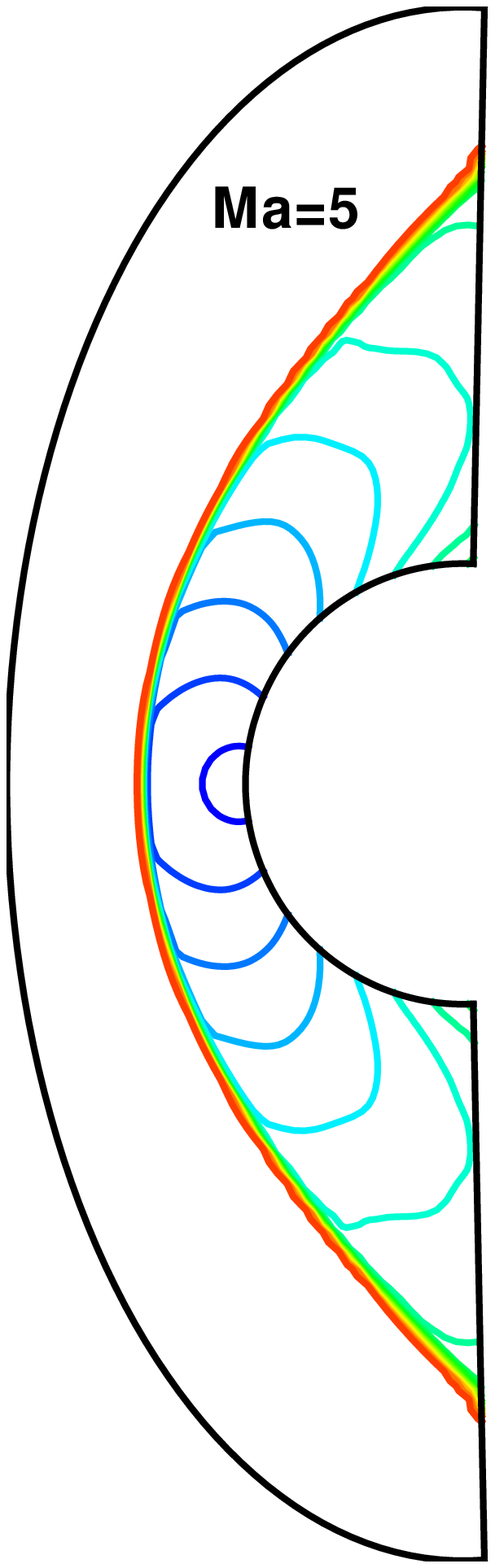}
\includegraphics[width=0.14\textwidth]{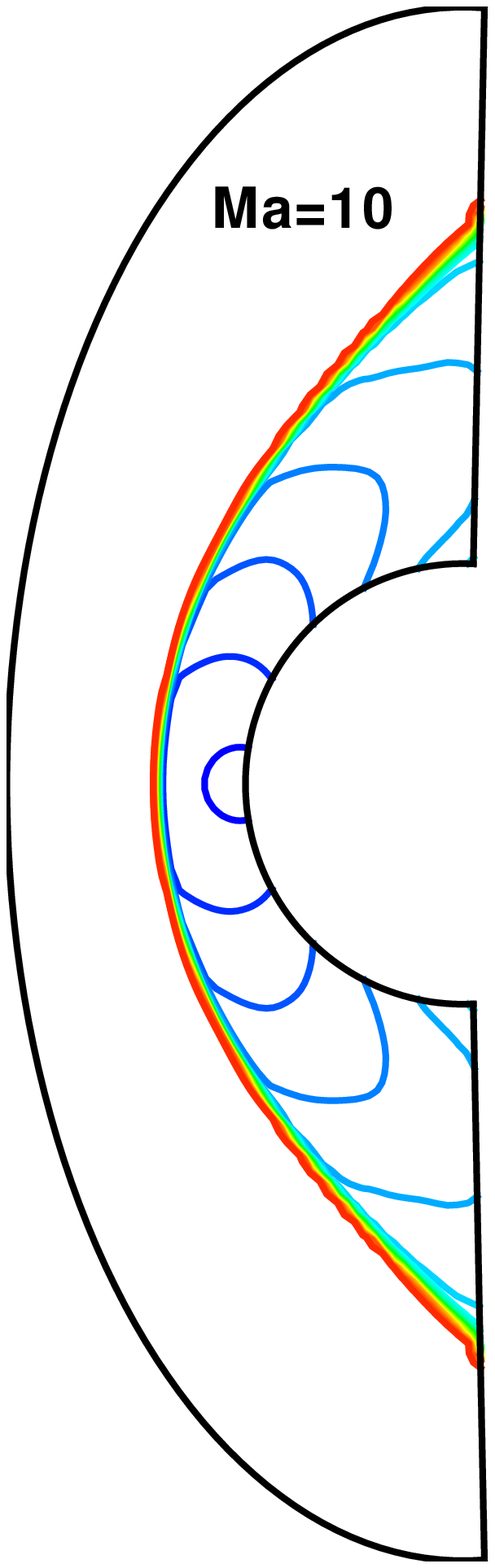}
\includegraphics[width=0.14\textwidth]{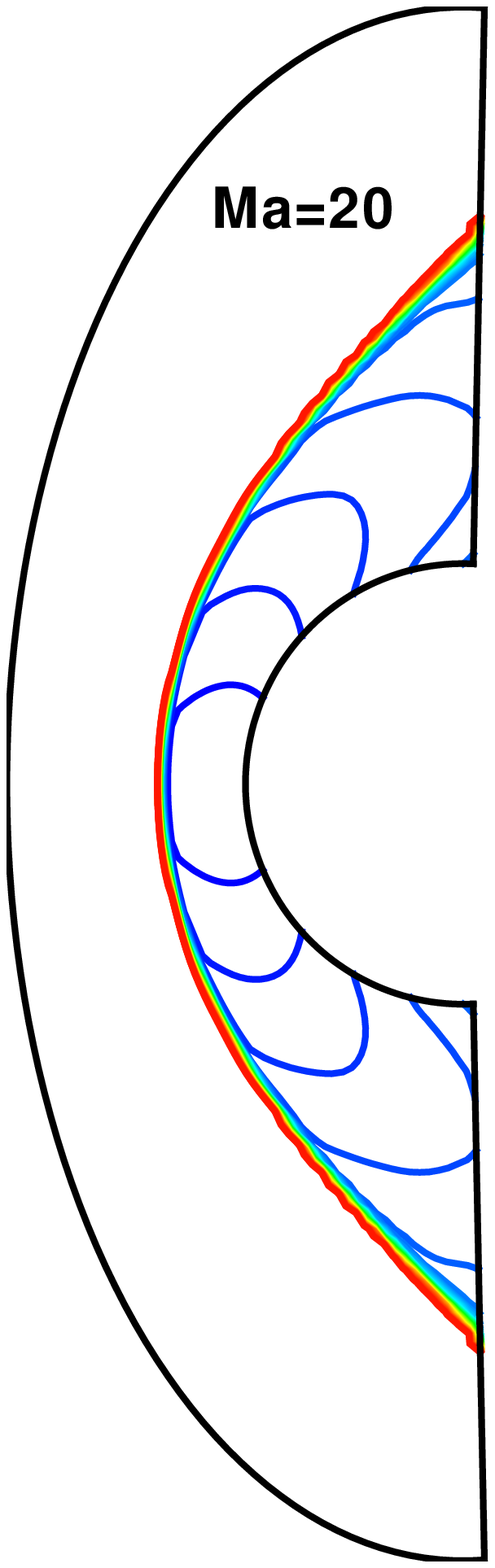}
\caption{\label{cylinder-1}Hypersonic inviscid flow past a cylinder:
the Mach number distributions for the flow with Mach number $Ma=5$,
$10$ and $20$. }
\end{figure}

\subsection{Hypersonic flow past a cylinder}
In this case, the hypersonic flows impinging on a unit cylinder are
tested to validate robustness of the current scheme. The first one
the the inviscid flow, which were also studied in \cite{GRP-GKS} as
comparison between GRP and GKS. This problem is initialized by the flow
moving towards to a cylinder with different Mach numbers. The
reflective boundary condition is imposed on the surface of cylinder,
and the outflow boundary condition is set on the right boundary. In the
computation, the $60\times100$ mesh points are used, which is shown in
Fig.\ref{cylinder-1}. The Mach number distributions for the flows
with $Ma=5, 10$, and $20$ are also presented in Fig.\ref{cylinder-1},
which show that the current scheme can capture strong shocks very well without
carbuncle phenomenon \cite{Case-Pandolfi}. The robustness of the scheme is well validated.

\begin{figure}[!h]
\centering
\includegraphics[width=0.18\textwidth]{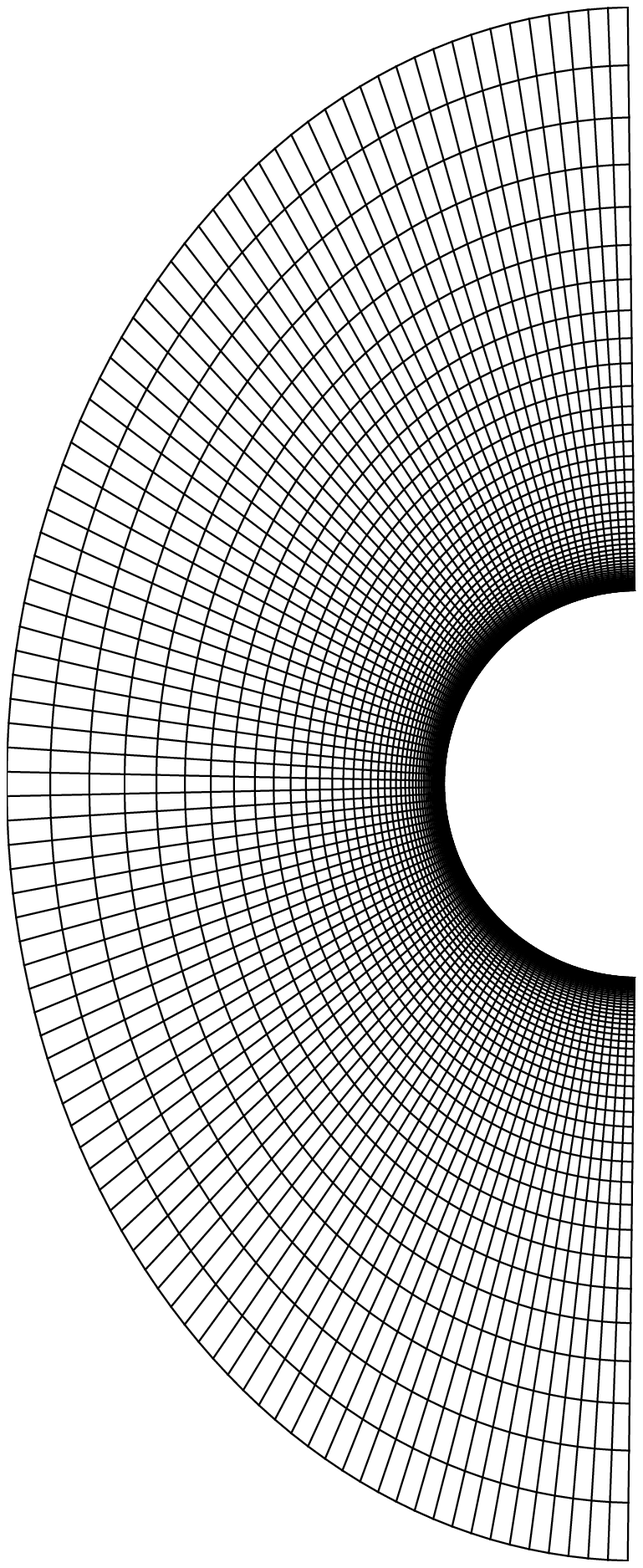}
\includegraphics[width=0.18\textwidth]{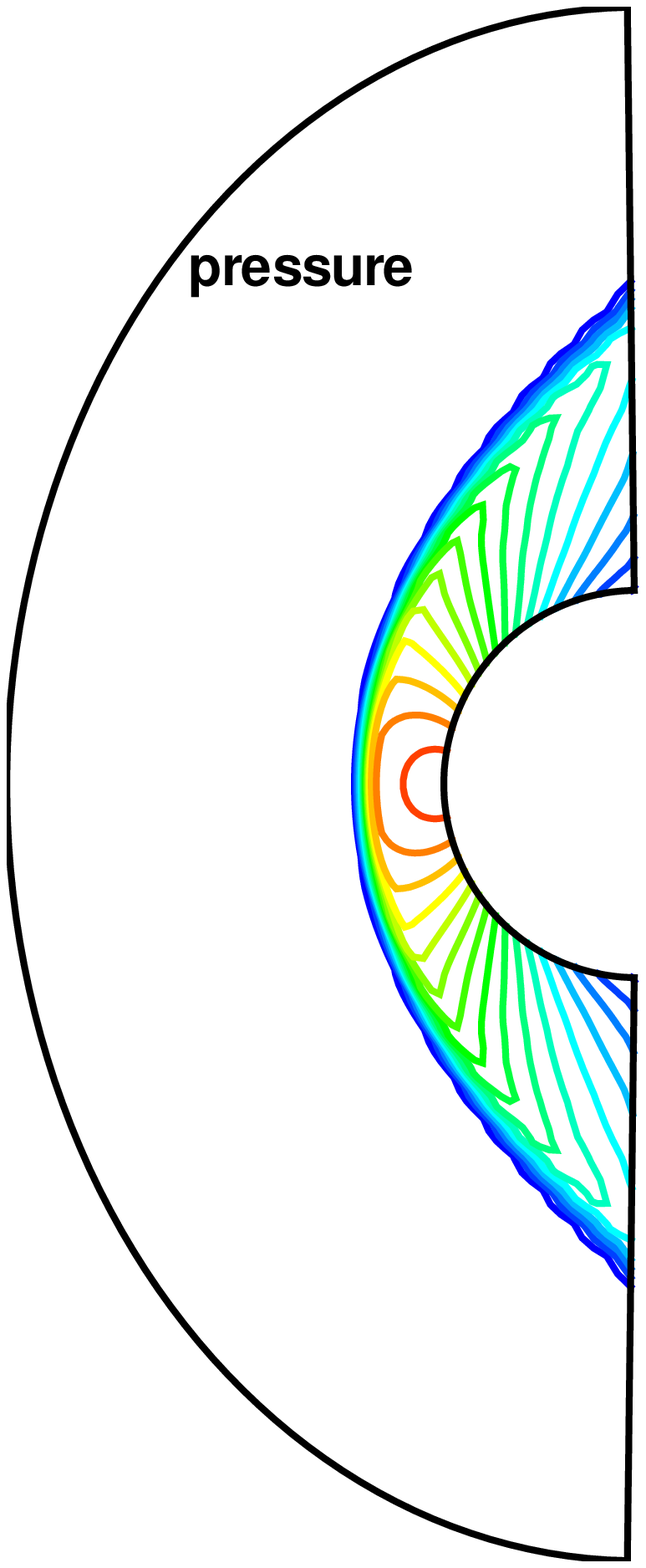}
\includegraphics[width=0.18\textwidth]{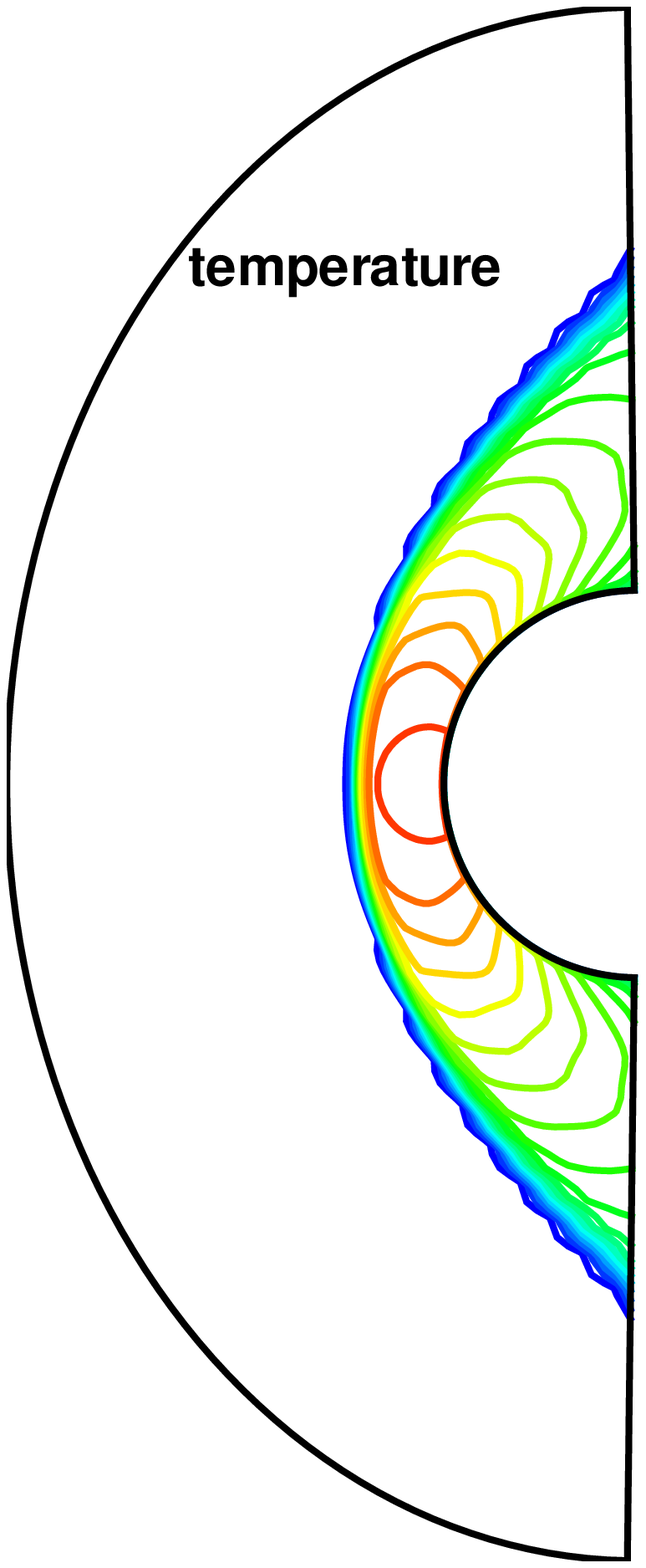}
\includegraphics[width=0.18\textwidth]{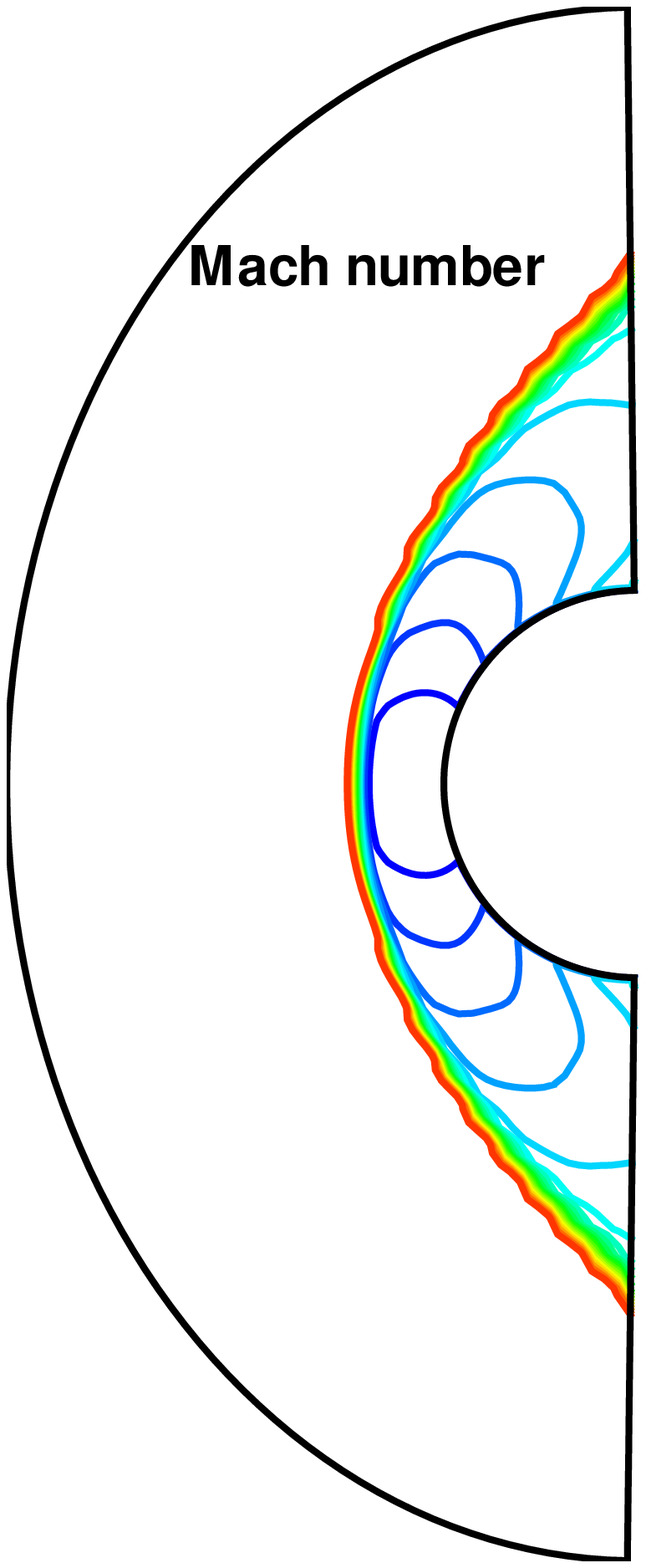}
\caption{\label{cylinder-2}Hypersonic viscous flow past a cylinder
with $Ma=8.03$: the mesh, pressure, temperature, and Mach number
distributions. }
\end{figure}

\begin{figure}[!h]
\centering
\includegraphics[width=0.42\textwidth]{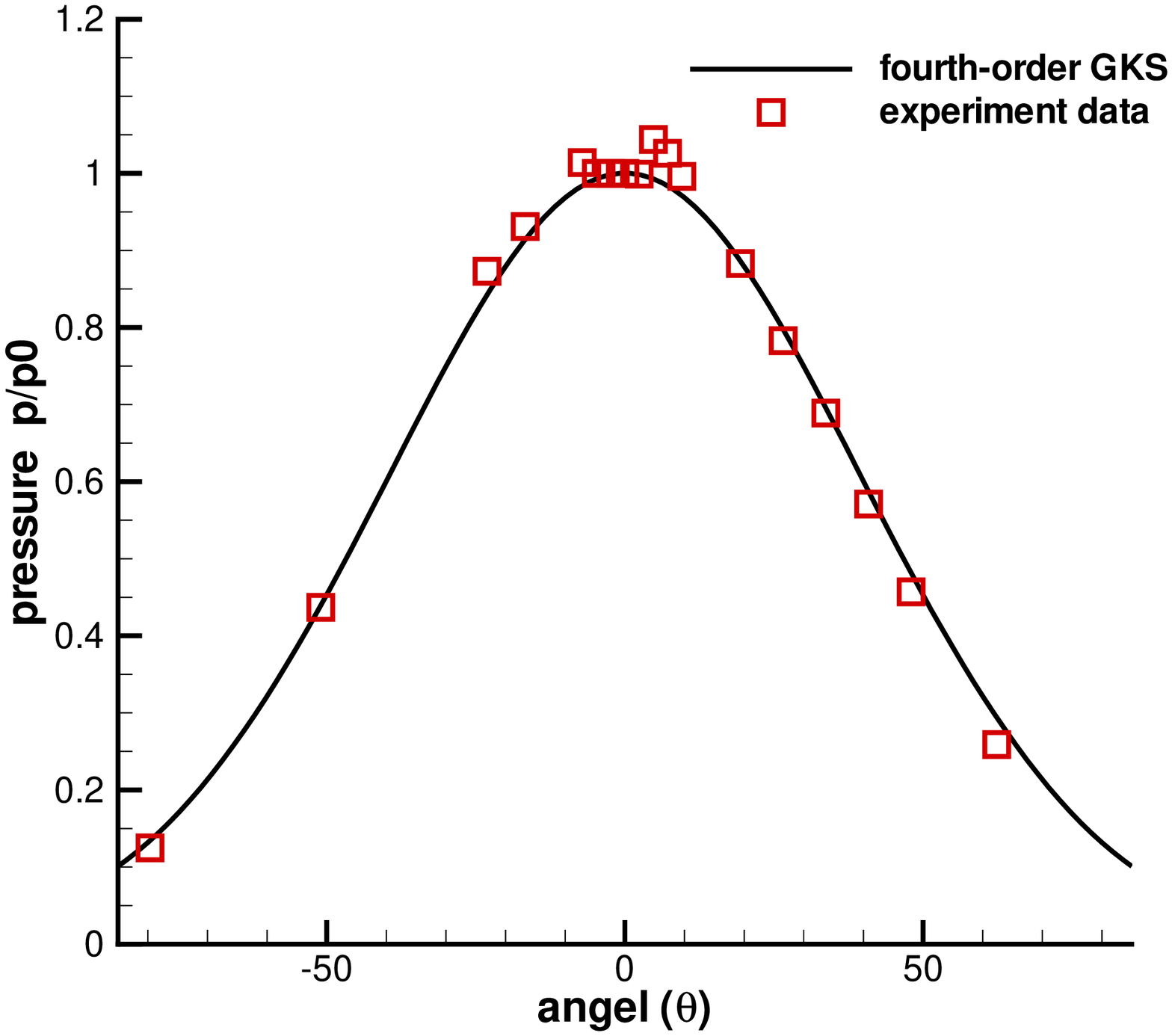}
\includegraphics[width=0.42\textwidth]{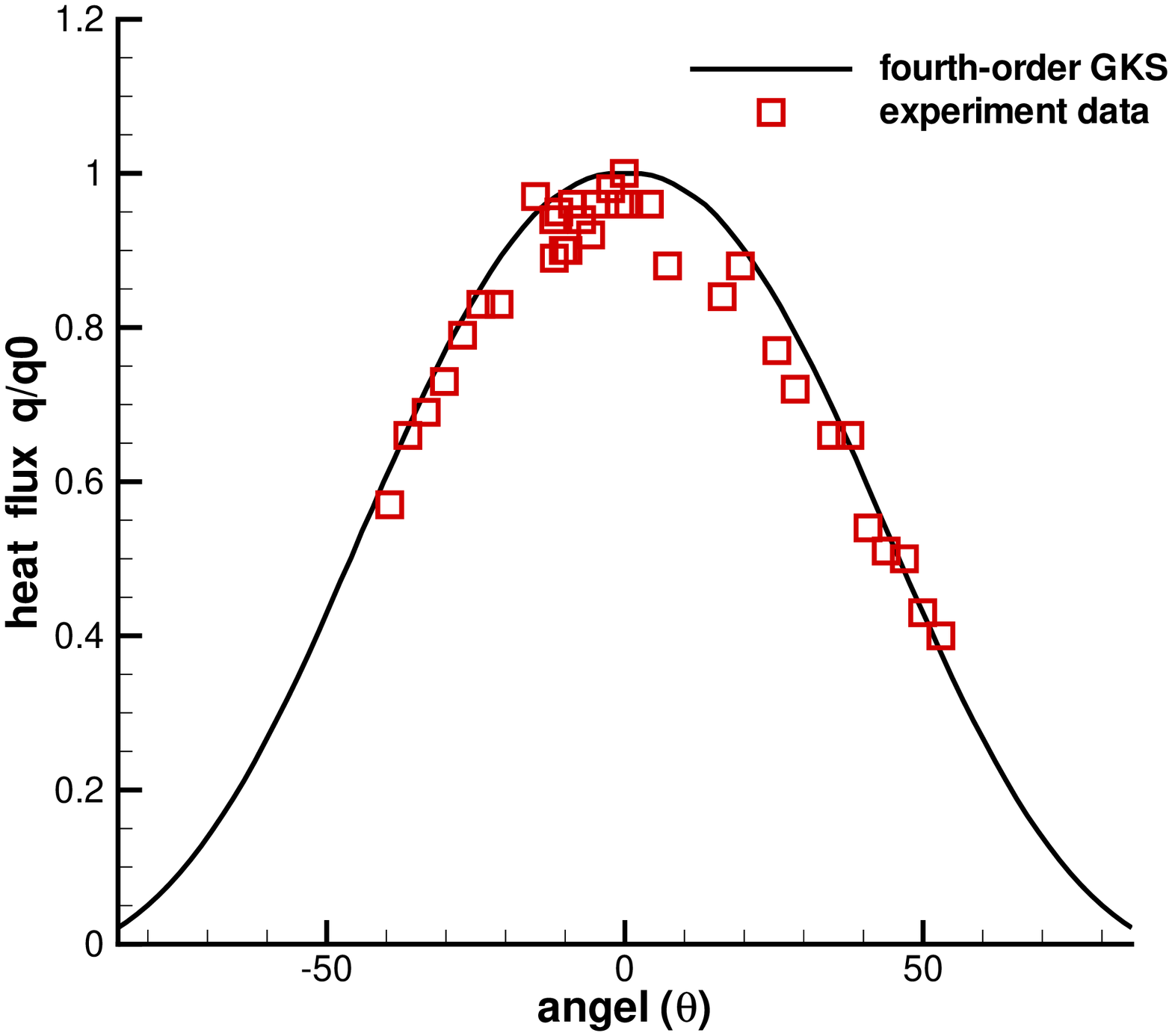}
\caption{\label{cylinder-3}Hypersonic viscous flow past a cylinder
with $Ma=8.03$. Comparison of the computed pressure and heat flux
along the cylindrical surface with the experimental data
\cite{Case-Wieting}. }
\end{figure}

The viscous and heat conducting case at high Mach number is also tested. The flow condition
is given as $Ma_\infty=8.03$, $T_\infty=124.94K$ for the far field,
the wall temperature is $T_W=294.44K$, and the Reynolds number is
$Re=1.835\times10^5$ with cylinder radius and the far field flow
parameters. This test case is taken from the experiment done by
Wieting \cite{Case-Wieting}.  A non-uniform mesh of $60\times160$
cells is used with the near-wall cell width of $1/2000$ to resolve
the boundary layer. The mesh, pressure, temperature, and Mach number
distributions are given in Fig.\ref{cylinder-2}. The pressure and
heat flux along the cylindrical surface are presented in
Fig.\ref{cylinder-3}, where the numerical results agree well with the
experimental data \cite{Case-Wieting}.

\subsection{Laminar boundary layer}
\begin{figure}[!h]
\centering
\includegraphics[width=0.4\textwidth]{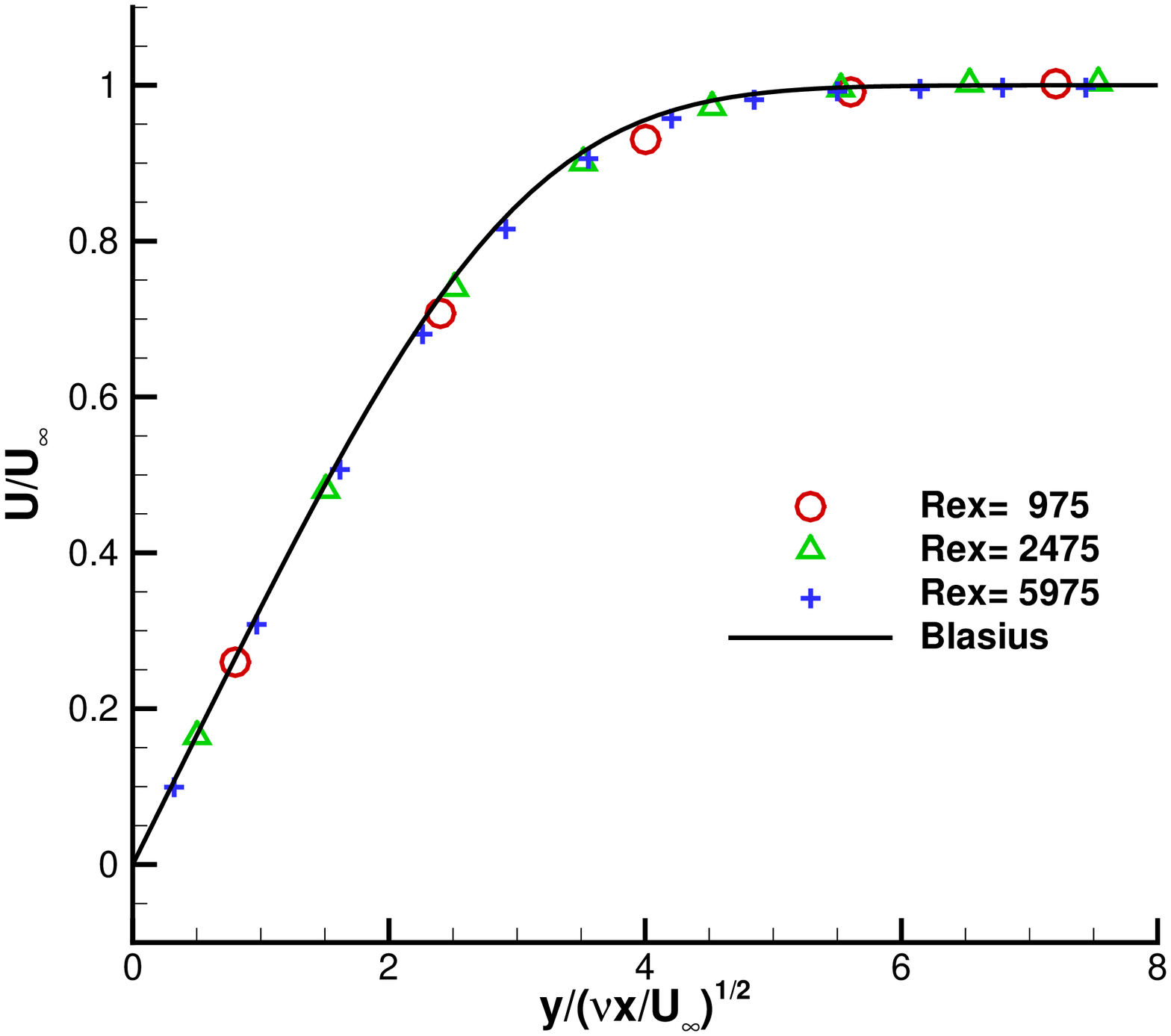}
\includegraphics[width=0.4\textwidth]{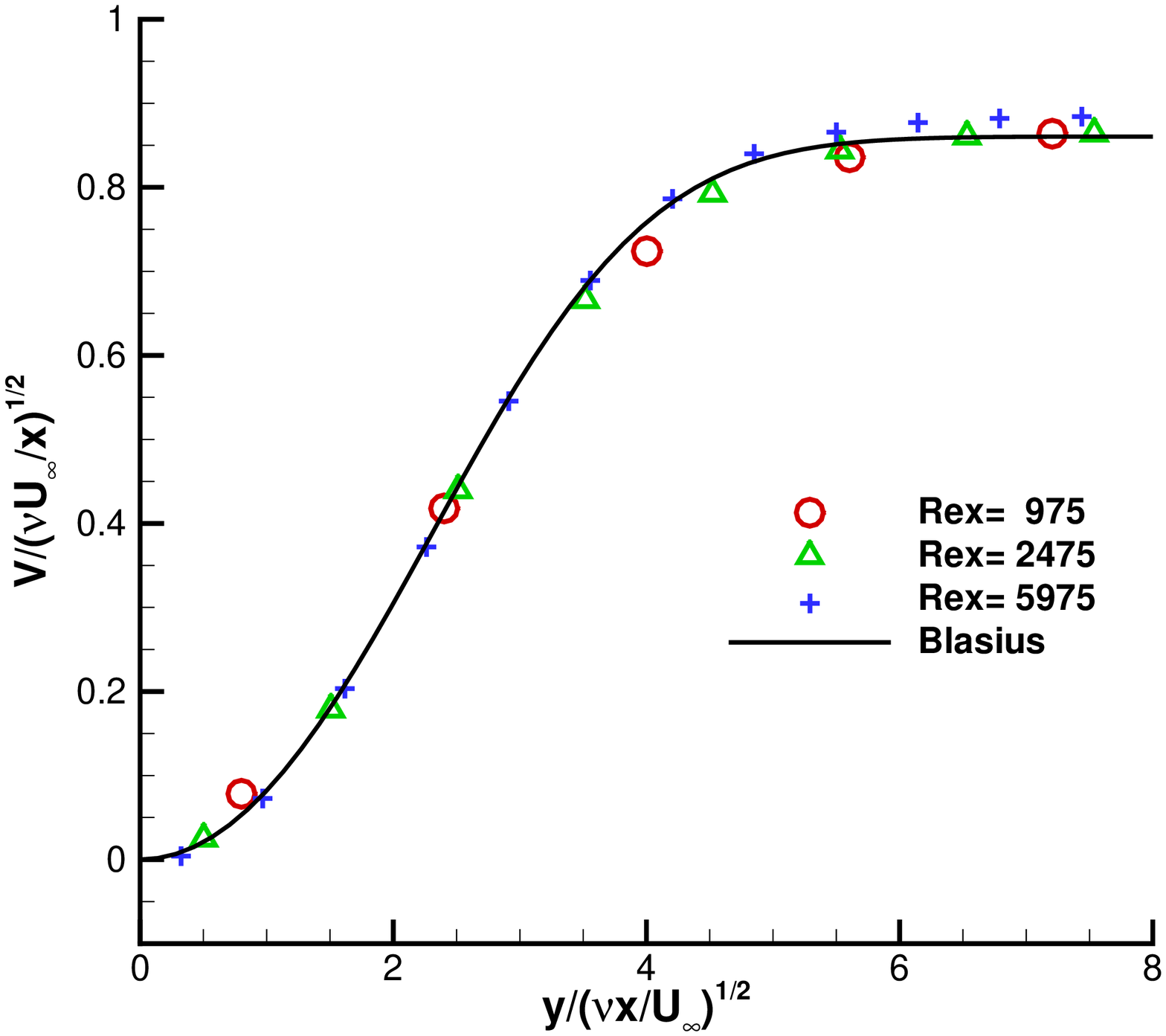}
\includegraphics[width=0.4\textwidth]{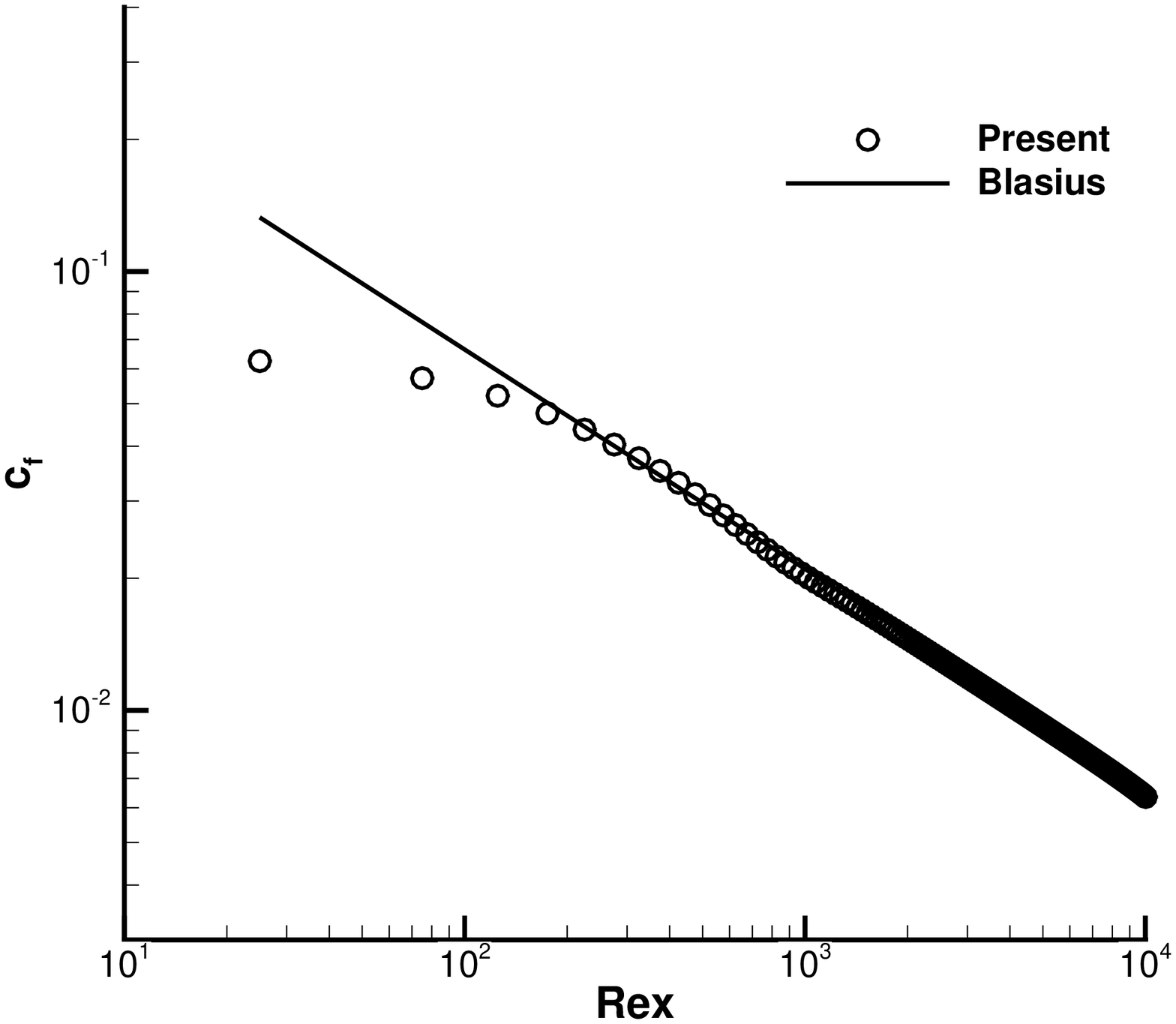}
\caption{\label{boundary}Laminar boundary layer: the $U$ and $V$
velocity profiles at different locations and wall friction
coefficient distribution.}
\end{figure}

A laminar boundary layer is tested over a flat plate. The Mach
number of the free-stream is $Ma=0.15$ and the Reynolds number is
$Re=U_{\infty}L/\nu=10^5$, $\nu$ is the viscous coefficient.  The
non-slip adiabatic boundary condition at the plate is used and a
symmetric condition is imposed at the bottom boundary before the flat
plate. The non-reflecting boundary condition based on the Riemann
invariants is adopted for the other boundaries.
A uniform mesh $260\times90$ points is adopted with $\Delta x =\Delta y
=1/200$, including $60\times90$ mesh points before the plate.
At steady state, the non-dimensional $U$ and $V$ velocity at different locations are presented in Fig.\ref{boundary}, as well as
the wall friction coefficient.  In all locations, the
numerical solutions match with the exact Blasius solution very
well. At the upstream location, the boundary layer profile can be
accurately captured with only four grid points within the layer.

\subsection{Lid-driven cavity flow}
\begin{figure}[!h]
\centering
\includegraphics[width=0.335\textwidth]{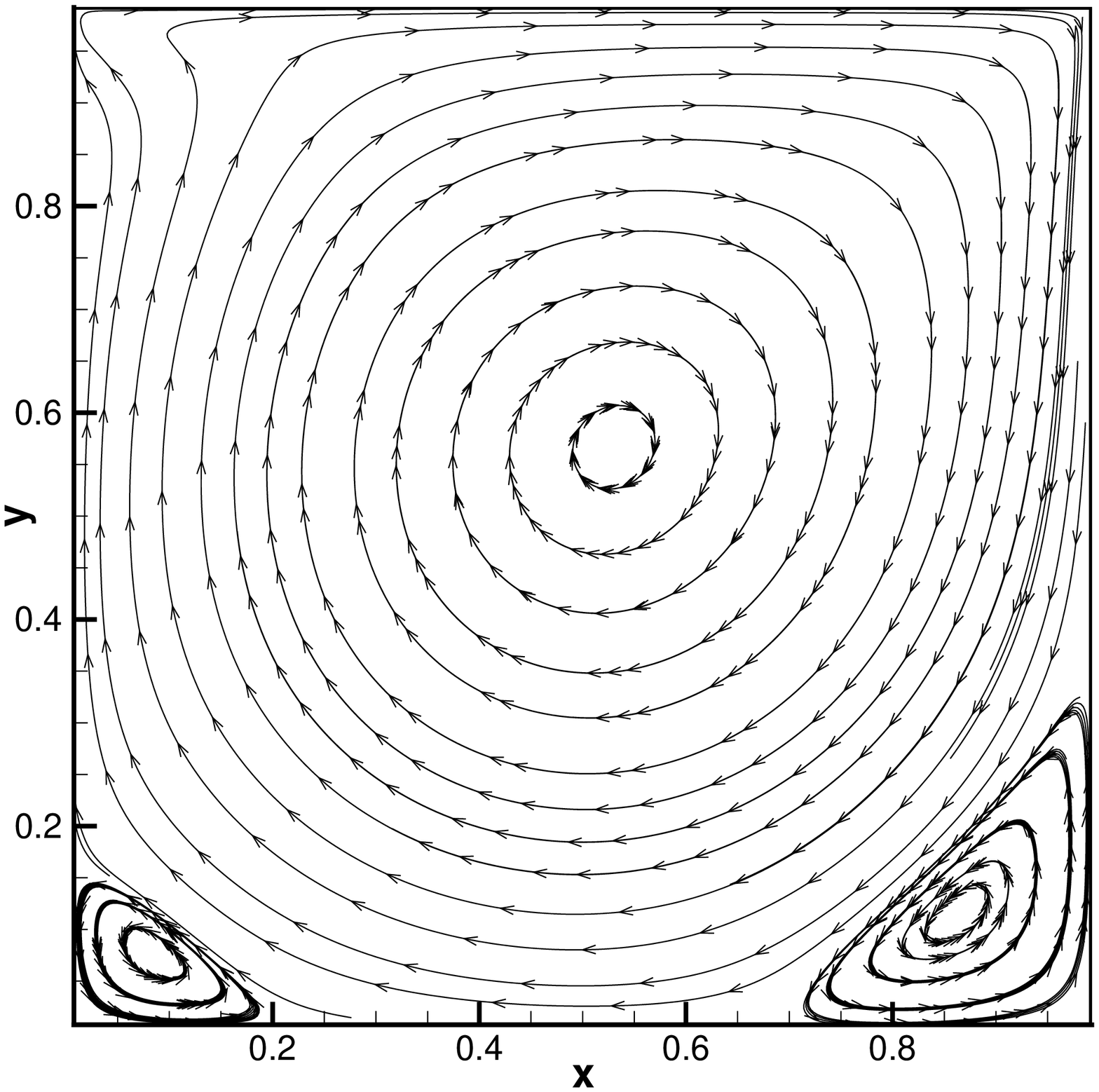}
\includegraphics[width=0.335\textwidth]{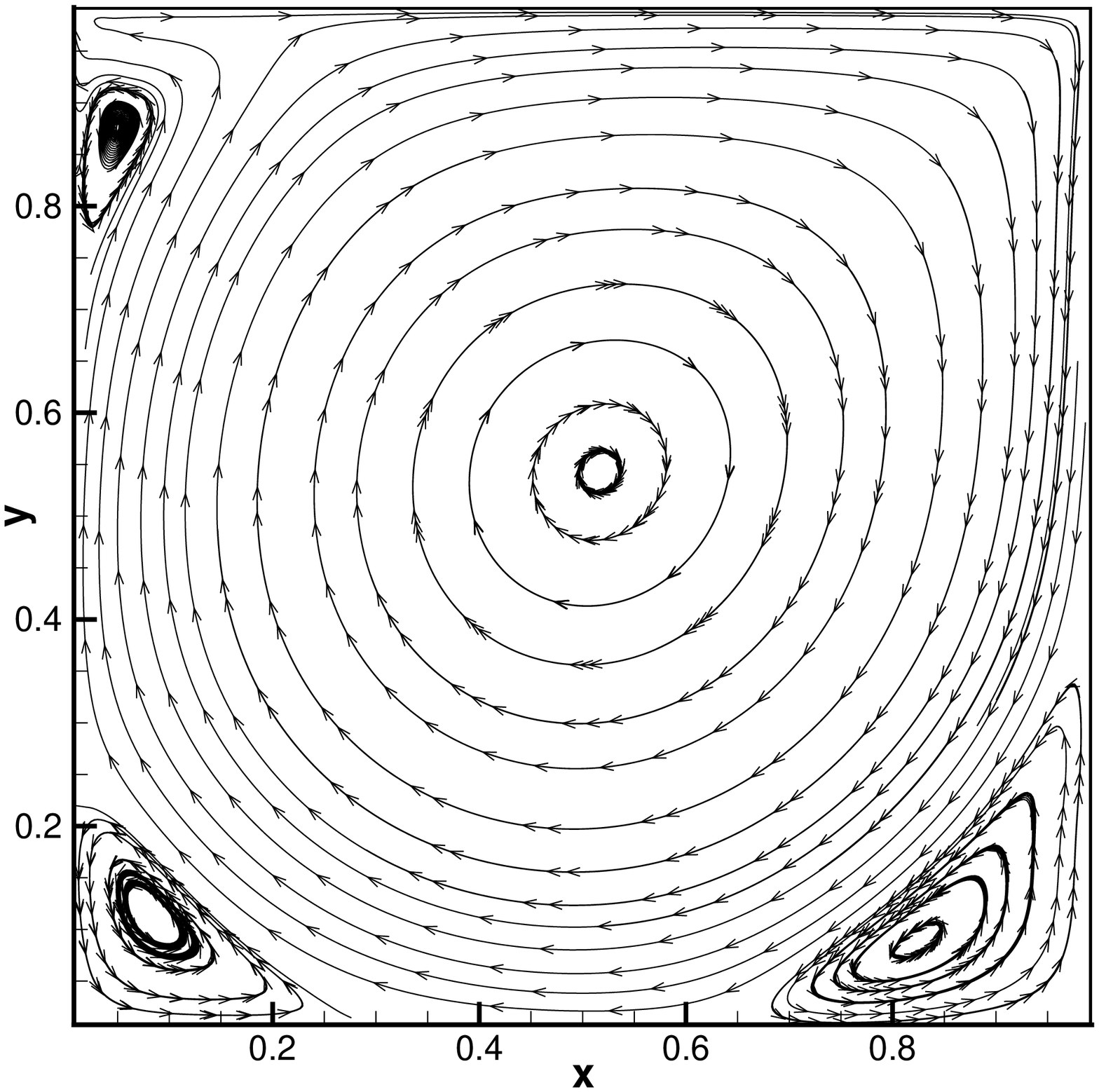}
\caption{\label{cavity-1} Lid-driven cavity flow: the streamlines with $65 \times 65 $ mesh points
from the fourth-order GKS at $Re=1000$ (left) and
$3200$ (right).}
\end{figure}

\begin{figure}[!h]
\centering
\includegraphics[width=0.4\textwidth]{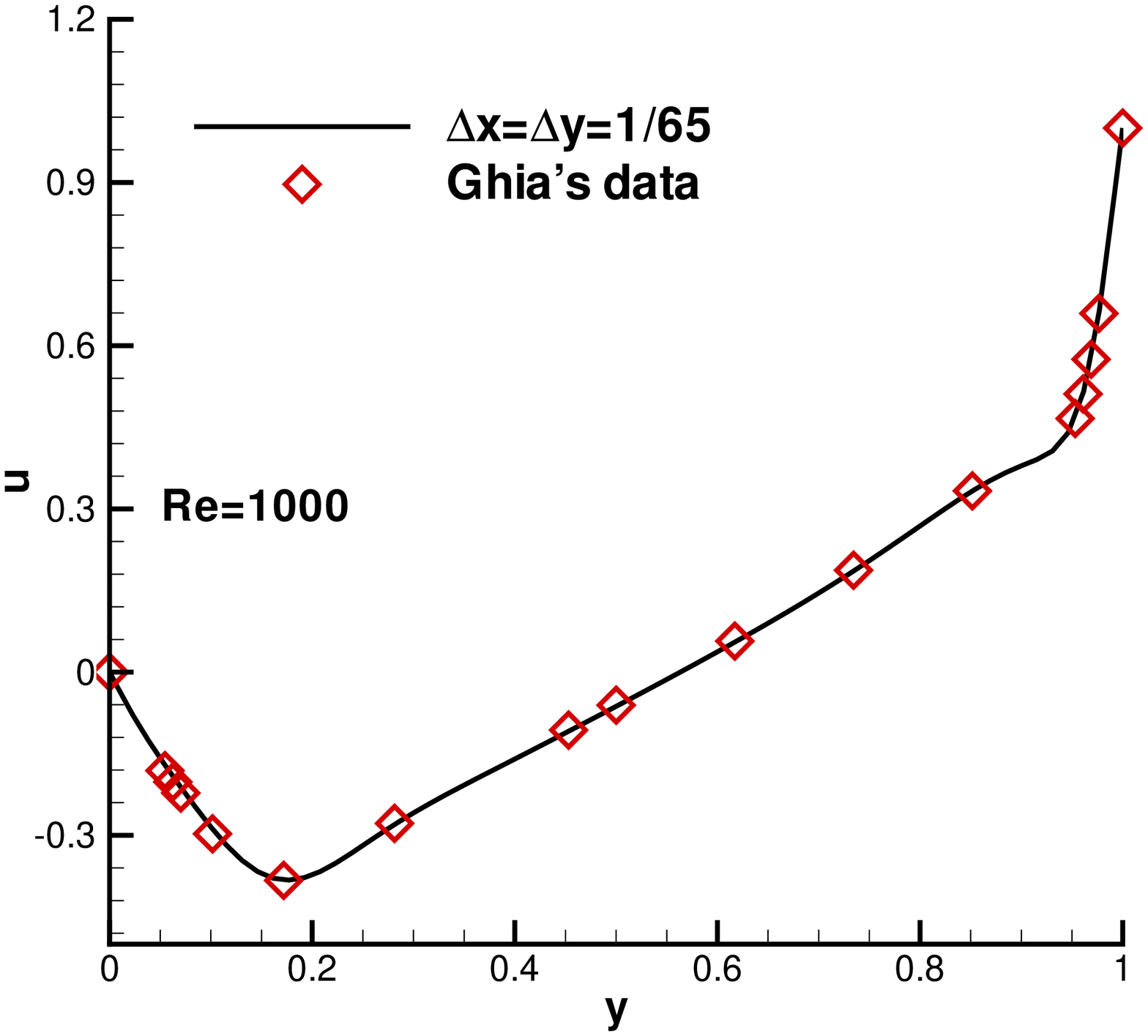}
\includegraphics[width=0.4\textwidth]{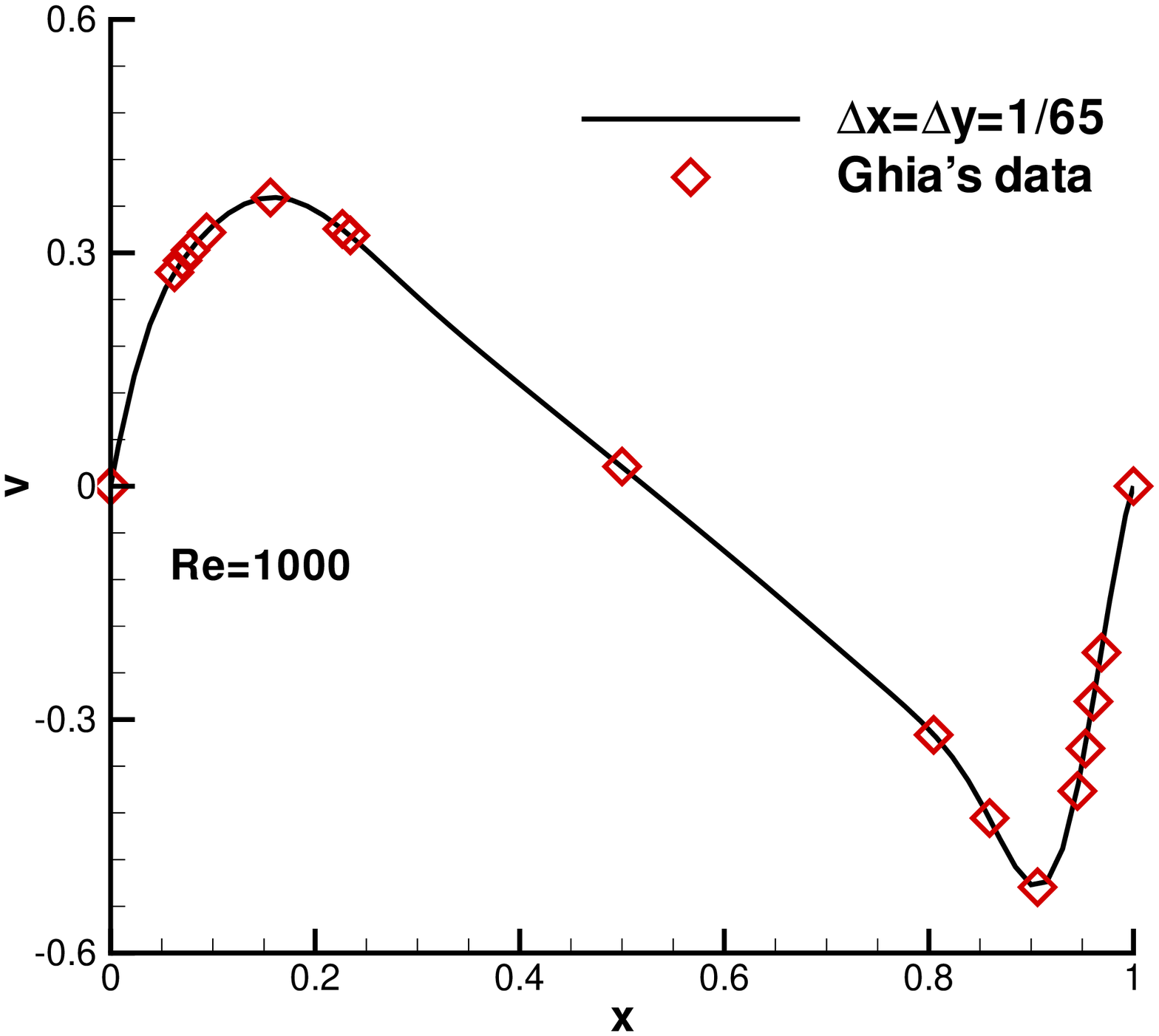}
\includegraphics[width=0.4\textwidth]{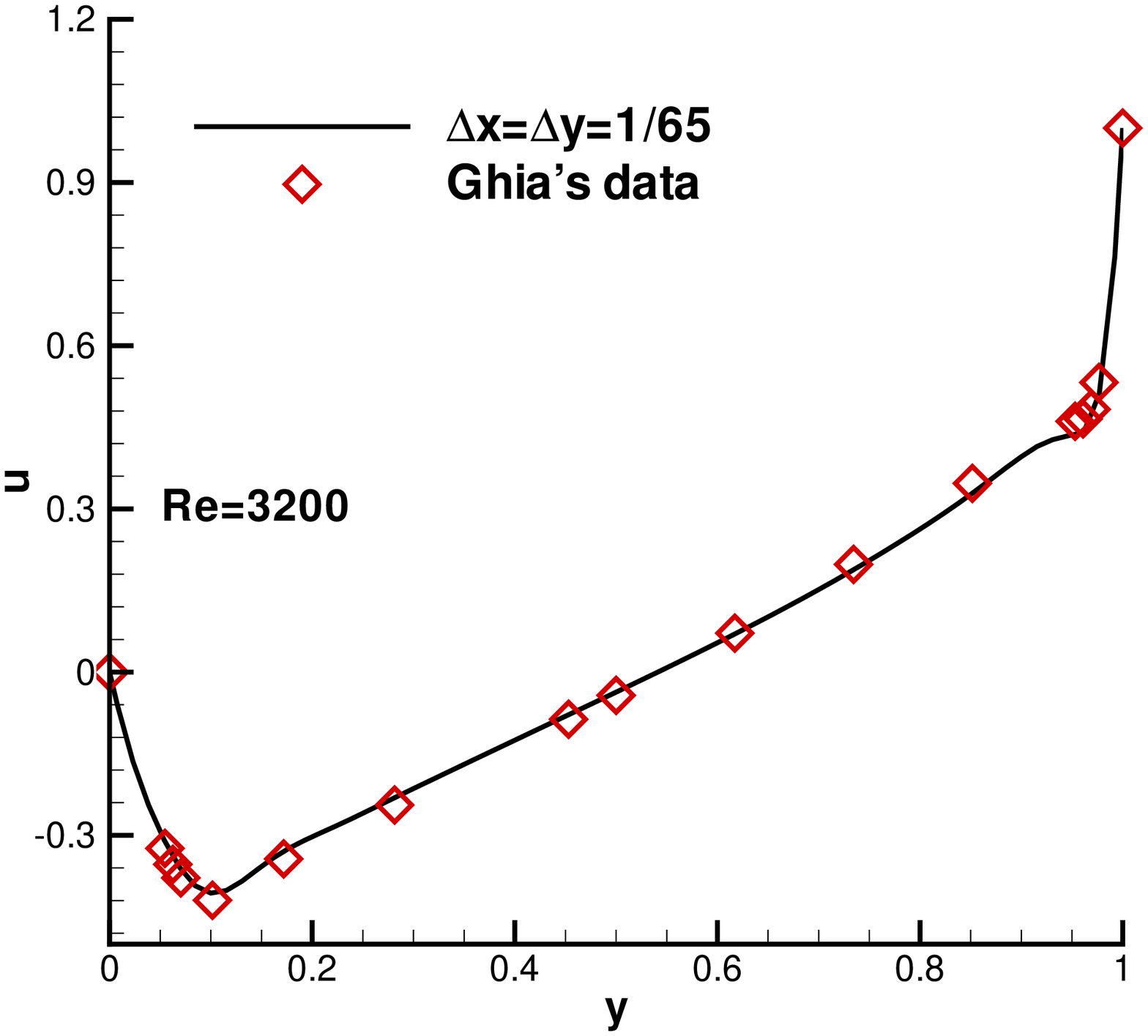}
\includegraphics[width=0.4\textwidth]{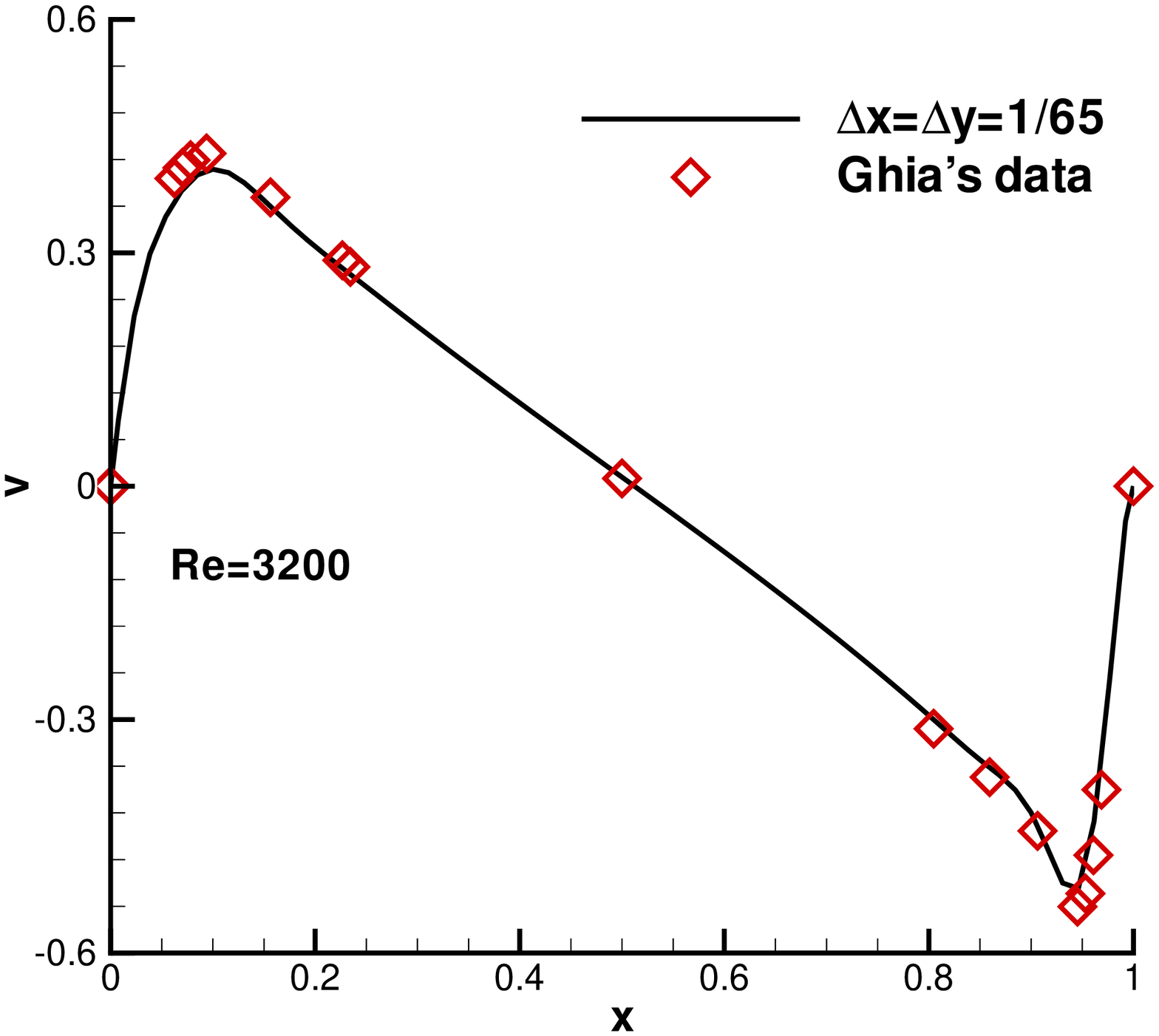}
\caption{\label{cavity-2} Lid-driven cavity flow:  $U$-velocity
along vertical centerline line and $V$-velocity along horizontal
centerline with $65$ mesh points in each direction from the fourth-order GKS at $Re=1000$ and
$3200$.}
\end{figure}

In order to further test the scheme in the capturing of vortex flow,
the lid-driven cavity problem is one of the most important
benchmarks for validating incompressible low speed Navier-Stokes
flow solvers. The fluid is bounded by a unit square and is driven by a
uniform translation of the top boundary. In this case, the flow is
simulated with Mach number $Ma=0.15$ and all boundaries are
isothermal and nonslip. The computational domain $[0, 1]\times[0,
1]$ is covered with $65\times65$ mesh points.
Numerical simulations are conducted for two different Reynolds numbers, i.e.,
$Re=1000$ and $3200$. The streamlines in Fig.\ref{cavity-1}, the
$U$-velocities along the center vertical line, and $V$-velocities along
the center horizontal line, are shown in Fig.\ref{cavity-2}. The
benchmark data \cite{Case-Ghia} for $Re=1000 $ and $ 3200$ are also presented,
and  the simulation results match well with these benchmark data.
The higher-order accuracy of the scheme is clearly demonstrated from the Reynolds number $3200$ case, where
only $65$ uniform mesh points are used in each direction.

\subsection{Viscous shock tube problems}
This problem was introduced to test the performances of different
schemes for viscous flows \cite{Case-Daru}. In this case, an ideal
gas is at rest in a two-dimensional unit box $[0,1]\times[0,1]$. A
membrane located at $x=0.5$ separates two different states of the
gas and the dimensionless initial states are
\begin{equation*}
(\rho,U,p)=\left\{\begin{aligned}
&(120, 0, 120/\gamma), \ \ \ &  0<x<0.5,\\
&(1.2, 0, 1.2/\gamma),  & 0.5<x<1,
\end{aligned} \right.
\end{equation*}
where $\gamma=1.4$ and Prandtl number $Pr=0.73$.

\begin{figure}[!h]
\centering
\includegraphics[width=0.535\textwidth]{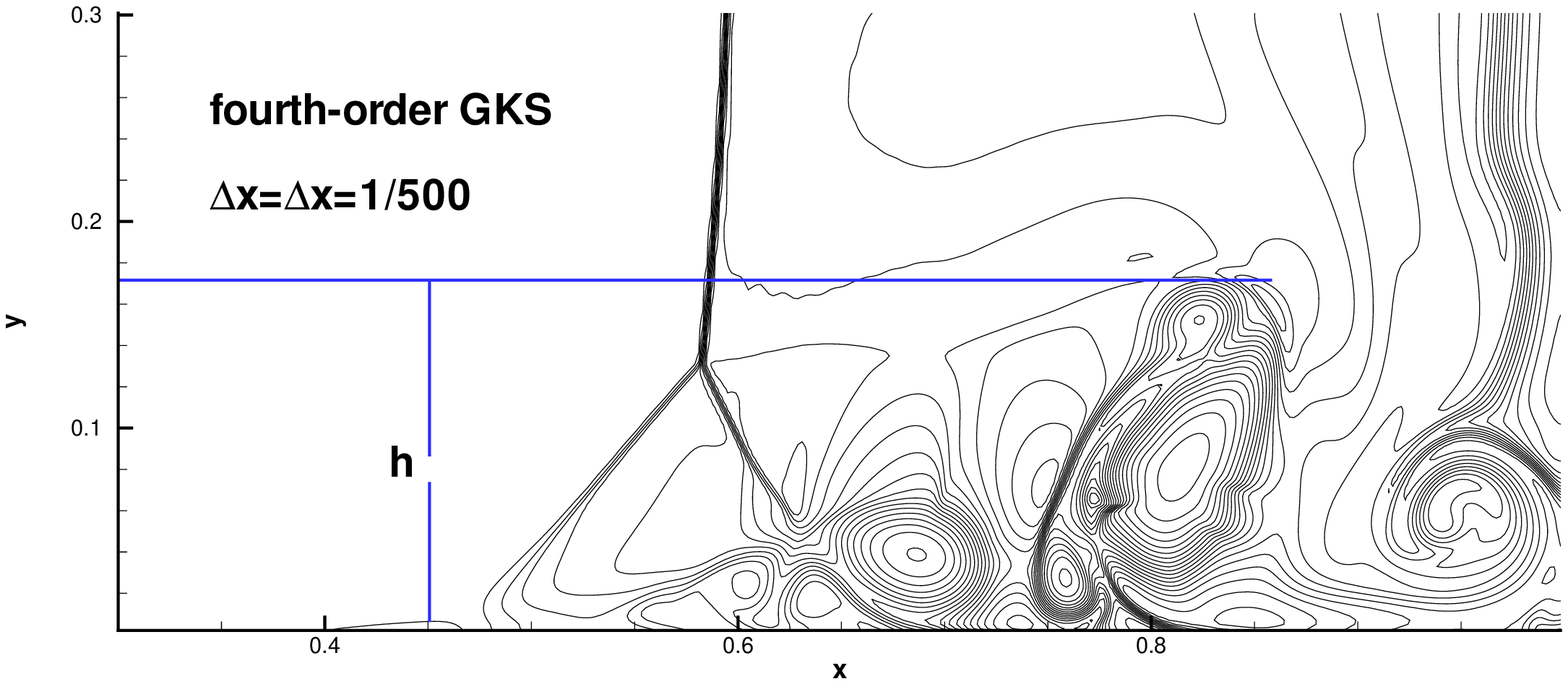}
\includegraphics[width=0.535\textwidth]{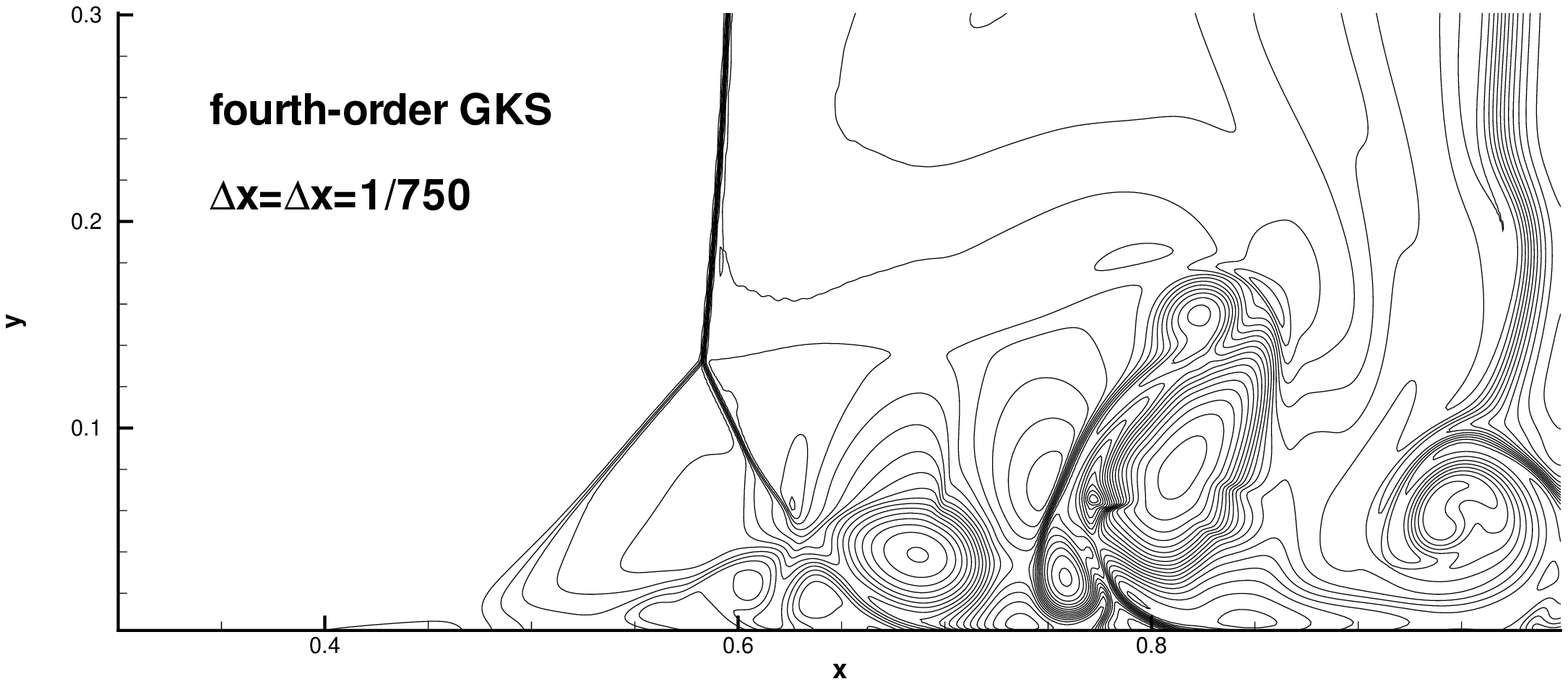}
\caption{\label{shock-boundary-200} Reflecting shock-boundary layer
interaction: density distributions at $t=1$ with $Re=200$ from the
fourth-order GKS with $\Delta x=\Delta y=1/500, 1/750$.}
\end{figure}

\begin{figure}[!h]
\centering
\includegraphics[width=0.55\textwidth]{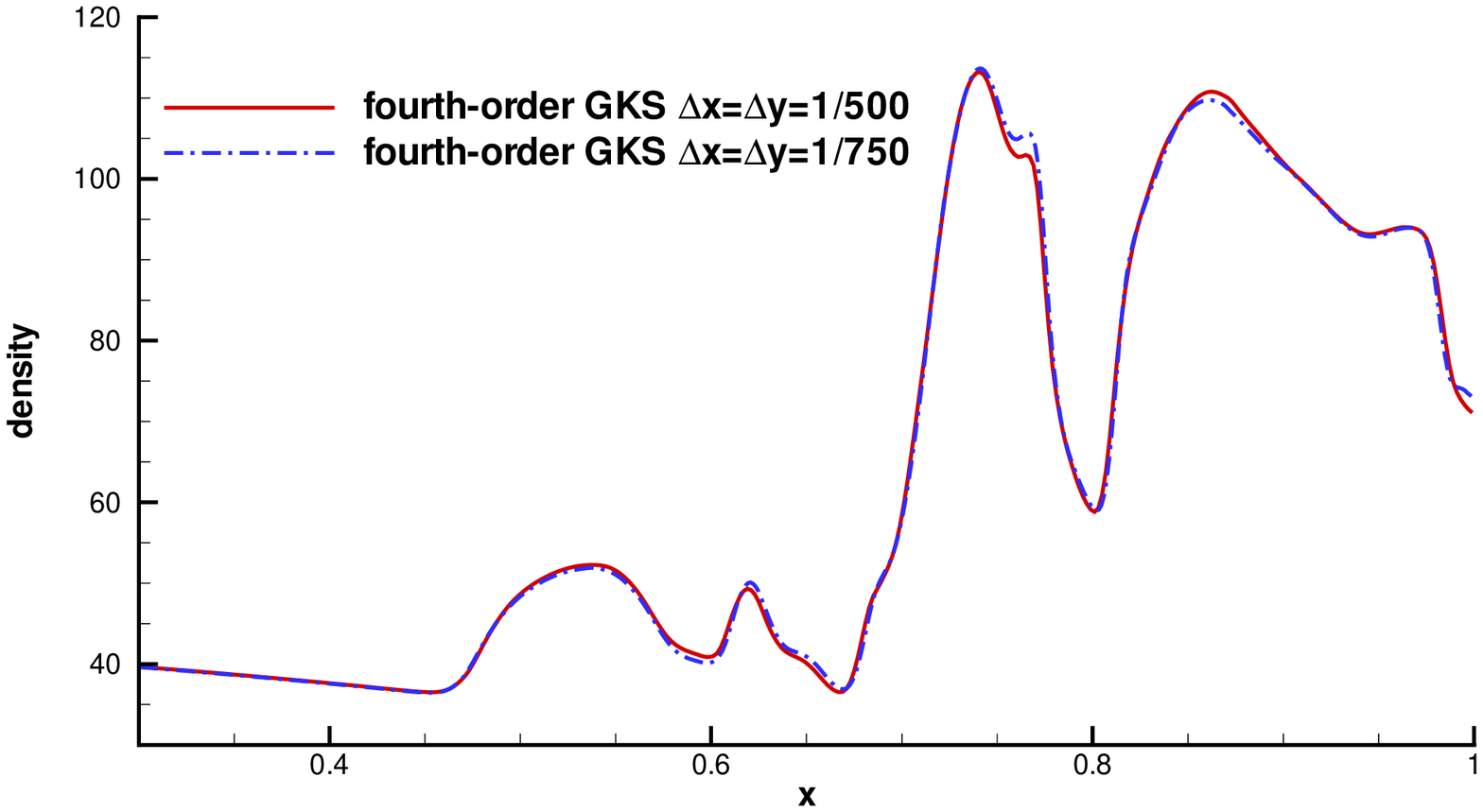}
\caption{\label{shock-boundary-200b} Reflecting shock-boundary layer
interaction: density distribution along the lower wall with
different mesh sizes for $Re=200$.}
\end{figure}

\begin{table}[!h]
\begin{center}
\def\temptablewidth{1\textwidth}
{\rule{\temptablewidth}{0.1pt}}
\begin{tabular*}{\temptablewidth}{@{\extracolsep{\fill}}ccccc}
Scheme  & AUSMPW+  &  M-AUSMPW+  &  fourth-order GKS\\
\hline height & 0.163& 0.168&  0.171
\end{tabular*}
{\rule{\temptablewidth}{0.1pt}}\caption{\label{height} Comparison of
the heights of primary vortex among gas kinetic scheme and
other reference methods \cite{Case-Kim} for the reflecting shock-boundary
layer interaction with $\Delta x=\Delta y=1/500$.}
\end{center}
\end{table}

\begin{figure}[!h]
\centering
\includegraphics[width=0.55\textwidth]{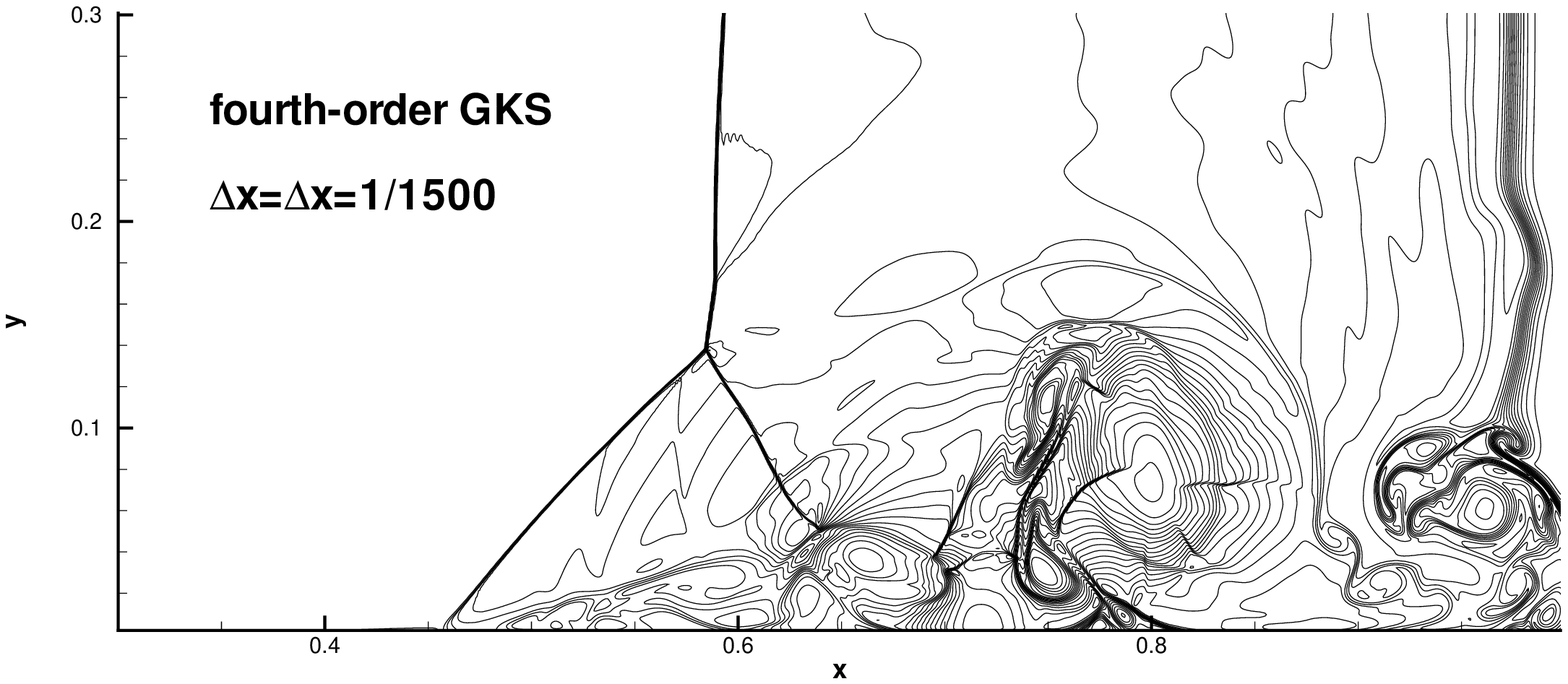}
\includegraphics[width=0.55\textwidth]{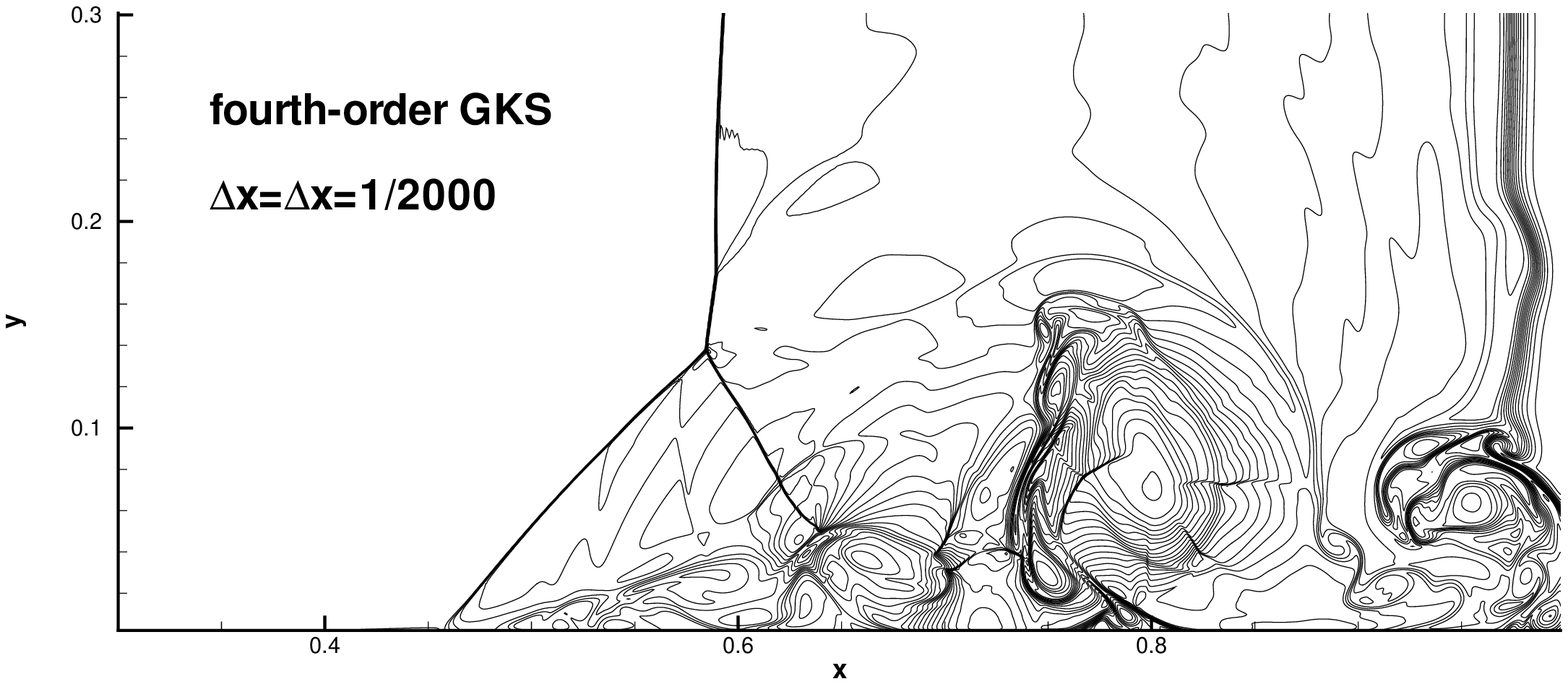}
\caption{\label{shock-boundary-1000} Reflecting shock-boundary layer
interaction. The density distribution at $t=1$ with $Re=1000$ with
$\Delta x=\Delta y= 1/1500$ and $1/2000$.} \centering
\includegraphics[width=0.55\textwidth]{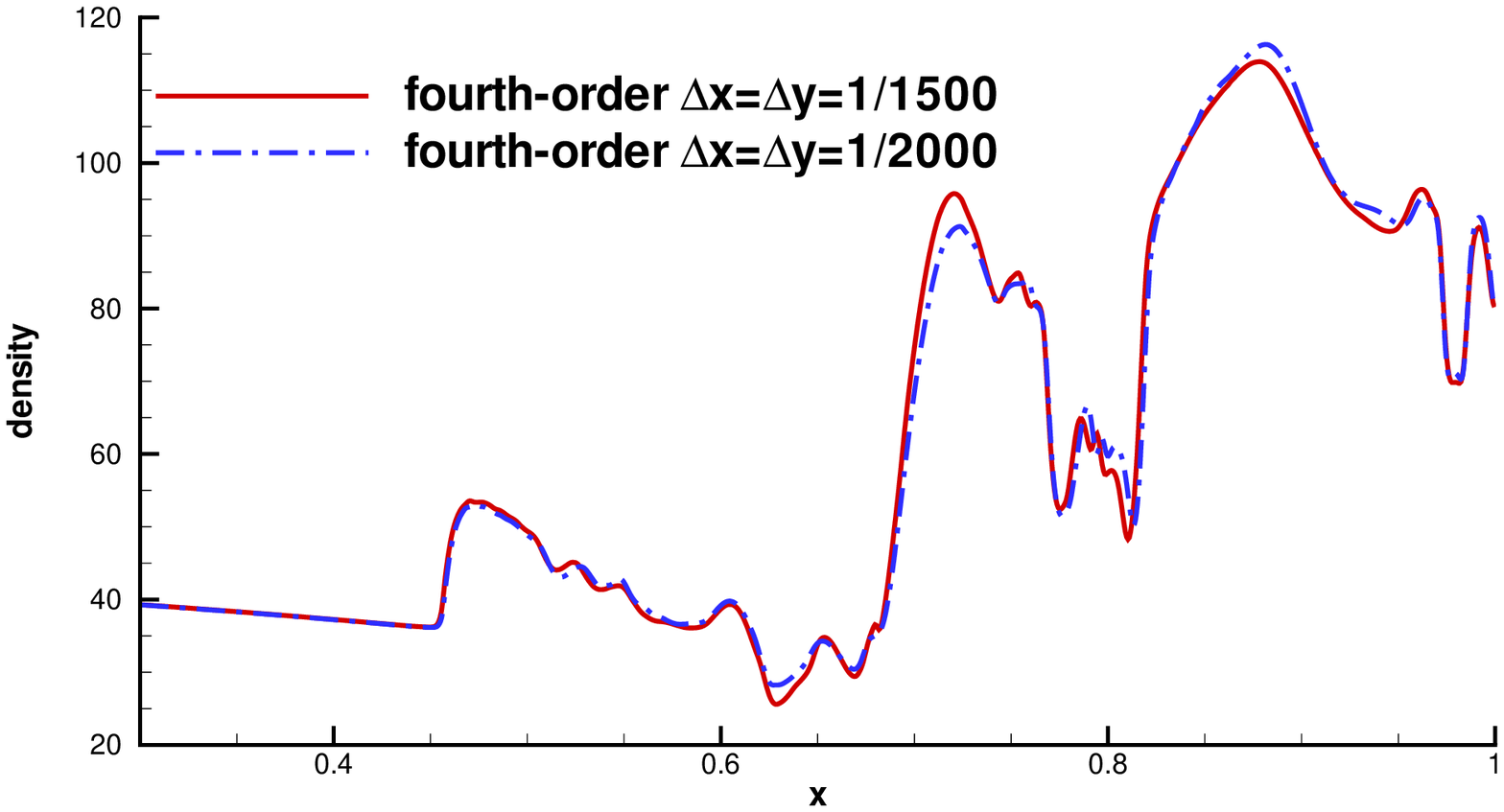}
\caption{\label{shock-boundary-200b} Reflecting shock-boundary layer
interaction: density distribution along the lower wall with
different mesh sizes for $Re=1000$.}
\end{figure}

The membrane is removed at time zero and wave interaction occurs. A
shock wave, followed by a contact discontinuity, moves to the right
with Mach number $Ma=2.37$ and reflects at the right end wall. After
the reflection, it interacts with the contact discontinuity. The
contact discontinuity and shock wave interact with the horizontal
wall and create a thin boundary layer during their propagation. The
solution will develop complex two-dimensional
shock/shear/boundary-layer interactions. This case is tested in the
computational domain $[0, 1]\times[0, 0.5]$, a symmetric boundary condition
is used on the top boundary $x\in[0, 1], y=0.5$, and non-slip
boundary condition, and adiabatic condition for temperature are
imposed at solid wall boundaries. The case with $Re=200$ is tested first.
The density distributions  are presented in Fig.\ref{shock-boundary-200} with two different mesh resolutions.
The results match well with each other. The density profiles along
the lower wall with $Re=200$ are also presented in
Fig.\ref{shock-boundary-200b}. A mesh-convergent solution
is observed for $Re=200$. As shown in Table.\ref{height}, the height
of primary vortex predicted by the current scheme agrees well with
the reference data \cite{Case-Kim}.
For the case with $Re=1000$, the
flow structure becomes more complicated. The density distributions
from the current scheme are given in Fig.\ref{shock-boundary-1000},
and the density profiles along the lower wall are presented in
Fig.\ref{shock-boundary-200b} with the mesh size $\Delta x=\Delta
y=1/1500$ and $1/2000$. The flow structure is complicated, and the
mesh convergence is basically obtained with the mesh size decreasing to
$1/1500$. The current results agree well with the
reference data very well. More studies for this problem can be found in
\cite{Case-Daru}.

\section{Conclusion}
In this paper, based on the two-stage time stepping method  a fourth-order gas-kinetic scheme is proposed
for both inviscid and viscous flow computations. With the fifth-order WENO
reconstruction, a GKS with a fifth order accuracy in space and a
fourth order accuracy in time is developed.
In comparison with the classical methods based on the first-order Riemann solver,
for a fourth-order accuracy in time the current GKS  only uses two stages instead of four stages in the conventional methods.
Therefore, the current GKS should be much more efficient than  these higher-order methods based on the Riemann solutions, especially for the NS solutions.
The current finite volume scheme can use a CFL number on the order of $0.5$.
The further development of the GKS to even higher-order accuracy can be achieved
with the inclusion of the second-order time derivative of the flux function, such as the fifth-order scheme presented in the Appendix.
The fourth-order GKS not only has the expected order of accuracy for the
smooth flow, but also has favorable shock capturing property for the
discontinuous solutions.
Most importantly, the numerical tests clearly demonstrate that the current fourth-order scheme is as robust as the second-order one.

By taking advantage of the time-accurate gas evolution model in the flux evaluation, an efficient and accurate fourth-order gas-kinetic scheme has been constructed.
The advantage of the time accurate evolution model can be further explored
for the construction of higher-order compact schemes \cite{GKS-high3,GKS-high4}, where not only the flux at the cell
interface, the conservative flow variables at the cell interface can be updated as well to construct compact stencils for the flow reconstruction in the next time level.
Based on the current study, we can conclude that the adaptation of  a higher-order gas evolution model for the flux evaluation has an indispensable advantage in the development of
 higher-order schemes. The real bottleneck which hinders the progress for the development of higher-order accurate and robust schemes for the compressible
flows is due to the use of first-order Riemann solver.

\section*{Acknowledgements}
The research of K.  Xu is supported by Hong Kong Research Grant Council (620813, 16211014, 16207715) and HKUST research fund
(PROVOST13SC01, IRS15SC29, SBI14SC11).
The research of Q. Li is partially supported by NSFC (11172154).  The work of J. Li is supported by NSFC (91130021, 11371063), and the doctoral program from the Education Ministry of China (20130003110004).

\section*{Appendix: Extension to higher order}

The key point for  developing a two-stage fourth-order temporal accurate
schemes  is the use of a time-dependent flux function.
The third-order GRP and GKS have
both first- and second-order time derivatives in the flux function \cite{Qian-Li, GKS-high1,GKS-high2,GKS-high3}.
Thus, with a two-stage temporal
discretization and the third-order GRP and GKS flux solvers, it is possible to
develop  schemes with fifth-order accuracy in time.
\vspace{0.2cm}

We consider the time-dependent equation Eq.\eqref{pde} with the
initial condition Eq.\eqref{pde2}.  Introducing an intermediate state
at $t_*=t_n+A\Delta t$,
\begin{align}\label{step-1}
\textbf{w}^*=\textbf{w}^n+A\Delta t\mathcal
{L}(\textbf{w}^n)+\frac{1}{2}A^2\Delta t^2\frac{\partial}{\partial
t}\mathcal{L}(\textbf{w}^n)+\frac{1}{6}A^3\Delta
t^3\frac{\partial^2}{\partial t^2}\mathcal{L}(\textbf{w}^n),
\end{align}
we will have  $\mathcal {L}(\textbf{u}^*)$,
$\displaystyle\frac{\partial }{\partial t}\mathcal
{L}(\textbf{w}^*)$ and $\displaystyle\frac{\partial^2 }{\partial
t^2}\mathcal {L}(\textbf{w}^*)$ at $t_*$.
Then, the update scheme can be written as
\begin{align}\label{step-2}
\textbf{w}^{n+1}=\textbf{w}^n&+\Delta t(B_0\mathcal
{L}(\textbf{w}^n)+B_1\mathcal {L}(\textbf{w}^*))+\frac{1}{2}\Delta
t^2\big(C_0\frac{\partial}{\partial
t}\mathcal{L}(\textbf{w}^n)+C_1\frac{\partial}{\partial
t}\mathcal{L}(\textbf{w}^*)\big)\nonumber\\&+\frac{1}{6}\Delta
t^3\big(D_0\frac{\partial^2}{\partial
t^2}\mathcal{L}(\textbf{w}^n)+D_1\frac{\partial^2}{\partial
t^2}\mathcal{L}(\textbf{w}^*)\big).
\end{align}
It can be proved that  Eq.\eqref{step-1} and
Eq.\eqref{step-2} can provide a fifth-order temporal accurate approximation to
the solution $\textbf{w}(t)$ at $t=t_n +\Delta t$ with the following
coefficients
\begin{align}\label{sol}
A=\frac{2}{5}, B_0=1, B_1=0, C_0=1, C_1=0, D_0=\frac{3}{8},
D_1=\frac{5}{8}.
\end{align}

To prove this proposition, the following equation needs to be
satisfied, using the same approach as in \cite{GRP-high},
\begin{align*}
\textbf{w}^{n+1}=\textbf{w}^n+\int_{t_n}^{t_n+\Delta
t}\mathcal{L}(\textbf{w}(t))dt+\mathcal{O}(\Delta t^6).
\end{align*}
According to the Taylor expansion of the operator $\mathcal{L}$ at
$t_n$, the integral can be expressed as
\begin{align}\label{exp}
\int_{t_n}^{t_n+\Delta t}\mathcal{L}(\textbf{w}(t))dt=\Delta
t\mathcal{L}+\frac{\Delta t^2}{2}\frac{\partial\mathcal{L}}{\partial
t} +\frac{\Delta t^3}{6}\frac{\partial^2\mathcal{L}}{\partial
t^2}+\frac{\Delta t^4}{24}\frac{\partial^3\mathcal{L}}{\partial
t^3}+\frac{\Delta t^5}{120}\frac{\partial^4\mathcal{L}}{\partial
t^4}+\mathcal{O}(\Delta t^6),
\end{align}
where the time derivatives for the operator $\mathcal{L}$ can be
given by the chain rule, for example
\begin{align*}
\displaystyle\frac{\partial\mathcal{L}}{\partial
t}=\mathcal{L}_{\textbf{w}}\mathcal{L},~~~
\frac{\partial^2\mathcal{L}}{\partial
t^2}=\mathcal{L}_\textbf{w}^2\mathcal{L}+\mathcal{L}_{\textbf{ww}}\mathcal{L}^2,
...
\end{align*}
Denote
$\displaystyle\mathcal{G}(\textbf{w})=\frac{\partial}{\partial
t}\mathcal{L}(\textbf{w})$ and
$\displaystyle\mathcal{H}(\textbf{w})=\frac{\partial^2}{\partial
t^2}\mathcal{L}(\textbf{w})$, and expand
$\displaystyle\mathcal{L}(\textbf{w})$,
$\displaystyle\mathcal{G}(\textbf{w})$ and
$\displaystyle\mathcal{H}(\textbf{w})$ in the neighboring of
$\textbf{u}^*$ to the corresponding order, we have
\begin{align*}
\mathcal{L}(\textbf{w}^*)&=\mathcal{L}(\textbf{w}^n)+\mathcal{L}_\textbf{w}(\textbf{w}^*-\textbf{w}^n)+\frac{\mathcal{L}_{\textbf{w}\textbf{w}}}{2}(\textbf{w}^*-\textbf{w}^n)^2
+\frac{\mathcal{L}_{\textbf{w}\textbf{w}\textbf{w}}}{6}(\textbf{w}^*-\textbf{w}^n)^3+\frac{\mathcal{L}_{\textbf{w}\textbf{w}\textbf{w}\textbf{w}}}{24}(\textbf{w}^*-\textbf{w}^n)^4,\\
\mathcal{G}(\textbf{w}^*)&=\mathcal{G}(\textbf{w}^n)+\mathcal{G}_\textbf{w}(\textbf{w}^*-\textbf{w}^n)+\frac{\mathcal{G}_{\textbf{w}\textbf{w}}}{2}(\textbf{w}^*-\textbf{w}^n)^2
+\frac{\mathcal{G}_{\textbf{w}\textbf{w}\textbf{w}}}{6}(\textbf{w}^*-\textbf{w}^n)^3,\\
\mathcal{H}(\textbf{w}^*)&=\mathcal{H}(\textbf{w}^n)+\mathcal{H}_\textbf{w}(\textbf{w}^*-\textbf{w}^n)+\frac{\mathcal{H}_{\textbf{w}\textbf{w}\textbf{w}}}{2}(\textbf{w}^*-\textbf{w}^n)^2,
\end{align*}
where obvious higher order terms are ignored,  $\textbf{w}^*$ is given by Eq.\eqref{step-1}. Substituting
$\mathcal{L}(\textbf{w}^*)$, $\mathcal{G}(\textbf{w}^*),$
$\mathcal{H}(\textbf{w}^*)$ into Eq.\eqref{step-2} and comparing the
coefficients of Eq.\eqref{exp}, a system of equations will be
obtained, and Eq.\eqref{sol} is its unique solution for this system.

To develop the gas-kinetic scheme with fifth-order temporal
accuracy, the time-dependent flux should by approximated by the
quadratic function, which is expressed as follows
\begin{align}\label{expansion}
F_{i+1/2}(W^n,t)=F_{i+1/2}^n+ \partial_t
F_{i+1/2}^nt+\frac{1}{2}\partial_{tt}F_{i+1/2}^nt^2.
\end{align}
In a time step, the coefficients $F_{i+1/2}^n$,
$\partial_tF_{j+1/2}^n$ and $\partial_{tt}F_{i+1/2}^n$ can be
determined as follows
\begin{align*}
F_{i+1/2}\Delta t+\frac{1}{2}\partial_t F_{i+1/2}\Delta
t^2+\frac{1}{6}\partial_{tt}
F_{i+1/2}\Delta t^3&=\mathbb{F}_{i+1/2}(W^n,\Delta t),\\
\frac{2}{3}F_{i+1/2}\Delta t+\frac{2}{9}\partial_t F_{i+1/2}\Delta
t^2+\frac{4}{81}\partial_{tt}
F_{i+1/2}\Delta t^3&=\mathbb{F}_{i+1/2}(W^n,2\Delta t/3),\\
\frac{1}{3}F_{i+1/2}\Delta t+\frac{1}{18}\partial_t F_{i+1/2}\Delta
t^2+\frac{1}{162}\partial_{tt} F_{i+1/2}\Delta
t^3&=\mathbb{F}_{i+1/2}(W^n,\Delta t/3),
\end{align*}
where
\begin{align*}
\mathbb{F}_{i+1/2}(W^n,\delta)=\int_{t_n}^{t_n+\delta}\int
uf(x_{i+1/2},t,u, v,\xi)dud\xi dt.
\end{align*}
The formulation for the gas distribution $f(x_{i+1/2},t,u, v, \xi)$
can be found in \cite{GKS-high1,GKS-high2,GKS-high3}. By solving the
linear system, all coefficients can be determined.

\end{document}